\documentclass[12pt]{article}

\usepackage{latexsym,amsfonts,amssymb}
\usepackage{amsmath}
\usepackage{amsthm}

\usepackage{psfrag}
\usepackage{epsfig}
\usepackage{pb-diagram}
\usepackage{amscd}

\usepackage{graphicx} 

\newtheorem{theorem}{Theorem}[section]

\newtheorem{claim}[theorem]{Claim}

\newtheorem{remark}[theorem]{Remark}

\newtheorem{definition}[theorem]{Definition}
\newtheorem{example}[theorem]{Example}
\newtheorem{exercise}[theorem]{Exercise}
\newcommand{\Z}{\mathbb{Z}}
\newcommand{\C}{\mathbb{C}}

\newcommand{\s}{\sigma}

\newcommand{\R}{\mathbb R}
\newcommand{\Q}{\mathbb Q}
\newcommand{\spa}{{\rm span}}

\newcommand{\dto}{\stackrel{d}{\to}}

\newcommand{\Hom}{{\rm Hom}}

\newcommand{\begeq}{\begin{eqnarray*}} 
\newcommand{\eneq}{\end{eqnarray*}} 
\newcommand{\mc}{\mathcal}
\newcommand{\lra}{\longrightarrow}

\newcommand{\oplusop}[1]{{\mathop{\oplus}\limits_{#1}}}
\newcommand{\oplusoop}[2]{{\mathop{\oplus}\limits_{#1}^{#2}}}
\newcommand{\Id}{\mathrm{Id}}

\title{Notes on Link Homology}
      \author{Marta Asaeda and Mikhail Khovanov}

\begin{document}

\date{April 8, 2008}
 
\maketitle
 
\begin{abstract}
This article consists of six lectures on the categorification of the Burau representation 
and on link homology groups which categorify the Jones and the 
HOMFLY-PT polynomial. The notes are based on the lecture course 
at the PCMI 2006 summer school in Park City, Utah. 
\end{abstract}

\newpage 

\tableofcontents

\newpage 

\section*{Introduction} 

These notes are based on lectures delivered by the second author at the PCMI 
summer school in Utah in the summer of 2006. The goal was to give an informal 
introduction at the graduate level to the ideas and constructions of combinatorial homology 
theories that categorify various quantum invariants of knots and links. We 
made an emphasis on the theories lifting the Jones polynomial and the HOMFLY-PT 
polynomial. Ideally, a link homology theory is a functor from the category of 
link cobordisms to some algebraic category, such as the category of abelian groups. 
A model example is worked out in Lectures 3-5. In Lecture 3 a categorification 
of the Jones polynomial to a bigraded link homology theory is sketched. In Lectures 4 and 5
we explain how to generalize this theory to tangles and tangle cobordisms. This 
generalization encodes a simple proof that the homology theory is functorial and 
extends to link cobordisms. Prior to that, in Lectures 1 and 2, we introduce a toy model of 
the story, a categorification of the Burau representation, which produces 
invariants of braids and braid cobordisms. In Lecture 6 we describe a triply-graded 
link homology theory categorifying the HOMFLY-PT polynomial. 

In the past few years link homology has become an extensively researched area, with 
a significant body of literature, which we won't try to fully survey in these short introductory 
lectures. Although neither knot Floer homology nor contact homology is discussed 
here, we refer the reader to the survey papers ~\cite{Ng2}, \cite{OS-M} on these 
topics and to~\cite{Tur} for another set of lecture notes on link homology. 

The final version of this work will appear in the AMS 
 lecture notes from the Graduate Summer School program on 
 Low Dimensional Topology  held in Park City, Utah, on June 25 - July 15, 2006. 

M.A. and M.K. would like to thank John Polking, the editor of the series,
Tom Mrowka and Peter Ozsv\'ath, the organizers of the summer school, for 
their encouragement and patience with the authors. 
The authors are indebted to Alexander Schumakovitch for providing knot homology 
tables for Lecture 3 and for enlightening discussions. M.A. was partially supported 
by the NSF grant DMS-0504199. M.K. was partially supported 
by the NSF grant DMS-0706924. This work was finished while the second author was 
at the Institute for Advanced Study, and he would like to thank the Institute and the 
NSF for supporting him during that time through the grant DMS-0635607. 


\newpage 

\section{A braid group action on a category of complexes} 

\subsection{Path rings}

\begin{definition}
An oriented graph $\Gamma$ (see a picture below) 
consists of finitely many vertices and oriented edges. For an 
edge $\xi$ let $s(\xi)$ and $t(\xi)$ be the source and the target vertices of $\xi$.
A path $\alpha$ is a concatenation of some edges $\xi_1, ... , \xi_k$,  so that 
$t(\xi_{i})=s(\xi_{i+1})$ for $i=1, ..., k-1$.  
We define $s(\alpha)$, $t(\alpha)$,  and the path length $|\alpha|$ to be 
$s(\xi_1)$, $t(\xi_k)$, and $k$, respectively. 
A path may be denoted by $(a_1 | a_2 | ... |a_k)$ where $a_i$'s are the vertices in 
the order that the path goes through, as long as there is only one such path.

 \psfig{figure=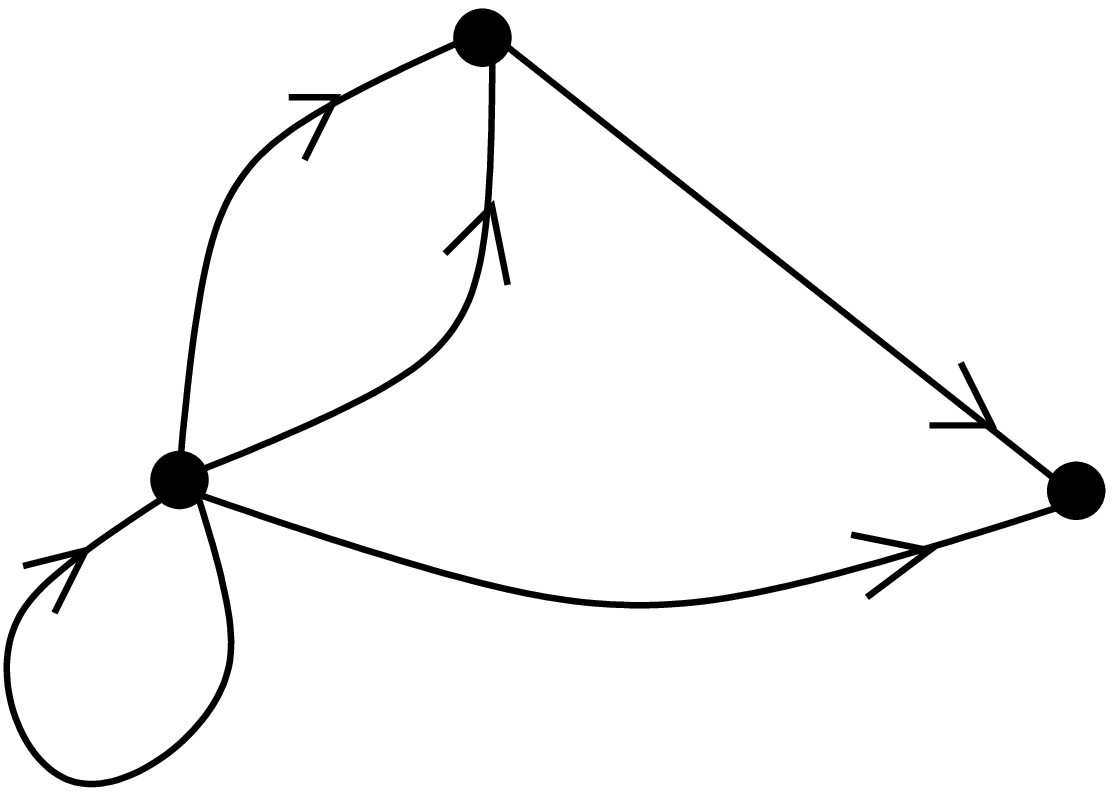,height=3cm}
 
The path ring $\Z \Gamma$ is a free abelian group with a basis given by all the 
paths in $\Gamma$, equipped with the following product: for paths 
$\alpha$ and $\beta$, $\alpha \beta$ is their concatenation if $t(\alpha)=s(\beta)$ 
and  zero otherwise.  We extend the product to $\Z \Gamma$ by linearity; this 
multiplication operation is associative. 
\end{definition}
\begin{example}
$\Gamma=\psfig{figure=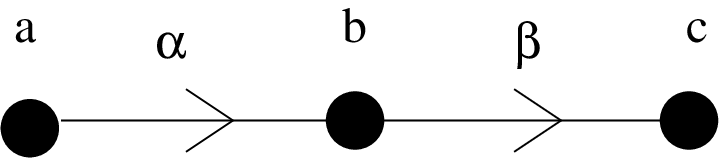,height=0.7cm}$. There are six paths in $\Gamma$: 
$\alpha$, $\beta$, $\alpha  \beta$,  and $(a)$, $(b)$, $(c)$, where the 
last three are the length zero paths consisting of a vertex. 
Note that $\beta  \alpha=0$, $(a)  (a) =(a)$, $(a)  \alpha =\alpha = \alpha (b)$, and so on. 
\end{example}
\begin{exercise}
Check that in the above example, $(a)+(b)+(c)$ is the unit element 
of the path ring $\Z\Gamma$. For any oriented graph $\Gamma$, the 
sum of the vertices (i.e. length zero paths) is the unit of $\Z \Gamma$. 
\end{exercise}
\begin{example}
$\Gamma=\raisebox{-0.3cm}{\psfig{figure=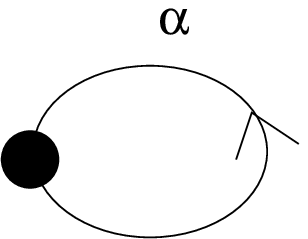,height=0.9cm}}$. Then 
$\Z \Gamma= \Z [\alpha]$ is the polynomial ring in one variable. 
\end{example}

Let ${\mathrm{mod}}$-$\Z\Gamma$ be the category of right $\Z\Gamma$-modules. 
Each object $M$ of ${\mathrm{mod}}$-$\Z\Gamma$  decomposes  into the direct 
sum  
$$M=\oplusop{a}M(a),$$ 
over the vertices $a$ of the graph, as  an abelian group. 
Multiplication by an edge $\xi$ with $s(\xi)=a$ and $t(\xi)=b$ is an abelian group 
homomorphism  $M(a)\to M(b)$. Vice versa, a right $\Z\Gamma$-module is determined 
by a collection of abelian groups, one for each vertex of the graph, and  
homomorphisms between these groups, one for each edge. 

\begin{exercise}
(1) Give a similar description of left $\Z\Gamma$-modules. \\
(1) Describe $\Z\Gamma$-module homomorphisms in this language.  
\end{exercise}

Converting $\Z$ into a field $k,$ we arrive at the notion of path algebra $k\Gamma.$ 
These algebras have homological dimension one (any submodule of a projective 
module is projective), just like rings of integers in number fields and rings of functions 
on smooth affine curves. Their representation theory is a spectacular story in progress; 
you can get a first taste of it from \cite{CB}. 

In this lecture we consider a very special quotient of a certain path ring. In general, 
if paths $\alpha_1, \dots, \alpha_m$ all have the same source vertex and the same 
target vertex, we can quotient the path ring by the relation 
$$ \lambda_1 \alpha_1 + \lambda_2 \alpha_2 + \dots + \lambda_m \alpha_m=0$$ 
for some $\lambda_1, \dots, \lambda_m \in \Z$. 


\subsection{Zigzag rings $A_n$}
For an $n>2$ consider the graph $\Gamma$ with vertices labelled from 1 to $n$ and oriented 
edges from $i$ to $i\pm 1$:
\psfrag{d1}{{}}
\psfrag{d2}{{}}
$$\raisebox{-0.3cm}{\psfig{figure=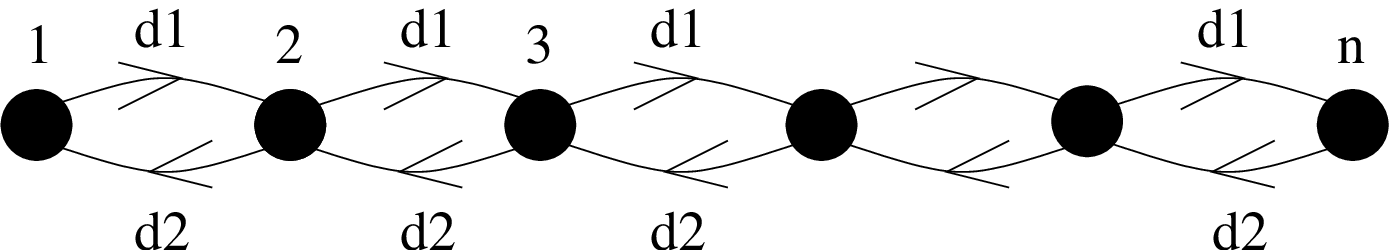,height=1cm}}$$
We define the ring $A_n$ as the quotient of $\Z\Gamma$ modulo the following 
relations
\begin{enumerate}
\item \psfig{figure=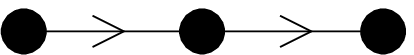,height=0.2cm}$=0$,  that is,  $(i|i+1|i+2)=0$;  
\item  \psfig{figure=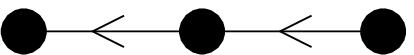,height=0.2cm}$=0$, that is,  $(i|i-1|i-2)=0$;  
\item  \raisebox{-0.2cm}{\psfig{figure=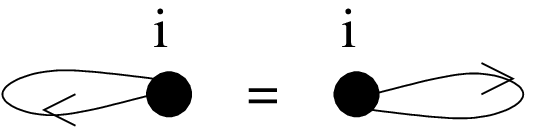,height=0.8cm}}:  
$(i|i-1|i)=(i|i+1|i)$. 
\end{enumerate}

If we label all edges pointed to the right (resp. left) by $\partial_1$ (resp. $\partial_2$), 
\psfrag{d1}{{\tiny $\partial_1$}}
\psfrag{d2}{{\tiny $\partial_2$}}
$$\raisebox{-0.3cm}{\psfig{figure=lec1.4.eps,height=1cm}}$$
the relations become
$$ \partial_1^2=0, \, \partial_2^2=0, \, \partial_1\partial_2=\partial_2\partial_1.$$
These are the relations for a bicomplex. In fact, the category of $A_n$-modules (either left or 
right) is equivalent to the category of mixed complexes of abelian groups, bounded above by 
$n$ \cite[2.5.13]{Lo}. The algebra $A_n$ is isomorphic to its opposite, hence its categories of 
left and right modules are equivalent. If we make $A_n$ graded, by assigning 
degree $1$ to all arrows going to the right and degree $0$ to all arrows going to the left, then 
the category of graded $A_n$-modules is isomorphic to the category of bicomplexes of abelian 
groups, with a suitable boundedness condition. 

We will use a different grading on $A_n$, by path length. Any path of length 
at least three is zero in $A_n$. Indeed, the first two relations imply that a non-trivial path 
should stay within some interval $[i-1,i]$. If its length is more than two, we can flip  a part of 
it, as illustrated below, to get zero.  
 $$ \psfig{figure=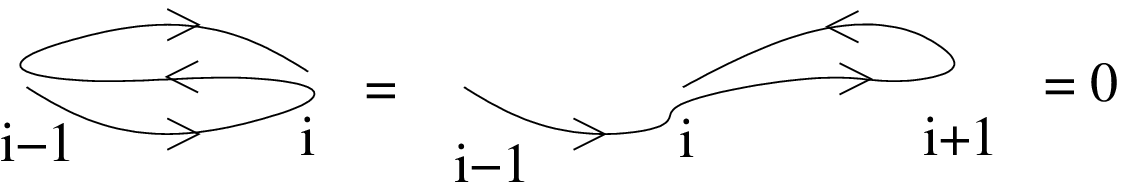,height=1cm} $$
It is easy to check that $A_n$ is a free abelian group with the basis of  
 \begin{enumerate}
 \item[$\bullet$]  length zero paths: $\left\{\raisebox{-.2cm}
{\psfig{figure=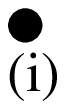,height=0.5cm}}   
\,  | i=1, ..., n \right\}$
 \item[$\bullet$] length one paths: $\left\{\raisebox{-.1cm}
{\psfig{figure=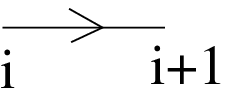,height=0.5cm}} 
, \raisebox{-.1cm}{\psfig{figure=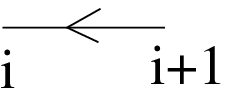,height=0.5cm}}   \,  | i=1, ..., n-1 \right\}$, 
 \item[$\bullet$]  length two paths:  $\left\{X_i := \raisebox{-0.2cm}
{\psfig{figure=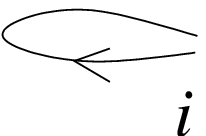,height=0.6cm}} \, {\rm or} \,  \raisebox{-.2cm}
{\psfig{figure=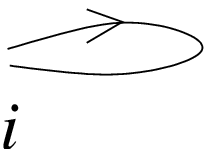,height=0.6cm}}  \,  | i=1, ..., n \right\}$
 \end{enumerate}
 

\subsection{A functor realization of the Temperley-Lieb algebra} 
In this section we make the Temperley-Lieb algebra act by functors on 
the category of $A_n$-modules.  
The Temperley-Lieb algebra $TL_{n+1}$ over the ground ring $R=\Z[q, q^{-1}]$ has 
generators $u_1, \dots , u_n$ and relations 
\begin{eqnarray*}
u_i^2&=&(q+q^{-1})u_i, \\ 
u_i u_{i\pm 1} u_i &=& u_i, \\
u_i u_j &=& u_j u_i, \ |i-j|>1. 
\end{eqnarray*}
The Temperley-Lieb algebra has a graphical interpretation, via the following assignment: 
$$u_i=\ \raisebox{-1cm}{\psfig{figure=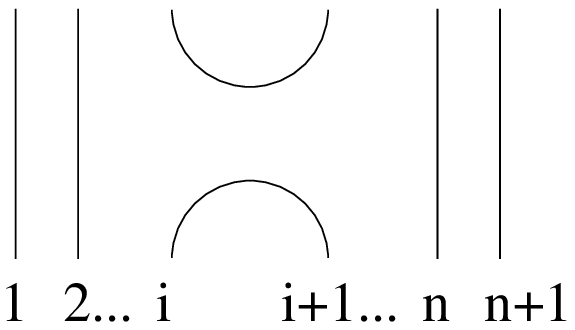,height=2.5cm}},$$
while the product of generators corresponds to concatenation  
$$ba=\ \raisebox{-1.2cm}{\psfig{figure=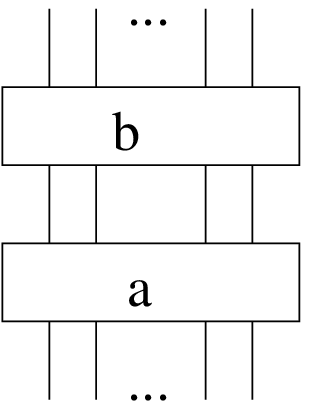,height=2.5cm}}$$
By setting the value of the closed loop to $q+q^{-1}$: 
$$\raisebox{-0.3cm}{\psfig{figure=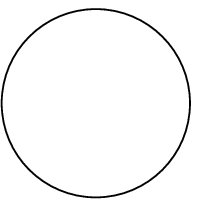,height=1cm}}=q+q^{-1},$$
and allowing arbitrary isotopies rel boundary, we obtain the relations in $TL_{n+1}$
(see \cite{KL} for more). 

$A_n$, as a left module over itself, decomposes into the direct sum 
$A_n= \oplusoop{i=1}{n} P_i$. 
Here 
$$P_i= A_n(i)={\rm span}_\Z \{ \raisebox{-0.2cm}
{\psfig{figure=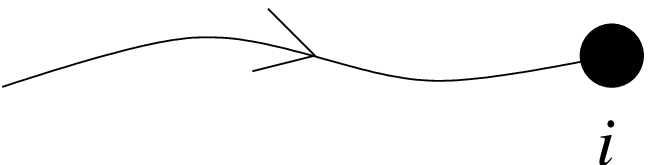,height=0.5cm}}\}.$$
 is a left projective $A_n$-module spanned over $\Z$ by paths that end in vertex $i.$ 
As an abelian group, $P_i$ is a free or rank $4$ with the basis 
$$\{(i), (i-1|i), (i+1|i), X_i\}$$ 
if $1<i<n$ and free of rank $3$ if $i=1,n$. 

Likewise, define the right projective $A_n$-module 
$$_i P :=(i)A_n={\rm span}_\Z \{ \raisebox{-0.2cm}
{\psfig{figure=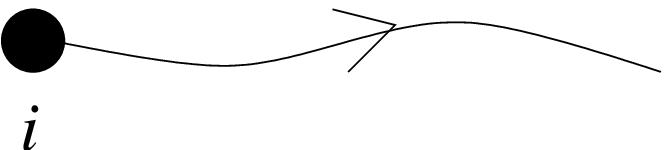,height=0.5cm}}\}.$$
\begin{exercise}
\label{ippj}
The following holds: 
$$_i P \otimes_{A_n} P_j=\left\{ 
\begin{array}{c}
0 \ \mbox{if} \ |i-j|>1, \\
\Z (i|j) \ \mbox{if} \ j=i\pm 1, \\
\Z(i) \oplus \Z X_i \ \mbox{if} \ i=j 
\end{array}
\right.
$$
(think of the LHS as spanned by paths that start in $i$ and end in $j$). 
\end{exercise}


Consider $A_n$-bimodules 
$$U_i := P_i \otimes_\Z {}_i P.$$
\begin{claim}
\label{Ui} There are bimodule isomorphisms
\begin{eqnarray*}
&& U_i \otimes_{A_n} U_i\cong U_i \oplus U_i, \\
&& U_i \otimes_{A_n} U_{i \pm 1} \otimes_{A_n} U_i\cong U_i, \\
&& U_i \otimes_{A_n}U_j =0, \ |i-j| >1.  \\
 \end{eqnarray*}
\end{claim}
\begin{proof}
We prove the first equality (the rest is equally easy to check). We use Exercise \ref{ippj}:
\begin{eqnarray*}
U_i \otimes_{A_n} U_i&\cong &P_i \otimes_\Z   ({}_i P \otimes_{A_n} P_i) 
\otimes_\Z  {}_i P \\
&\cong &(P_i \otimes_{\Z} \Z (i) \otimes_\Z  {}_i P) 
\oplus (P_i \otimes_{A_n} \Z X_i \otimes_\Z {}_i P)  \\
&\cong & (P_i \otimes_{\Z}  {}_i P) \, \oplus (P_i \otimes_\Z  {}_i P) \cong U_i \oplus U_i. \\
\end{eqnarray*}
\end{proof}
The equalities immediately remind us of the relations in $TL_{n+1}$ at $q=1$. The bimodule 
$U_i$ plays the role of the generator $u_i,$ tensor product of bimodules is analogous to 
the multiplication in $TL_{n+1}$, direct sum of bimodules lifts addition in the Temperley-Lieb 
algebra, etc. Due to the degenerate nature of our example, the tensor product 
$U_i \otimes_{A_n}U_j =0$ when  $|i-j| >1$ rather than just being isomorphic to 
the opposite tensor product. Thus, we get a bimodule realization of the quotient of 
$TL_{n+1}$ at $q=1$ by the relations $u_iu_j=0$ if  $|i-j| >1$ 
(we will construct a non-degenerate example in Lecture 4). 
The unit element $1$ of the Temperley-Lieb algebra corresponds to $A_n$, viewed  as 
a bimodule over itself. The canonical isomorphism $A_n \otimes_{A_n} M\cong M$, functorial 
in a bimodule $M$, lifts the identity $1m=m$ for $m\in TL_{n+1}$. 

We now bring $q$ into the play. Recall that $A_n$, $P_i,$ ${}_iP$ and $U_i$ are graded 
by path length. We work with graded modules and bimodules and denote by $\{m\}$ the 
grading shift up by $m$. Redefine $U_i$ by shifting its grading down by 1: 
$$ U_i = P_i \otimes_{\Z} {}_i P \{ -1\}.$$ 
For instance, the element $(i)\otimes (i)$ of $U_i$ now sits in degree $-1$.   
It is easy to see that there are isomorphisms of graded bimodules
\begin{eqnarray*}
&& U_i \otimes_{A_n} U_i\cong U_i\{1\} \oplus U_i\{-1\}, \\
&& U_i \otimes_{A_n} U_{i \pm 1} \otimes_{A_n} U_i\cong U_i, \\
&& U_i \otimes_{A_n}U_j =0, \ |i-j| >1.  \\
 \end{eqnarray*}
In this way, multiplication by $q$ becomes the grading shift $\{1\}$. 

To interpret the meaning of the minus sign in our bimodule realization 
of the Temperley-Lieb algebra we need to work with complexes of modules and 
bimodules, and take a small detour in the next subsection to review their basics. 

  
\subsection{The homotopy category of complexes}
Let $\mc{A}$  be an abelian category (for instance, the category of modules over 
some ring). Denote by $\mathrm{Kom}(\mc{A})$ the  category with objects--complexes 
of objects of $\mc{A}$ and morphisms--homomorphisms of complexes. A morphism $t$ from 
an object $M=\{\dots \to M^{i-1}\to M^i \to M^{i+1}\to \dots\}$ to 
$N = \{\dots \to  N^{i-1}\to N^i \to N^{i+1}\to \dots\}$ is a collection of morphisms 
$t_i: M^i \to N^i$ that make the following diagram commute 
$$
\begin{diagram}
\dgARROWLENGTH=1.5em
\node{M} \node{\cdots}\arrow{e,t}{d} \node{M^i}
\arrow{s,l}{t_i}\arrow{e,t}{d}\node{M^{i+1}}
\arrow{s,l}{t_{i+1}}\arrow{e,t}{d} \node{M^{i+2} }
\arrow{s,l}{t_{i+2}}\arrow{e,t}{d} \node{ \cdots} \\
\node{N} \node{ \cdots}\arrow{e,t}{d}\node{N^i} \arrow{e,t}{d}
\node{N^{i+1}}\arrow{e,t}{d}\node{N^{i+2}} \arrow{e,t}{d} \node{\cdots} 
\end{diagram}
$$
 $\mathrm{Kom}(\mc{A})$ is still an abelian category. Recall that a chain map 
$t$ is null-homotopic (we write $t \sim 0$) if there are maps $h_i:M^i \to N^{i-1}$
such that $t=dh+hd$ (in more detail, $t_i = d_N h_i + h_{i+1}d_M$). 
$$
\begin{diagram}
\dgARROWLENGTH=1.5em
\node{M} \node{\cdots}\arrow{e,t}{d} \node{M^i}\arrow{sw,l}{h}\arrow{s,l}{t}
\arrow{e,t}{d}\node{M^{i+1}} \arrow{sw,l}{h}\arrow{s,l}{t}\arrow{e,t}{d} 
\node{M^{i+2} }\arrow{sw,l}{h}\arrow{s,l}{t}\arrow{e,t}{d} \node{ \cdots} \\
\node{N} \node{ \cdots}\arrow{e,t}{d}\node{N^i} \arrow{e,t}{d}
\node{N^{i+1}}\arrow{e,t}{d}\node{N^{i+2}} \arrow{e,t}{d} \node{\cdots} 
\end{diagram}
$$
We define $\mathrm{Com}(\mc{A})$ as the quotient category of 
$\mathrm{Kom}(\mc{A})$ by the ideal of null-homotopic morphisms. 
\begin{exercise}
Check that null-homotopic morphisms constitute an ideal in $\mathrm{Kom}(\mc{A})$.
First you need to define the notion of an ideal in an abelian or an additive category.   
\end{exercise}
The quotient category has the same objects as $\mathrm{Kom}(\mc{A})$ but fewer 
morphisms: 
$$\Hom_{\mathrm{Com}(\mc{A})}(M,N)= \Hom_{\mathrm{Kom}(\mc{A})}(M,N)/\sim.$$ 
Two morphisms $f,g$ become equal in $\mathrm{Com}(\mc{A})$ if their difference is 
null-homotopic. 

Although we did not change objects when forming the quotient category, 
there are now more relations between them. 
\begin{exercise}
Check that for any nontrivial object $K$ of $\mc{A}$ complexes  
$$0 \to K \stackrel{{\rm id}}{\to} K\to 0 \ \mathrm{and}  \ 0 \to 0 \to 0 \to 0$$ 
are isomorphic in  $\mathrm{Com}(\mc{A})$ but not in  $\mathrm{Kom}(\mc{A})$.  
\end{exercise}
The category $\mathrm{Com}(\mc{A})$ is no longer abelian but triangulated 
(see \cite{GM}, \cite{Wei}) and comes with the following operations. \\
(1) Shift. For $M \in {\rm Ob}\ {\mathcal C}(A)$ define $M[j]$ to be the chain complex 
obtained from $M$ by shifting it $j$ steps to the left,   $M[j]^{i}=M^{i+j}$, and 
multiplying the differential by $(-1)^j$.  \\
(2) Cone of a morphism $f: M \to N$. The mapping cone of $f$ is the chain 
complex $C(f):=M[1]\oplus N$ with the differential $D:= -d_M + f + d_N$. 
Note that $C(f)^i= M^{i+1} \oplus N^i$. 
$$
\begin{diagram}
\dgARROWLENGTH=1.5em
\node{M^i}\arrow{e,t}{-d_M} \arrow{se,t}{f}\node{M^{i+1}} \arrow{se,t}{f}
\arrow{e,t}{-d_M}\node{M^{i+2} } \\
\node{N^{i-1}}\arrow{e,t}{d_N} \node{N^{i} } \arrow{e,t}{d_N} \node{N^{i+1}}
\end{diagram}
$$
From here on we specialize to categories of modules and bimodules. For a ring 
$A$ we denote by $\mc{C}(A)$ the category $\mathrm{Com}(A-\mathrm{mod})$ 
of complexes of $A$-modules up to chain homotopies. If $A$ is graded and we're working 
with graded modules, we'll use the same notation  $\mc{C}(A)$ for the category of 
complexes of graded $A$-modules up to chain homotopies. The differential of a complex 
of graded modules must preserve the grading. We can view a complex of graded 
$A$-modules as a bigraded $A$-module with a differential of bidegree $(1,0)$ which 
commutes with the action of $A$. 

Denote by $A^e=A\otimes A^{op}$ the tensor product of $A$ and its opposite ring. 
$A$-bimodules can also be described as left or right $A^e$-modules. We denote 
by $\mc{C}(A^e)$ the category of complexes of $A$-bimodules up to chain homotopy 
(and the category of complexes of graded $A$-bimodules, whenever necessary). 

Tensoring with a given $A$-bimodule is an endofunctor in the category of 
$A$-modules and  
an endofunctor in the category of $A$-bimodules. Likewise, tensoring with a complex 
of $A$-bimodules is an endofunctor in $\mc{C}(A)$ and $\mc{C}(A^e).$ 
The tensor product of $M,N\in \mc{C}(A^e)$ is the complex of bimodules given 
by placing the bimodule $M^i\otimes_A N^j$ into the $(i,j)$-node of the plane and 
then collapsing the grading onto the principal diagonal, so that the degree $k$ term of 
$M\otimes_A N$ is the direct sum 
$$ \oplusop{i\in \Z} M^i \otimes_{A} N^{k-i}, $$ 
with the differential combining those of $M$ and $N$: 
$$ d(m\otimes n) = d(m) \otimes n + (-1)^i m \otimes d(n), \hspace{0.2in} m\in M^i.$$  


\subsection{Braid group representation}
There exists a representation 
$\pi: Br_{n+1}\longrightarrow TL_{n+1}^{\ast}$ of the braid group on $n+1$-strands 
into the group of invertible elements in the Temperley-Lieb  algebra given on the 
standard generators of the braid group by $\pi(\sigma_i)=1-qu_i$. Graphically,  
$$\pi \left(  \raisebox{-0.5cm}{\psfig{figure=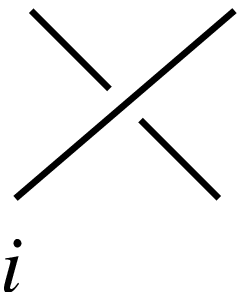,height=1.1cm}} \right)=
\raisebox{-0.5cm}{\psfig{figure=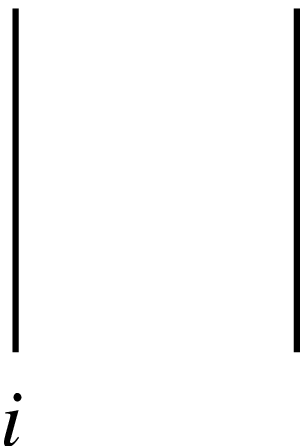,height=1.1cm}} -q \; 
\raisebox{-0.5cm}{\psfig{figure=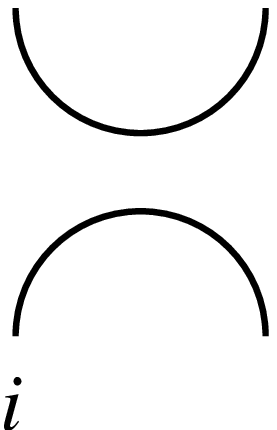,height=1.1cm}}. $$
What would the meaning of $1-qu_i$ be in our bimodule interpretation of the Temperley-Lieb 
algebra quotient? It should become the ``difference'' of graded bimodules $A_n$ and 
$U_i\{1\}=P_i\otimes_{\Z}{}_i P,$ which we interpret as the complex 
$$ 0 \lra P_i \otimes_{\Z}{}_i P\stackrel{\beta_i}{\lra} A_n \lra 0 $$ 
for the bimodule homomorphism $\beta_i$ which takes 
$x\otimes y \in P_i \otimes_{\Z}{}_i P$ to $xy\in A_n.$ This grading-preserving map
composes a path which ends in $i$ with a path which starts in $i$: 
$$  \raisebox{-0.3cm}{\psfig{figure=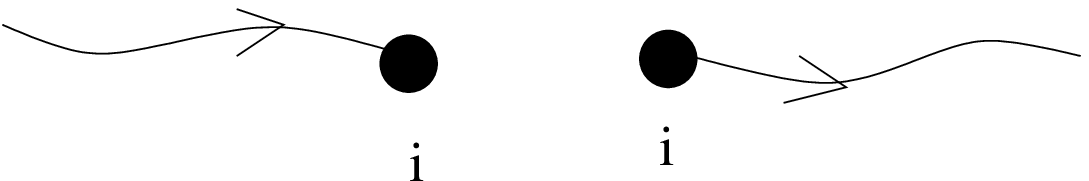,height=0.6cm}}\quad \mapsto \quad
   \raisebox{-0.4cm}{\psfig{figure=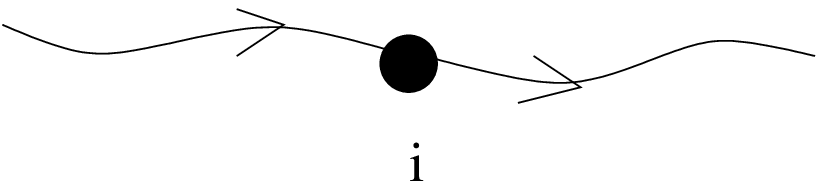,height=0.6cm}}$$
Thus, 
 $$ \sigma_i \raisebox{-0.5cm}{\psfig{figure=lec1.19.eps,height=1cm}} \mapsto 
(0 \to U_i\{1\} \stackrel{\beta_i}{\to} A_n \to 0).$$ 
We denote this complex of graded bimodules by $R_i$ and normalize it so that $A_n$ 
sits in cohomological degree $0$.  

The homomorphism $\pi:Br_{n+1}\to TL_{n+1}^{\ast}$ takes $\sigma_i^{-1}$ to 
$1- q^{-1}u_i.$ To interpret this difference we consider the complex $R_i'$ given by 
$$ 0 \lra A_n \stackrel{\gamma_i}{\lra} U_i\{-1\} \lra 0, $$ 
with the bimodule map $\gamma_i$ determined by the condition 
$$\gamma_i(1)=(i-1|i) \otimes (i|i-1) + (i+1|i) \otimes (i |i+1) + 
X_i \otimes (i) + (i) \otimes X_i$$
(for $1< i < n$; for $i=1,n$ omit one of the terms in the sum),
and $A_n$ placed again in cohomological degree $0$. 
\begin{theorem}
\label{reidem}
There are isomorphisms in $\mc{C}(A_n^e)$ of complexes of graded bimodules: 
\begin{eqnarray} 
R_i \otimes R_{i \pm 1} \otimes R_i  & \cong &  R_{i \pm 1} 
\otimes R_i  \otimes R_{i \pm 1}, \\
R_i \otimes R_j &  \cong & R_j \otimes R_i, \ |i-j|>1, \\
R_i \otimes R_i' & \cong & A_n \cong R_i'\otimes R_i 
\end{eqnarray} 
\end{theorem} 
The last relation tell us that $R_i$ and $R_i'$ are mutually inverse complexes of bimodules. 
The first two are the braid relations. The middle relation holds already in the abelian category 
of complexes of bimodules, before modding out by homotopies, but not the other two. 
The proof can be found in~\cite{KS} (also see the next lecture). 
This action was independently discovered by R.~Rouquier and A.~Zimmermann~\cite{RZ}; 
its algebraic geometry counterparts are studied in~\cite{ST}. 

The theorem implies that there is a braid group action on $\mc{C}(A_n)$ and on 
$\mc{C}(A_n^e)$ in which the generator $\sigma_i$ of the braid group $Br_{n+1}$ 
acts on a complex $M$ of graded $A_n$-modules (or bimodules) 
by tensoring it with $R_i$: 
$$\sigma_i(M) \ := \  R_i \otimes_{A_n} M, $$
while 
$$ \sigma_i^{-1}(M) \ :=  \ R_i'\otimes_{A_n} M .$$ 
More precisely, the theorem is about a group action in the \emph{weak} sense. A \emph{weak} 
action of a group $G$ on a category $\mc{C}$ assigns an invertible functor 
$F_g: \mc{C}\lra\mc{C}$ 
to each element of $G$ such that $F_{gh} \cong F_g F_h.$ For a weak action to be an action 
requires a specific choice of isomorphisms  $F_{gh} \cong F_g F_h$ for all $g,h\in G$ subject 
to the associativity relation that the diagram below is commutative for all $g,h,k\in G$:  
$$
\begin{CD}
F_{ghk} @>{\cong}>> F_{gh}F_k  \\
@V{\cong}VV     @V{\cong}VV     \\
F_g F_{hk} @>{\cong}>> F_g F_h F_k 
\end{CD}
$$ 

P.~Deligne~\cite{Del} gave a simple criterion for when a weak action of a braid group on 
a category can be upgraded to an action. His criterion holds in our 
case, and the weak action described above lifts to an actual action of $Br_{n+1}$ on 
$\mc{C}(A_n)$ and $\mc{C}(A_n^e)$. Furthermore, we have:
\begin{theorem}
The above action of the braid group $B_n$ on ${\mathcal C}(A_n)$ is faithful. 
\end{theorem}
We say that an action of the group $G$ on a category $\mc{C}$ is \emph{faithful} if 
the functors $F_g$ are not isomorphic for different $g\in G.$ 
See the next lecture for a sketch of a proof of the last theorem.


\section{More on braid group actions}


\subsection{Invertibility of $R_i$}
We begin the lecture with a sketch of isomorphisms: 
$$R_i \otimes R_i' \cong A_n  \cong R_i' \otimes R_i$$
from Theorem~\ref{reidem}. 
The double complex corresponding to $R_i \otimes R_i'$ has the form
$$
\begin{array}{ccccccc}
&&0&&0&& \\
&&\uparrow&&\uparrow&& \\
0& \to & U_i \otimes U_i & \to& A_n \otimes U_i\{-1\} & \to & 0 \\
&&\uparrow&&\uparrow&& \\
0& \to & U_i\{1\} \otimes A_n & \to &  A_n  \otimes A_n & \to & 0 \\
&&\uparrow&&\uparrow&& \\
&&0&&0&& 
\end{array}
$$
Noting that $U_i \otimes U_i \cong U_i\{1\} \oplus U_i\{-1\}$ from Claim \ref{Ui}, we obtain 
the total compex
$$R_i \otimes R_i' =( 0 \to U_i\{1\} \stackrel{d}{\to} A_n \oplus U_i\{1\} \oplus U_i\{-1\} 
\stackrel{d}{\to} U_i\{-1\} \to 0).$$
\begin{exercise} Write an explicit formula for the differential $d$ 
above and check that the complex decomposes into a direct sum 
$$
(0 \to U_i\{1\} \stackrel{1}{\to} U_i\{1\} \to 0) \oplus 
(0 \to A_n \to 0) 
\oplus (0 \to U_i\{-1\} \stackrel{1}{\to} U_i\{-1\} \to 0).
$$
\end{exercise}
The first and the last summands are null-homotopic, implying that 
$$R_i \otimes R_i' \cong (0 \to A_n \to 0) = A_n.$$


\subsection{Braid group action on complexes of projective modules $P_i$ and topology 
of plane curves}
In this section we explain how to prove that the braid group action on the homotopy 
category ${\mathcal C}(A_n)$ is faithful. For simplicity, we write $\sigma(P)$ for 
the object $F_{\sigma}(P)$ given by applying the functor $F_{\sigma}$ to 
$P\in{\mathcal C}(A_n)$. We will give a geometric presentation of  $\sigma(P_i)$
for all $\sigma$ in the braid group $Br_{n+1}$.  

First, let's look at a couple of easy examples. \\ 
\begin{example} \label{example2.1}
$$\sigma_i(P_i)=R_i \otimes P_i \cong (0 \to P_i \otimes_\Z {}_i P \otimes_{A_n} P_i  
\lra P_i \to 0)$$ 
The subcomplex 
$$ 0 \lra P_i \otimes (i) \otimes (i) \stackrel{1}{\lra} P_i \lra 0 $$ 
is contractible and the quotient complex is $0 \lra P_i\otimes X_i \otimes (i) \lra 0$
with the nontrivial term in cohomological degree $-1.$ The degree of $X_i$ is two and 
hence $P_i\otimes X_i \otimes (i) \cong P_i\{2\}$ as a graded module. 
thus $\sigma_i(P_i)\cong P_i [1]\{2\}$.  
\end{example} 
\begin{example}\label{example2.2} By induction on $m>0$ one can check that 
$$
\sigma_i^m(P_{i+1}) \cong 
(0 \to P_1\{2m-1\} \stackrel{X_i}{\lra} ... \stackrel{X_i}{\lra} P_1\{3\}\stackrel{X_i}{\lra} 
P_i\{1\}\stackrel{(i|i+1)}{\lra} P_{i+1} \to 0).
$$
\end{example}

We set up the following ingredients. Consider a disk with $n+1$ marked 
points aligned on a line as below. 
 $$\raisebox{-0.5cm}{\psfig{figure=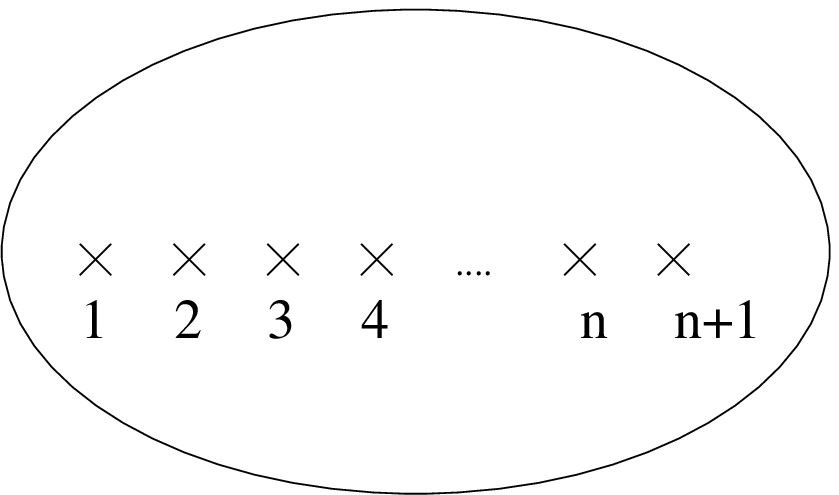,height=3cm}}$$
The braid group $B_{n+1}$ is isomorphic to the  the mapping class group of the disk that 
fixes the boundary and permutes the marked points. In particular, $B_{n+1}$ acts on isotopy classes 
of simple curves in the disk which have marked points as their endpoints and don't contain 
marked points in their interior. We assume that the generator 
$\s_i$ acts on the disk by permuting the vertices $i$ and $i+1$ counterclockwise. 
We fix a chain of curves $c_1, \dots , c_n$ as follows 
\psfrag{c1}{{ $c_1$}}
\psfrag{c2}{{ $c_2$}}
\psfrag{cn}{{ $c_n$}}
$$\psfig{figure=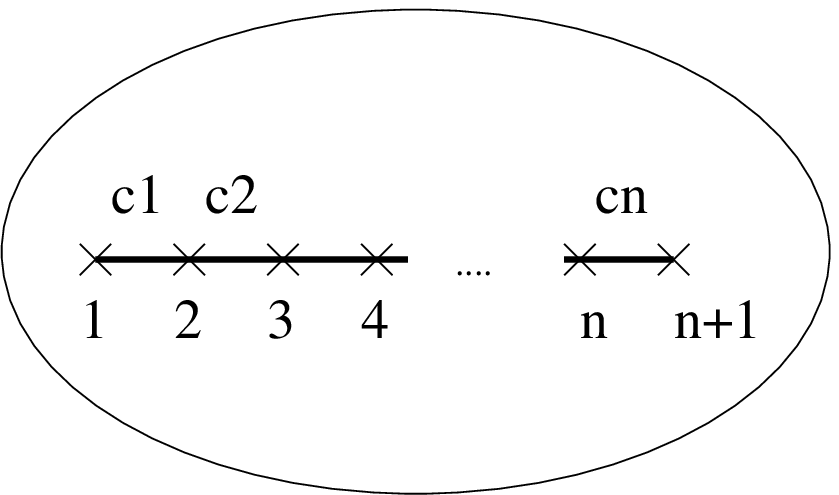,height=3cm}$$
The curve $c_i$ connects marked points $i$ and $i+1$. 
 The braid group action on the disc induces a braid group action on the isotopy 
classes of unoriented arcs that connects pairs of marked points. Any such isotopy class has 
the form $\sigma (c_i)$ for any $i$ and some braid $\sigma$. Curves $c_1, \dots , c_n$ 
represent some of these isotopy classes.  We would like to relate the braid group 
action on our category of complexes with the braid group action on the isotopy 
classes of curves. 

Consider vertical dotted lines $e_1, \dots , e_n$ orthogonal to $c_1, \dots, c_n$. 
\psfrag{c1}{{$e_1$}}
\psfrag{c2}{{$e_2$}}
\psfrag{cn}{{$e_n$}}
$$\psfig{figure=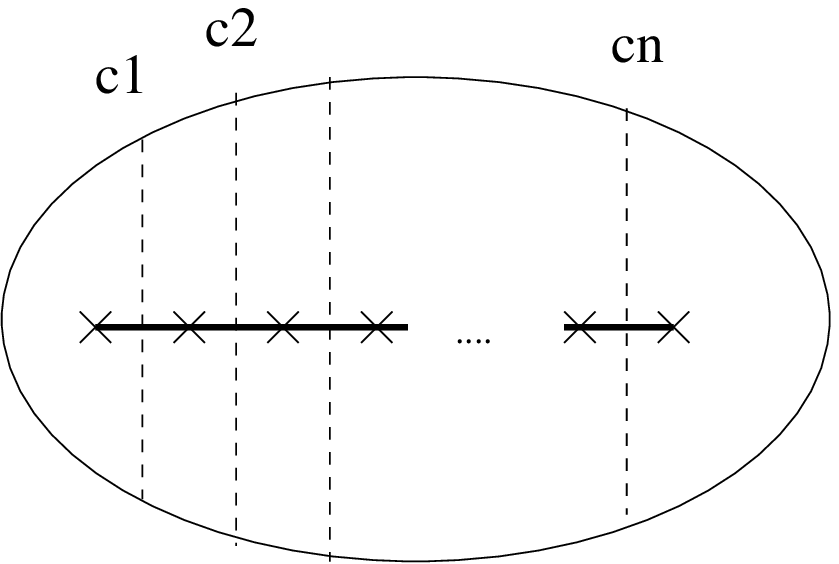,height=3cm}$$
Given an isotopy class $c$ of an arc in the disc with marked endpoints, we 
can choose a representative $c'$ in the minimal position relative to the system of 
intervals $e_1, \dots, e_n$, in the sense that the number of intersection points of 
$c$ with each of $e_i$ is the minimal possible among curves in the isotopy class $c$. 
Such representative is unique in the appropriate sense and can be obtained 
from any generic diagram in $c$ by a sequence of simplifications 
\psfrag{Arrow}{{$\Longrightarrow$}}
$$\psfig{figure=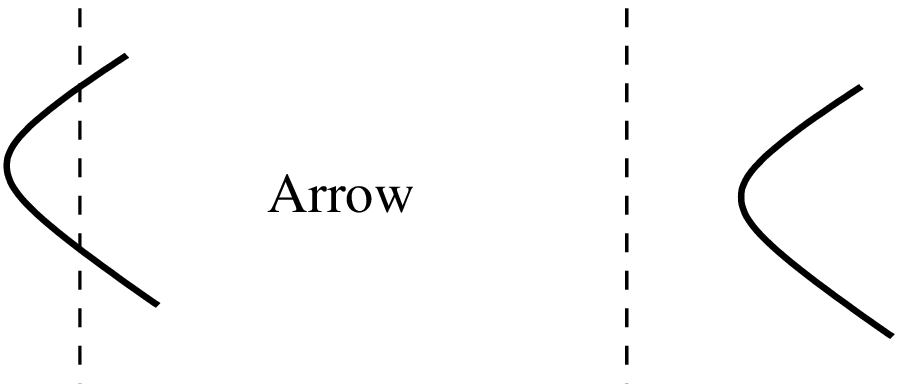,height=2cm}$$
Here's an example of the minimal representative. 
\psfrag{fa1}{{$e_1$}}
\psfrag{fa2}{{$e_2$}}
\psfrag{fa3}{{$e_3$}}
\psfrag{fa4}{{$e_4$}}
$$\psfig{figure=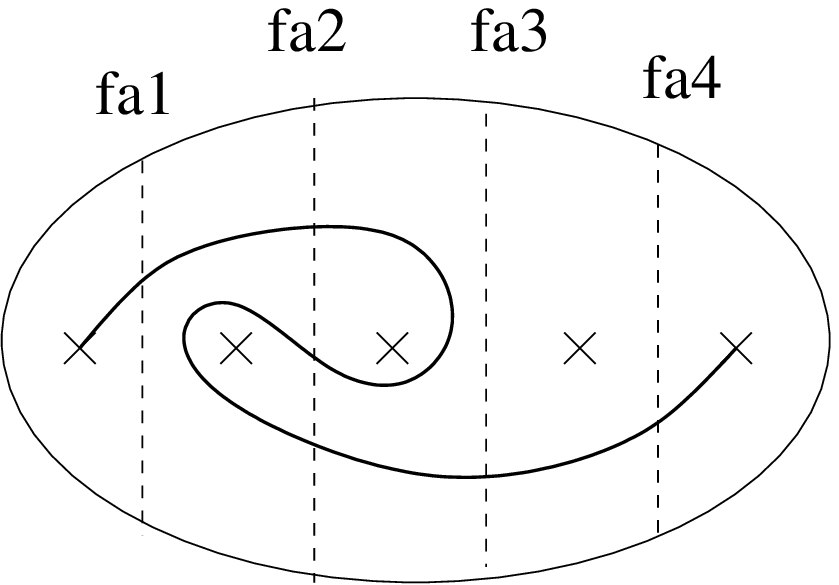,height=3cm}$$

To a isotopy class $c$ of arc we assign a complex $P(c)$ 
of projective $A_n$-modules as follows. Let $c'$ be the minimal 
representative of $c$. The vertical lines cut $c'$ into segments.
Discard two segments containing the endpoints of $c'$ and orient   
each remaining segment bounded by vertical lines $e_i$ and $e_{i+1}$ 
clockwise around the marked point $i+1$, the only marked point between 
these vertical lines.
\psfrag{fa1}{{$e_1$}}
\psfrag{fa2}{{$e_2$}}
\psfrag{fa3}{{$e_3$}}
\psfrag{fa4}{{$e_4$}}
$$\psfig{figure=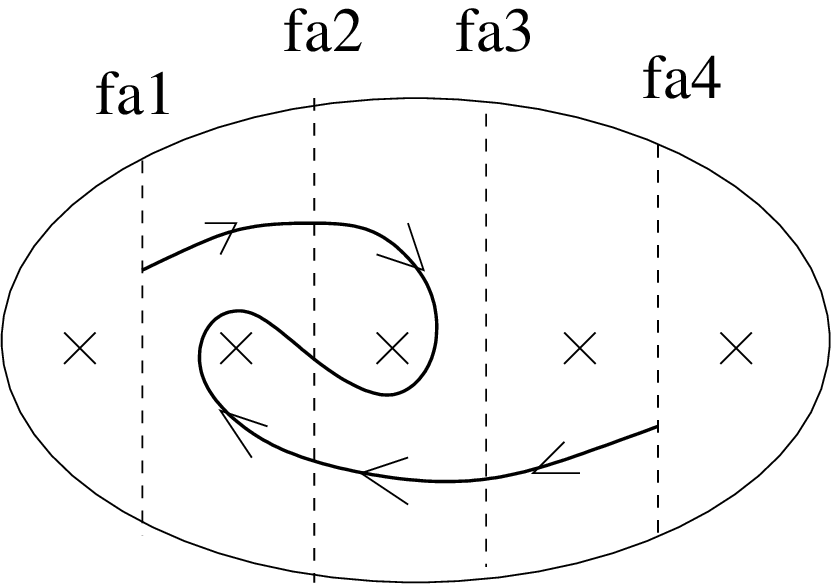,height=3cm}$$
To each intersection point of the curve with the vertical line $e_i$ assign number $i$. 
Now, pull the curve with these additional decorations out of the disk  
and draw  it on the plane  so that all orientations look to the right. 
$$\psfig{figure=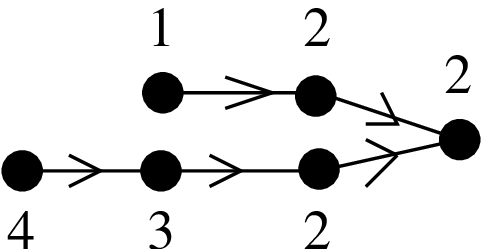,height=2cm}$$
We next put  $P_i$ at each vertex labeled $i$ and take the direct sum vertically. 
Define the differential as the sum of contributions from each arrow. To an arrow 
from $i$ to ${i\pm 1}$ assign the module homomorphism 
$P_i \lra P_{i\pm 1}$ which takes $a$ to $a(i|i\pm 1)$. To an arrow 
from $i$ to $i$ assign the homomorphism  $P_i \lra P_i$ of multiplication by 
$X_i$. We obtain a chain complex; for the example above it has the form 
$$ 
\begin{diagram}
\dgARROWLENGTH=0.7em
\node[4]{} \node{P_1}\arrow[2]{e,t}{ (1|2)} \node{} \node{P_2}
\arrow[2]{se,l}{X_2} \node{} \node{} \\
 \node[4]{} \node{}\node{} \node{}\node{}\node{} \\
\node[4]{} \node{\oplus}  \node{} \node{\oplus}  \node{} \node{P_2} \\
 \node[4]{} \node{}\node{} \node{}\node{}\node{} \\
\node{} \node{} 
\node{P_4}\arrow[2]{e,t}{(4|3)} \node{} \node{P_3}
\arrow[2]{e,t}{ (3|2)} \node{} \node{P_2}\arrow[2]{ne,r}{X_2}\node[2]{}
\end{diagram}
$$
which we can also write as 
$$ 0 \to P_4 \lra P_1\oplus P_3 \lra P_2 \oplus P_2 \lra P_2 \to 0. $$
In this way to an isotopy class $c$ of arcs we assign a complex $P(c)$ of 
projective $A_n$-modules. We did not specify the overall grading shift for $P(c)$; 
the reader can find this and other information in \cite{KS}. It is also possible to 
keep track of the internal grading and view $P(c)$ as a complex of  
graded $A_n$-modules. For $c=c_i$ the complex $P(c_i)=(0 \lra P_i \lra 0)$. 

\begin{theorem} \label{thm-braid} 
For any braid $\sigma$ and any number $i$ between $1$ and $n$ 
the complex $P(\sigma c_i)$ is homotopy equivalent to $\sigma P_i$.  
\end{theorem}

This theorem \cite{KS} tells us how a braid $\sigma$ acts on projective modules $P_i$, 
and that the action can be read off the braid group action on isotopy classes 
of arcs. 

\begin{exercise} Rethink Example~\ref{example2.2} via this theorem. \end{exercise} 

To prove  that the braid group action on ${\mathcal C}(A_n)$ 
is faithful it suffices to check that for any nontrivial braid $\sigma$ we 
can find some $i$ so that $\sigma P_i$ is not isomorphic to $P_i$ in the 
homotopy category. Our description of $\sigma P_i$ implies that if it is isomophic to $P_i$ 
then $\sigma c_i = c_i.$ If $\sigma c_i = c_i$ for all $i$ then $\sigma$ is central 
and is a multiple of the full twist. But it is easy to compute that the full twist takes $P_i$ 
to $P_i[j]$ for some $j\not= 0$ (compare with Example~\ref{example2.1}). 
\begin{exercise} Find this $j$. \end{exercise}  
The faithfullness of the action follows, modulo Theorem~\ref{thm-braid}, not proved 
in these notes.


\subsection{Reduced Burau representation} 
A braid takes a complex of (graded) projective $A_n$-modules to 
a complex of (graded) projective $A_n$-modules. One can check that any 
finitely-generated projective $A_n$-module is isomorphic to a direct sum of $P_i$'s, 
and the multiplicity of $P_i$ in this decomposition is an invariant of the projective 
module. Likewise, any finitely-generated projective graded $A_n$-module is 
isomorphic to a direct sum of $P_i\{j\},$ and the multiplicity of $P_i\{j\}$ is 
an invariant of the module. For the rest of this section all modules are 
assumed to be left, graded and finitely generated.  We introduce a
formal symbol $[P]$ of each projective module $P$. Let $K_0(A_n)$ 
be the $\Z[q,q^{-1}]$-module generated by these symbols 
subject to relations 
$$ [P\oplus Q] = [P] + [Q], \ \  [P \{ j \} ] = q^j [P]. $$
The direct sum decomposition property mentioned above implies that 
$K_0(A_n)$ is a  free $\Z[q,q^{-1}]$-module 
generated by symbols of indecomposables $[P_1], \dots, [P_n]$. 
The reader familiar with K-theory will recognize $K_0(A_n)$ as the group 
$K_0$ of the category of graded finitely-generated $A_n$-modules (see \cite{Ros} 
for an excellent introduction to algebraic K-theory). 

Given a bounded complex $P$ of  projective $A_n$-modules 
$$ 0 \to \dots \to P^i \to P^{i+1}\to \dots \to 0, $$ 
we define its Euler characteristic as 
$$ \chi(P) = \sum_i  (-1)^i [P_i] \in K_0(A_n).$$ 
If two complexes are chain homotopy equivalent, they have 
equal Euler characteristic; shifting the complex by $1$ adds a minus 
sign to the Euler characteristic,  $\chi(P[1]) = - \chi(P).$ 

The action of the braid group on the category of complexes of projective
$A_n$-modules (a subcategory of $\mc{C}(A_n)$) 
descends to a $\Z[q,q^{-1}]$-linear action of the braid group on $K_0(A_n).$ 
To determine this action, we write how $\sigma_i$ acts on $P_j$: 
\begin{eqnarray*} 
 \sigma_i(P_i) & \cong & P_i[1]\{2\} \\
 \sigma_i(P_{i\pm 1}) & \cong &  0 \to P_i\{1\} \to P_{i\pm 1} \to 0 \\
 \sigma_i(P_j) & \cong & P_j \ \mathrm{if} \ |i-j|>1, 
\end{eqnarray*} 
and pass to the Euler characteristic ($\sigma[P]=[\sigma(P)]$ by definition): 
\begin{eqnarray*} 
 \sigma_i[P_i]& = & -q^2 [P_i] \\
 \sigma_i[P_{i\pm 1}] & = & [P_{i\pm 1}] - q  [P_i] \\
 \sigma_i[P_j] & = & [P_j] \ \mathrm{if} \ |i-j|>1. 
\end{eqnarray*} 
The resulting action is isomorphic to the reduced Burau representation 
of the braid group. Hence, we can view the braid group action on the category of complexes 
of projective $A_n$-modules as a categorification of the Burau representation. 

To elements of the braid group we assign functors and it turns out that this 
assignment can be extended to braid cobordisms. A braid cobordism is a
surface in $\R^4$ that goes from one braid to the other. The condition that 
braids have no critical points when projected onto the $z$-axis extends to the 
condition that a braid cobordism is a simple branch covering when projected 
onto the $(z,w)$-plane in $\R^4.$ To a braid cobordism between braids 
$\sigma$ and $\tau$ it is possible to assign a natural tranformation between 
functors $F_{\sigma}$ and $F_{\tau}$ (modulo the issue of the overall sign) 
in a consistent way which respects compositions of braid cobordisms, see \cite{KT}. 
This example is the simplest way to get an algebraic invariant of a toy sector 
of four-dimensional topology (braid cobordisms) from homological algebra (of 
complexes of $A_n$-modules).


\section{ A Categorification of the Jones polynomial}

\subsection{The Jones polynomial}
The Jones polynomial \cite{Jones} 
is an isotopy invariant of oriented links in $\R^3$ that takes values in
 $\Z[q,q^{-1}]$ and satisfies the following skein relation for link diagrams:
\begin{equation}\label{skein} 
q^2 J \left( \raisebox{-0.3cm}{\psfig{figure=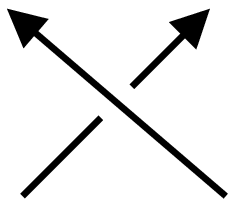,height=0.8cm}}\right) -
q^{-2} J \left( \raisebox{-0.3cm}{\psfig{figure=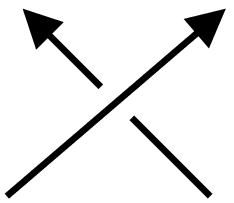,height=0.8cm}}\right) 
= (q-q^{-1})J\left( \raisebox{-0.3cm}{\psfig{figure=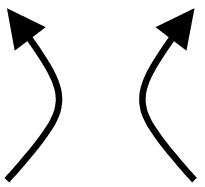,height=0.8cm}}\right). 
\end{equation} 
This equation implies 
\begin{equation} 
J\left(L \sqcup  \raisebox{-0.3cm}{\psfig{figure=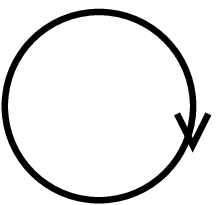,height=0.8cm}}
\right) =  (q+q^{-1}) J(L),
\end{equation}
that is, the Jones polynomial of the disjoint union of a link $L$ and the unknot is the 
Jones polynomial of $L$ times $q+q^{-1}$. Adding the unknot to a link multiplies the 
Jones polynomial by $q+q^{-1}.$ Hence, it is convenient to normalize
the invariant to take this value on the unknot:    
\begin{equation} 
J\left(  \raisebox{-0.3cm}{\psfig{figure=lec3.4.eps,height=0.8cm}}
\right) =  q+q^{-1}.
\end{equation}
We also extend the invariant to the empty link, by setting $J(\emptyset)=1$. 

Inductive simplification via the  skein relation (\ref{skein}) implies that 
the Jones polynomial is uniquely determined by this relation and its 
value on the unknot. A simple way to prove existence was found by 
Louis Kauffman \cite{Kau}. Define the Kauffman bracket 
polynomial $\langle D \rangle$ of an \emph{unoriented} link projection $D$ by expanding 
every crossing 
\begin{equation} \label{cr-res}
\left< \raisebox{-0.3cm}{\psfig{figure=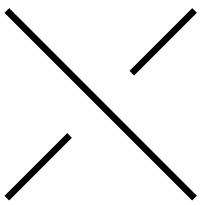,height=0.8cm}}  \right>=
\left< \raisebox{-0.3cm}{\psfig{figure=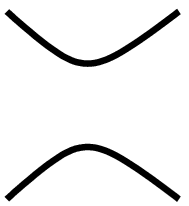,height=0.8cm}}  \right> -q^{-1} 
\left< \raisebox{-0.3cm}{\psfig{figure=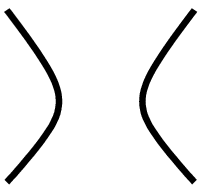,height=0.8cm}}  \right> 
\end{equation} 
and requiring that $\langle D \rangle= (q+q^{-1})^k$ if $D$ is a crossingless 
diagram with $k$ circles (the skein relation above differs slightly from Kauffman's 
original one, which has more symmetry). The Kauffman bracket 
of a planar diagram is invariant under the Reidemeister moves up to rescaling by 
plus or minus a power of $q$. To get rid of these scaling factors, 
define the Kauffman bracket of an \emph{oriented} planar diagram by 
\begin{equation}\label{orkau} 
K(D):=(-1)^{x(D)}q^{2 x(D)- y(D)} \langle D\rangle ,
\end{equation}
where $x(D)$, respectively $y(D)$, is the number of negative
\raisebox{-0.2cm}{\psfig{figure=lec3.1.eps,height=0.5cm}} , respectively 
positive \raisebox{-0.2cm}{\psfig{figure=lec3.2.eps,height=0.5cm}}, crossings of $D$, 
and $D$ is viewed as an unoriented diagram in the rightmost term of (\ref{orkau}).   

\begin{exercise}
Show that $K(D)=J(L)$ for any planar diagram $D$ of $L$. 
\end{exercise}

Thus, the Kauffman bracket is a link invariant, equal to the Jones polynomial. 
The two modifications
of a crossing that appear on the right hand side of (\ref{cr-res}) will be called 
\emph{resolutions}. 
We call modification $ \raisebox{-0.2cm}{\psfig{figure=lec3.5.eps,height=0.5cm}}  
\to  
 \raisebox{-0.2cm}{\psfig{figure=lec3.7.eps,height=0.5cm}}$ the {\it 0-resolution},  
and modification $ \raisebox{-0.2cm}{\psfig{figure=lec3.5.eps,height=0.5cm}}  \to  
 \raisebox{-0.2cm}{\psfig{figure=lec3.6.eps,height=0.5cm}}$ the {\it 1-resolution}.  
 

\subsection{Categorification and a bigraded link homology theory}

In one of its manifestations, \emph{categorification}, a term introduced 
by Louis Crane and Igor Frenkel \cite{CF}, 
lifts natural numbers to vector spaces or free abelian group. 
Going in the opposite direction (decategorifying), to a finite-dimensional vector 
space $V$ we assign its dimension $\dim(V)$ 
and to a finitely-generated free abelian group $V$ its rank  $\mathrm{rk}(V)$. 
Operations on vector spaces or free abelian groups mirror those on natural 
numbers. Direct sum of vector spaces corresponds to the sum of numbers, tensor 
product to multiplication: 
$$\dim(V\oplus W) =\dim(V) + \dim(W), \hspace{0.2in} 
  \dim(V\otimes W) = \dim(V) \dim(W).$$ 
Thus, we have an informal correspondence 
\begin{eqnarray*}
n &\leftrightarrow & \Z^n \  \mbox{or}\  k^n , 
\  \mathrm{for \, a \, field} \  k,   \\
n+m & \leftrightarrow & V \oplus W, \ \mbox{where} \  {\rm rk}(V)=n, 
{\rm rk}(W)=m, \\
nm &\leftrightarrow & V \otimes W.
\end{eqnarray*}
Lifting negative numbers and differences $n-m$ requires stepping beyond 
the category of vector spaces and considering the category of complexes of 
vector spaces or free abelian groups. 
The analogue of the dimension of a vector space is the Euler characteristic 
of a complex. In the simplest instance, if positive 
integers $n$ and $m$ have become vector spaces $V$ and $W$ of 
dimension $n$ and $m$, then the difference $n-m$ is the Euler 
characteristic of the complex $0 \to V \stackrel{d}{\lra} W \to 0$
for any linear map $d$, with $V$ sitting in even cohomological degree. 
More generally, if we have already lifted integers $n$ and $m$ to complexes $V$ and 
$W$, then $n-m$ can be interpreted as the Euler characteristic of the 
complex $\mathrm{Cone}(f)$ for some map $f: W\lra V$ of complexes
(alternatively, as the Euler characteristic of $\mathrm{Cone}(g)[\pm 1]$ 
for a map $g: V \lra W$). 

The standard example of categorification is passing from the Euler characteristic 
$\chi(X)$ of a topological space $X$ to its homology groups 
$\mathrm{H}_{\ast}(X, \Z)$. 
We recover the Euler characteristic by taking the alternating sum of ranks 
\begin{equation}\label{euler-char} 
\chi(X) = \sum_{i} (-1)^i \mathrm{rk}\ \mathrm{H}_i(X, \Z).
\end{equation}  
Homology can be built from the Euler characteristic by a lifting as above, 
starting with a CW-decomposition $X'$ of $X,$ taking the formula for $\chi(X)$ as 
the alternating sum of the number of $i$-dimensional cells,  lifting each 
$\pm 1$ term in the sum to the complex $0 \lra \Z \lra 0$ in the corresponding 
degree, and forming the complex $C(X')$ with the homology $\mathrm{H}_{\ast}(X, \Z)$.
Notice the multitude of benefits that the homology of $X$ provides compared to 
the Euler characteristic of $X$: 
\begin{itemize} 
\item The invariant is not just an integer but a graded abelian group, 
encoding more information about $X$. 
\item Homology extends to functor from the category of topological spaces and 
continuous maps modulo homotopies to the category of graded abelian groups and 
grading-preserving homomorphisms. Thus, it provides information about continuous 
maps as well, associating to $f:X\lra Y$ the homomorphism 
$$\mathrm{H}_{\ast}(f): \mathrm{H}_{\ast}(X, \Z)\lra \mathrm{H}_{\ast}(Y, \Z).$$
The Euler characteristic does not give any information about continuous maps. 
\item Homology groups (singular homology) are defined for any topological 
space. The Euler characteristic, in its naive version, is only defined for topological 
spaces admitting finite CW-decomposition. Once homology becomes available 
the Euler characteristic can be defined for a wider range of spaces via 
equation (\ref{euler-char}), assuming that $\mathrm{H}_{\ast}(X, \Z)$ 
has finite rank. Still, for many spaces ($\mathbb{CP}^{\infty}$ or a discrete 
infinite space are the simplest examples) the formula (\ref{euler-char}) does not 
help due to  $\mathrm{H}_{\ast}(X, \Z)$  having infinite rank. For such $X$ the 
Euler characteristic cannot be defined but the homology still makes sense. 
\item Relatives of homology groups such as the cohomology groups of $X$ and  
the $K$-theory of $X$ provide even more information via the multiplication in 
cohomology and in $K$-theory, cohomological operations, etc. We get 
a highly sophisticated theory called algebraic topology.  
\end{itemize} 

Now take a math book, find there a structure which is manifestly 
integral, with all the structure constants, coefficients, etc. being integers, 
and try to categorify it. This means consistently lifting integers to 
vector spaces or complexes of vector spaces. A book on combinatorics is a 
good place to start (you're likely to have less luck with a book on analysis 
since integral structures are not common there). The end result of 
your categorification efforts must be a richer and more beautiful structure 
then the one you started with, living one level above the original. Normally, 
in many cases your attempts will fall apart or the result will look artificial or 
shallow, but eventually you might be onto something. 
 
Today we will look at a successful case of categorification--a categorificaton 
of the Jones polynomial. The Jones polynomial $J(L)$ of  a link $L$ takes values 
in $\Z[q,q^{-1}]$, so its coefficients are integers: 
 $$J(L) =\sum_{j} a_j(L)  q^j, \  \ \  a_j(L)\in \Z.$$
We will realize each coefficient as the Euler characteristic of a $\Z$-graded 
link homology theory. Taking all cofficients together, we'll get bigraded 
homology groups associated to a link 
  $$H(L)=\oplusop{i,j} H^{i,j}(L)$$
  so that 
  $$a_j(L) = \sum_i (-1)^i {\rm rk} H^{i,j}(L), \ \  j\in \Z, \ \  
   J(L) = \sum_{i,j}(-1)^i q^j {\rm rk} H^{i,j}(L).$$ 
The homology will be constructed by lifting the Kauffman bracket formula 
for the Jones polynomial to complexes. To a diagram $D$ we will assign 
a complex $C(D)$ of graded free abelian groups 
$$ \cdots \lra C^{i-1}(D) \stackrel{d}{\lra} C^i(D) \stackrel{d}{\lra} 
C^{i+1}(D)\stackrel{d}{\lra} \cdots $$ 
with a grading-preserving differential;  $C^i(D) = \oplusop{j} C^{i,j}(D)$. 
 For a given degree $j$ the complex restricts to a complex of free abelian groups 
$$ \cdots \lra C^{i-1,j}(D) \stackrel{d}{\lra} C^{i,j}(D) 
\stackrel{d}{\lra} C^{i+1,j}(D)\stackrel{d}{\lra} \cdots , $$ 
with $a_j(L)$ being its Euler characteristic. 

It is convenient to imagine groups $C^{i,j}(D)$ as sitting in the $(i,j)$-square on 
the plane and differential going one step to the right. The grading shift $\{1\}$ 
moves everything one step up, and the shift $[1]$ moves the diagram one step to 
the left. We refer to the $i$-degree as horizontal/cohomological degree and 
the $j$-degree as vertical/internal degree and also as $q$-degree. 

First, we directly categorify the inductive formula for $\langle D \rangle$ for 
an unoriented diagram $D$ and turn $\langle D \rangle$ into the Euler characteristic 
of a complex $\overline{C} (D).$ Next, orienting $D$, we'll define 
\begin{equation}\label{add-shift} 
 C(D) := \overline{C} (D) [x(D)]\{2x(D)-y(D)\}, 
 \end{equation} 
mirroring the formula (\ref{orkau}). 

We start with the simplest diagrams. For the empty diagram $\langle \emptyset\rangle =1$,  
and we define $\overline{C} (\emptyset)=\Z$ in bidegree $(0,0).$ 
For a single circle diagram
$\langle\raisebox{-0.1cm}{\psfig{figure=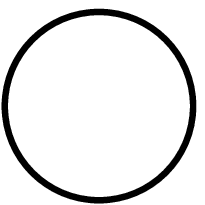,height=0.5cm}}\rangle=q+q^{-1}.$
Let $A=\Z{\bf 1}\oplus \Z X$ be a free graded abelian group 
with the basis $\{{\bf 1},X\}$ such that $\deg({\bf 1})=-1$ and $\deg(X)=1$ (the 
reason for notation ${\bf 1}$ will soon become clear). 
The graded rank of $A$ is $q+q^{-1}$ and we declare that 
$$\overline{C} (\raisebox{-0.1cm}{\psfig{figure=lec3.8.eps,height=0.5cm}})=A,$$ 
viewed as a complex of graded abelian groups $0 \lra A \lra 0$ sitting in cohomological 
degree $0$ (necessarily with the trivial differential).
In general, consider an arbitrary plane diagram $D$ without crossings.  Such 
diagram $D$ consists of $k$ disjoint circles embedded into the plane, possibly 
in a nested way. To such $D$ we assign the complex $\overline{C} (D):=A^{\otimes k}$ 
with the trivial differential and $A^{\otimes k}$ sitting in the cohomological 
degree $0$. The graded rank of $A^{\otimes k}$ is $(q+q^{-1})^k = \langle D\rangle $. 

Next, we need to tackle diagrams with crossings and interpret the relation 
(\ref{cr-res}) in our framework. Assuming that complexes for the two diagrams 
on the right hand side of (\ref{cr-res}) have already been defined, we could look for a 
homomorphism of complexes 
\begin{equation} \label{cr-hom}
f: \overline{C} ( \raisebox{-0.2cm}{\psfig{figure=lec3.6.eps,height=0.6cm}}) 
\lra \overline{C} ( \raisebox{-0.2cm}{\psfig{figure=lec3.7.eps,height=0.6cm}})  
\end{equation} 
and define $\overline{C} (\raisebox{-0.2cm}{\psfig{figure=lec3.5.eps,height=0.6cm}})$ 
as the cone of $f$ shifted one degree to the right (compare with~\cite{}). 
The relation (\ref{cr-res}) would then hold 
for the Euler characteristics of these three complexes. Here's how it works in the 
simplest cases. 

 \noindent
{\bf Example 1:} a kinked diagram $D= \raisebox{-0.3cm}
{\psfig{figure=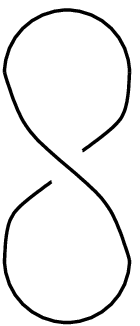,height=1cm}}$ of the trivial knot. Resolutions of the 
crossing produce two circles, respectively one circle: 
\psfrag{A2}{{$A^{\otimes 2}$}}
\psfrag{A-1}{{$A\{-1\}$}}
 $$
 {\psfig{figure=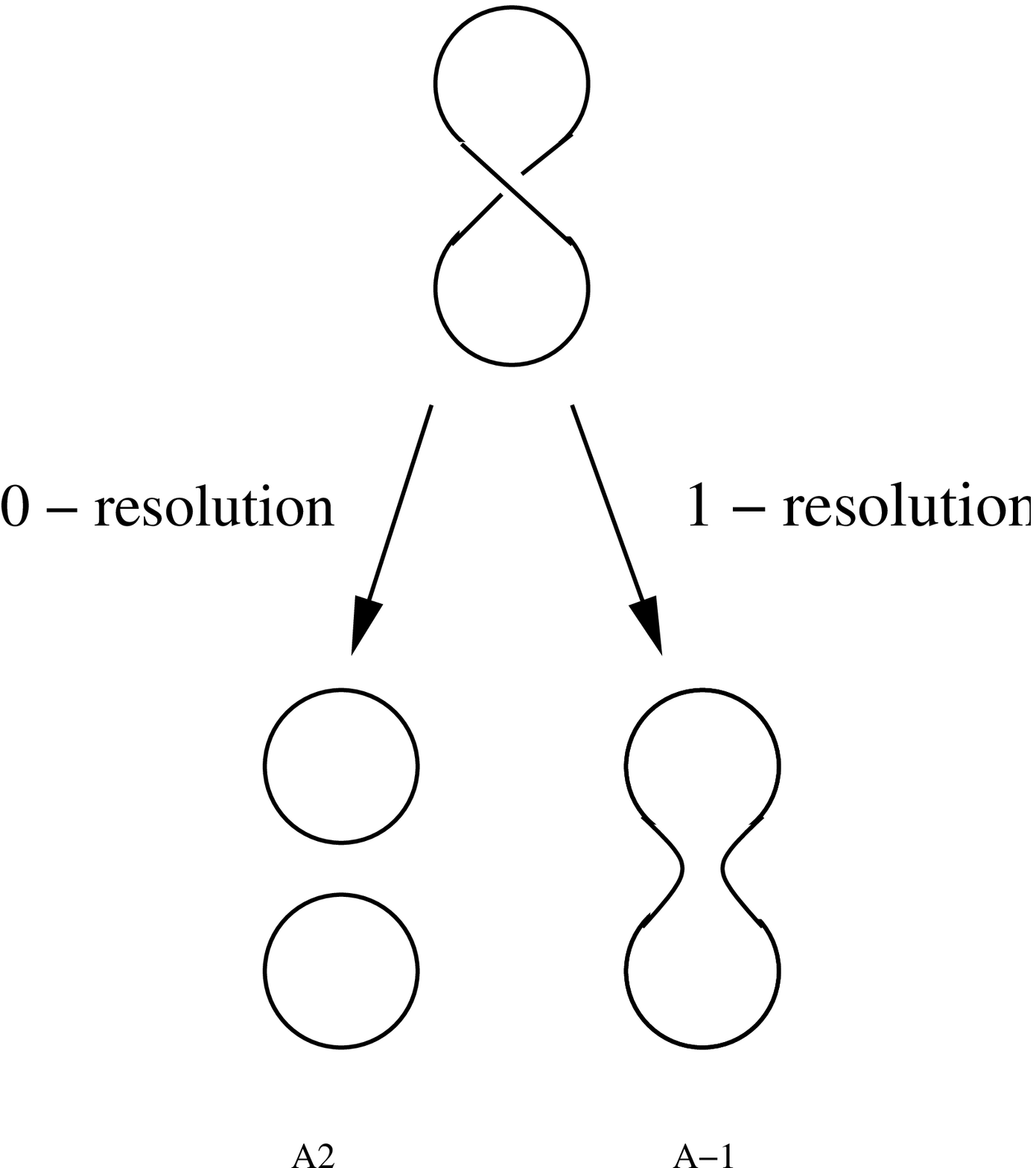,height=4cm}}$$
To the 0-resolution we assign $A^{\otimes 2},$ to the 1-resolution we assign $A\{-1\}$.  
The shift $\{-1\}$ mirrors multiplication by $q^{-1}$ in the formula 
(\ref{cr-res}). The minus sign indicates that these two terms should live in 
cohomological degrees of different parity, and the simplest guess gives us 
the complex 
\begin{equation} \label{eq-3-1} 
0 \lra A^{\otimes 2} \stackrel{m}{\lra} A\{-1\} \lra 0, 
\end{equation}
where we placed the first term in cohomological degree $0$. We denoted 
the differential in the complex by $m$ since it looks like a multiplication map. 
This map must preserve internal grading and the cohomology of the complex 
should be $\Z \oplus \Z,$ those of the unknot, since we want our theory to 
give an invariant of links and not just their diagrams. With these restrictions, 
there is very little choice available to us. We define 
$$ m({\bf 1}\otimes a) = m(a\otimes {\bf 1})= a, \ \mbox{for} \ \ a\in A, \ 
\ m(X \otimes X) = 0 .$$ 
This makes $A$ into an associative commutative unital algebra. If we 
shift the grading of $A$ up by $1,$ the multiplication becomes grading-preserving 
and we can identify $A$ with the integral cohomology ring of the 2-sphere. 
With this choice of $m$ the cohomology of the complex (\ref{eq-3-1}) is the subgroup of 
$A^{\otimes 2}$ spanned by $X\otimes 1 - 1\otimes X$ and $X\otimes X.$ 
Thus, up to overall grading shift (which we'll take care via equation (\ref{add-shift})), 
the cohomology is isomorphic to $A$. 

 \noindent
 {\bf Example 2:} the opposite kink $D= \raisebox{-0.3cm}
{\psfig{figure=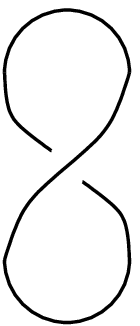,height=1cm}}$. 
Similarly to example 1 take two resolutions: 
 \psfrag{A}{{$A $}}
\psfrag{A2-1}{{$A^{\otimes 2} \{-1\}$}}
 $$
 {\psfig{figure=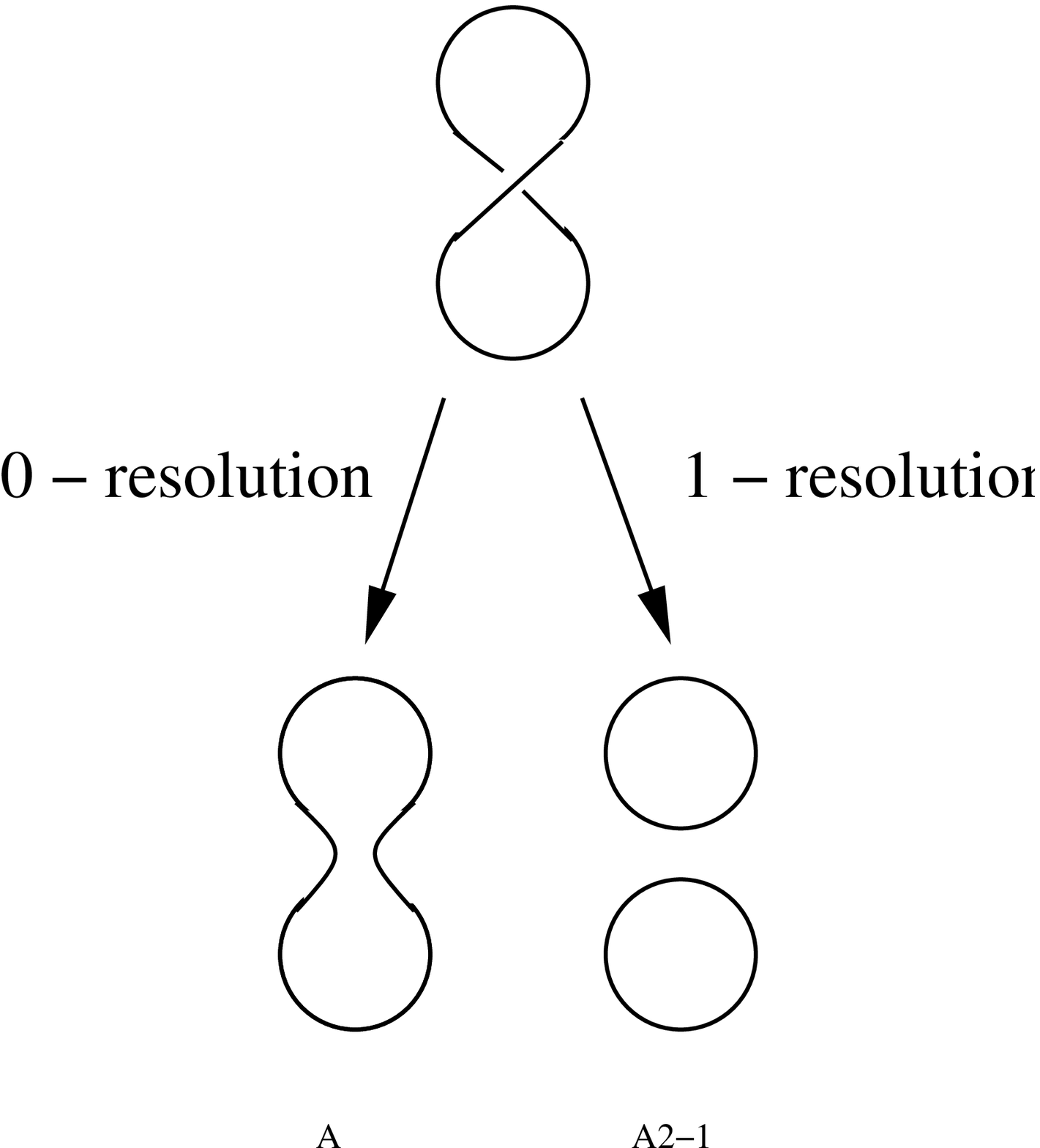,height=4cm}}$$
 and form the chain complex 
 $$0 \to A  \stackrel{\Delta}{\lra} A^{\otimes 2} \{-1\}  \to 0.$$
 It is suggestive to call the differential $\Delta$, which is the usual symbol for 
comultiplication. The bases in each chain group, sorted by the internal degree $j$, are
 $$
 \begin{array}{r|cccc}
 j  & A &    A^{\otimes 2}\{-1\} \\
 \hline
 1 & X & X \otimes X   \\
 -1 & {\bf 1} & X \otimes {\bf 1}, {\bf 1} \otimes X  \\
 -3 &0 & X \otimes {\bf 1}  
 \end{array}
 $$
and we choose $\Delta$ to be 
$$\Delta({\bf 1})= {\bf 1}\otimes X + X\otimes {\bf 1}, \ \ 
 \Delta(X) = X\otimes X.$$
The cohomology of the resulting complex is isomorphic to $A$, as a bigraded 
group, up to overall grading shift. 

We can now guess the definition of $\overline{C}(D)$ for an arbitrary $D$ with $m$ 
crossings. Each crossing has two resolutions, and the number of complete 
resolutions of $D$ is $2^m.$ Each complete resolutions is a crossingless 
diagram and has $A^{\otimes k}$ assigned to it, where $k$ is the number of 
circles. If $r_0,r_1$ are two complete resolutions that differ only in one place (near one 
crossing), with $r_0$, resp. $r_1$ being the $0$-resolution, resp. $1$-resolution there, 
then two things can happen. Either two circles of $r_0$ become one circle in 
$r_1$ or vice versa. If the first case we have a natural map $\overline{C}(r_0)\lra \overline{C}(r_1)$ 
which is $m: A^{\otimes 2} \lra A$ on the two $A$'s corresponding to these 
two circles times the identity map 
$\mathrm{Id}: A^{\otimes (k-1)}\lra A^{\otimes (k-1)}$
on the tensor product of the copies of $A$ corresponding to circles that don't change 
as we go from $r_0$ to $r_1$. Thus, the map is the composition
$$ \overline{C}(r_0) \cong A^{\otimes (k+1)} \stackrel{m\otimes \mathrm{Id}}{\lra} 
 A^{\otimes k} \cong \overline{C}(r_1).$$
In the second case, when $r_1$ has more circles than $r_0$, we have a similar 
map $\overline{C}(r_0)\lra \overline{C}(r_1)$  using $\Delta$ in place of $m$. 

Given an unoriented diagram $D$ with $m$ crossings we associate to it an $m$-dimensional 
cube with graded abelian groups $A^{\otimes k}$ (plus a grading shift) 
written in its vertices and maps $m\otimes \mathrm{Id}$, $\Delta\otimes \mathrm{Id}$ 
assigned to its edges. 

\begin{exercise} Check that every square facet of this cube is a commutative diagram. 
\end{exercise} 

We add signs to some of the edge maps so that each facet anticommutes and collapse 
the $m$-dimensional cube into a complex of graded abelian groups. The terms in the 
complex are given by direct sums of graded abelian groups sitting in vertices contained 
in a given hyperplane perpendicular to the main diagonal. We place the first term 
in cohomological degree $0$ and the last in cohomological degree $m$. The result 
is a complex of graded abelian groups, denoted $\overline{C}(D)$, with a grading-preserving 
differential. Let us look at an example. 

\noindent 
{\bf Example 3} A 3-crossing diagram
 $D=\raisebox{-0.8cm}{\psfig{figure=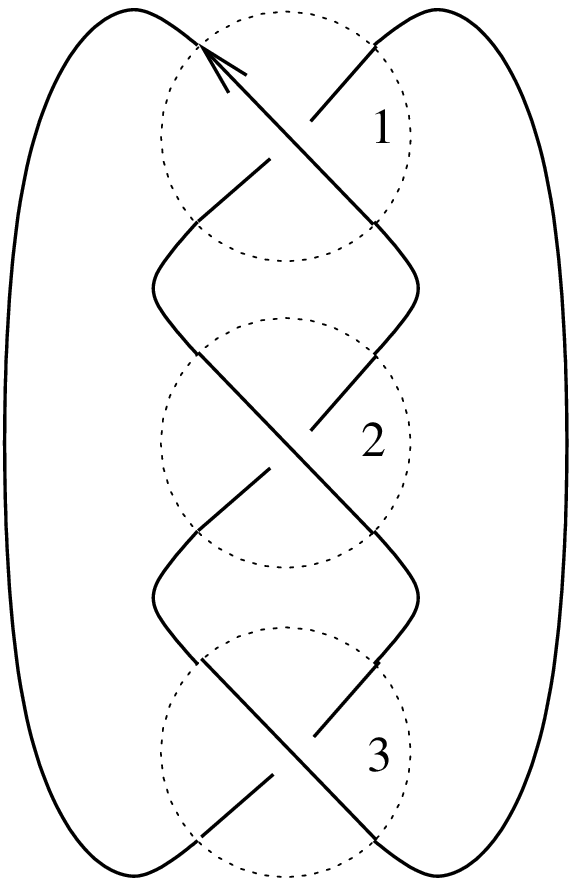,height=2cm}}$ of a trefoil. 
We ignore the orientation of $D$ and label the crossings by $1, 2,3 $; 
the resolutions are as follows. Arrows parallel the arrow labelled $d_i$ 
correspond to modifications at the $i$-th crossing. The sequence $010$ 
written above the bottom left diagram indicates that the first and the third 
crossings are modified via the $0$-resolution and the second crossing--via 
the $1$-resolution, etc. 
$$
\begin{diagram}
\node{}\node{\stackrel{001}{\raisebox{-0.6cm}{\psfig{figure=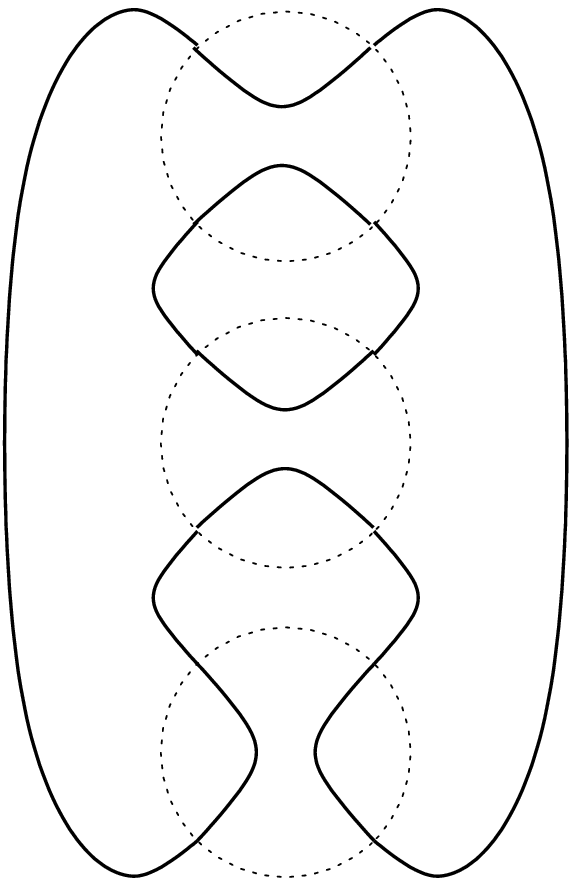,height=1.5cm}}}} 
\arrow[2]{e} \arrow{s,-} \node{}
\node{\stackrel{101}{\raisebox{-0.6cm}{\psfig{figure=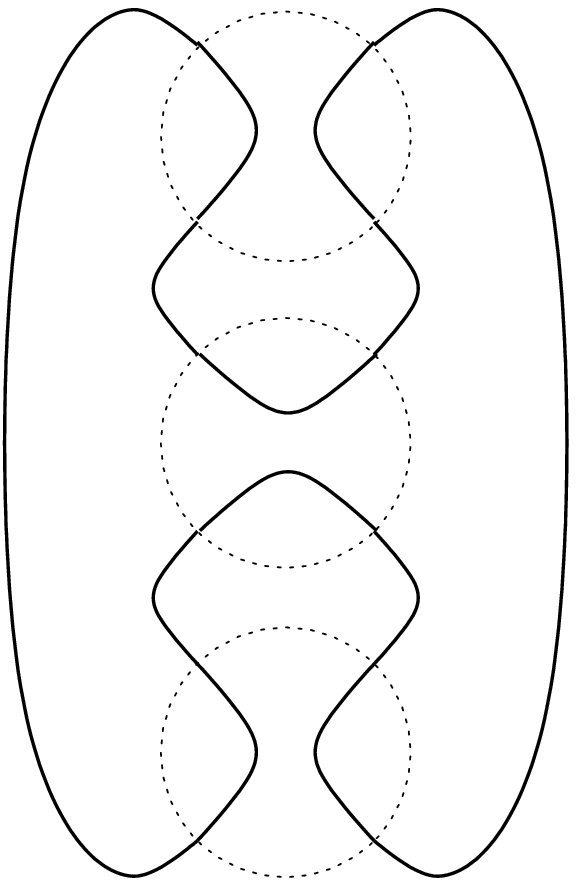,height=1.5cm}}}}
\arrow[2]{s}  \\
\node{\stackrel{000}{\raisebox{-0.6cm}{\psfig{figure=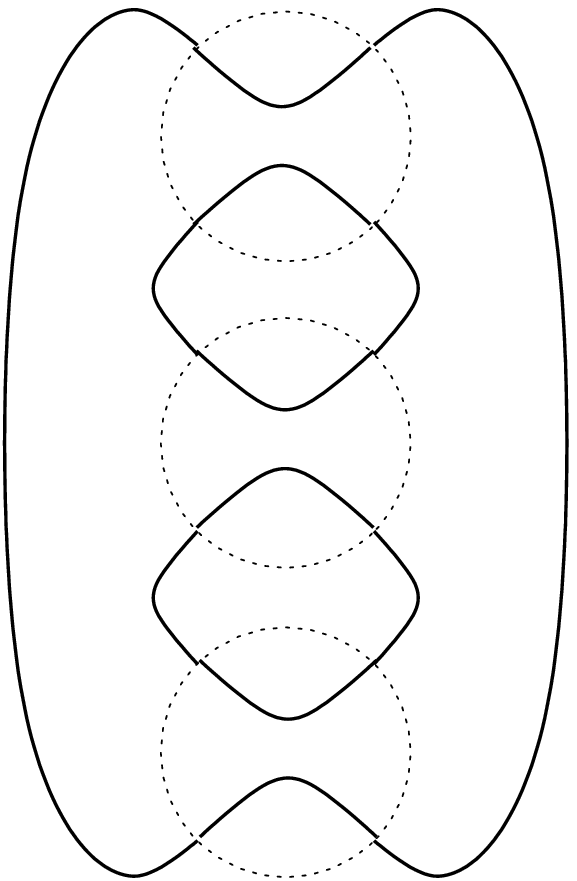,height=1.5cm}}}} 
\arrow[2]{e,t,1}{d_1} \arrow{ne,l}{d_3} 
\arrow[2]{s,l}{d_2} \node{}\arrow{s} \node{\stackrel{100}{\raisebox{-0.6cm}
{\psfig{figure=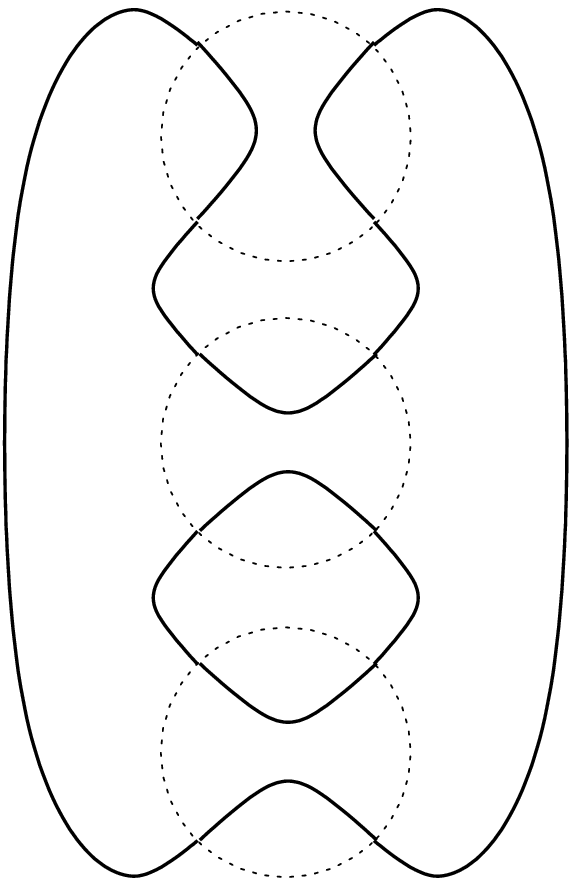,height=1.5cm}}}} 
\arrow{ne} \arrow[2]{s} \node{} \\
\node{} \node{\stackrel{011}{\raisebox{-0.6cm}
{\psfig{figure=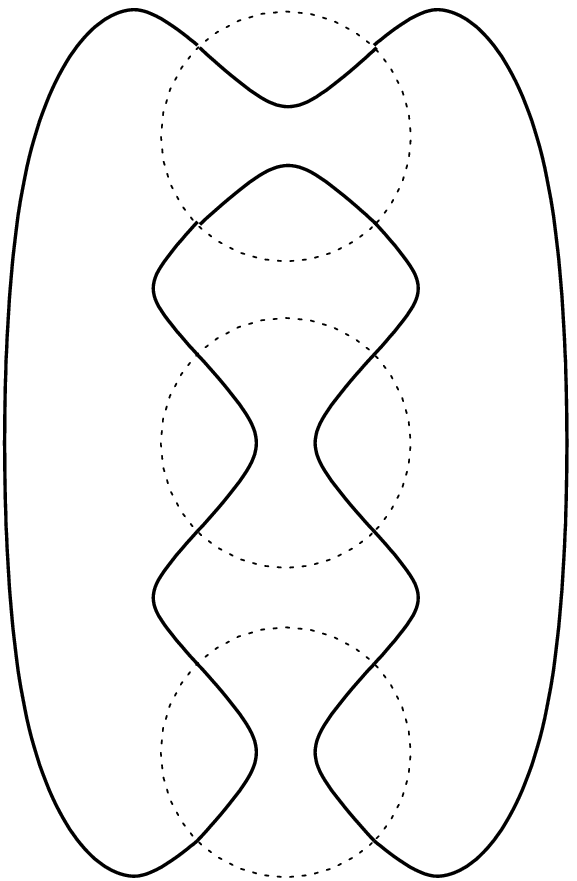,height=1.5cm}}}} 
\arrow{e,-} \node{} \arrow{e}
\node{\stackrel{111}{\raisebox{-0.6cm}{\psfig{figure=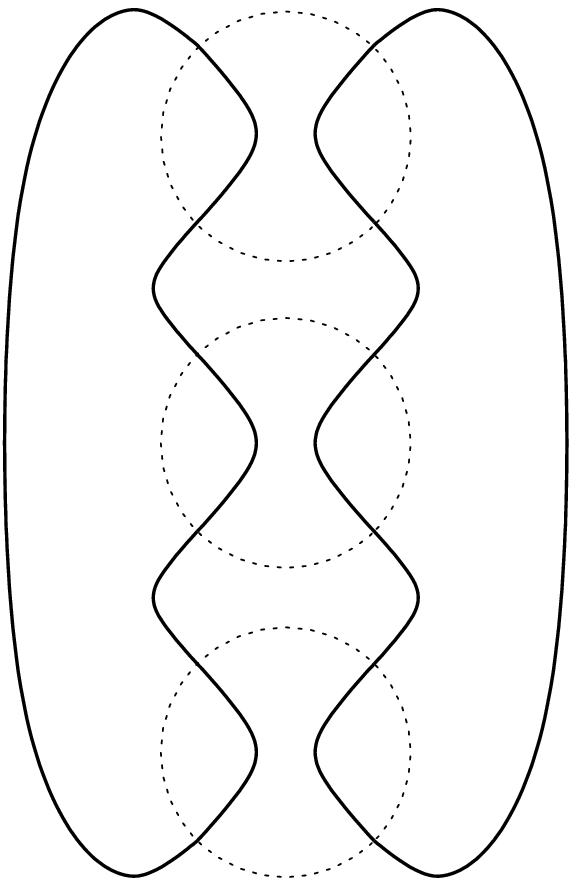,height=1.5cm}}}} \\
\node{\stackrel{010}{\raisebox{-0.6cm}{\psfig{figure=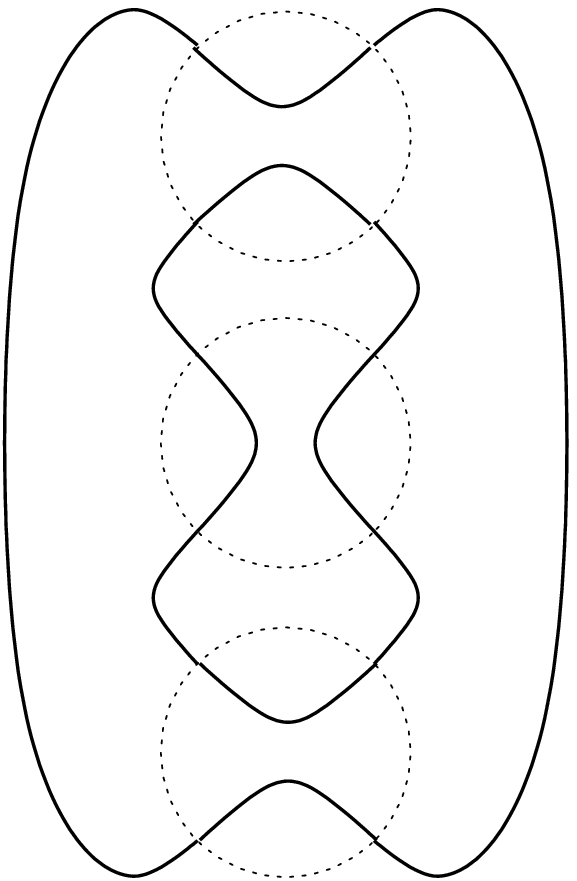,height=1.5cm}}}} 
\arrow[2]{e} \arrow{ne} \node{}
\node{\stackrel{110}{\raisebox{-0.6cm}{\psfig{figure=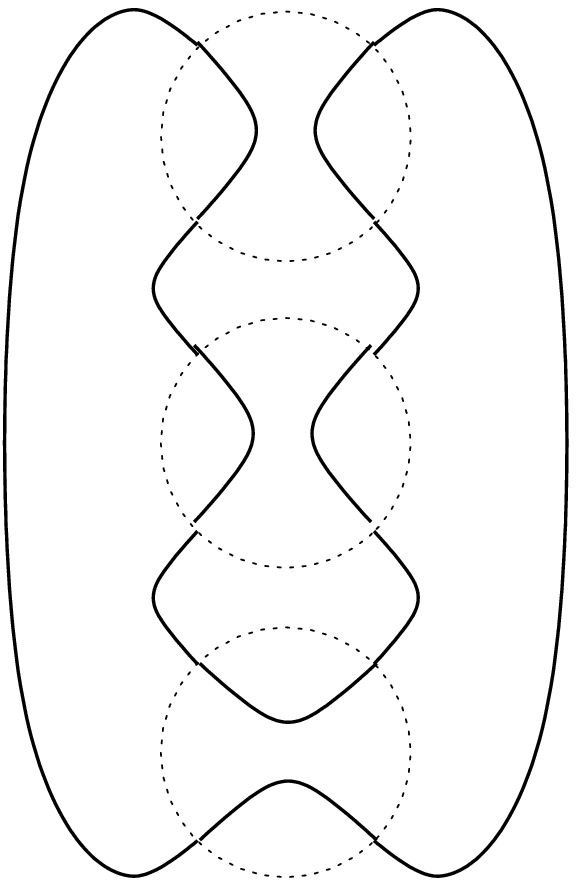,height=1.5cm}}}} 
\arrow{ne} \node{} 
\end{diagram}
 $$
 The corresponding groups and maps are (we write $m$ instead of $m\otimes 
\mathrm{Id}$ on top left arrows)  
 $$
 \begin{diagram}
 \node{}\node{A^{\otimes 2}\{-1\}} \arrow[2]{e,t}{m} \arrow{s,-} \node{}
\node{ A\{-2\}}\arrow[2]{s,l}{\Delta}  \\
\node{A^{\otimes 3}} \arrow[2]{e,t,1}{m} \arrow{ne,l}{m} 
\arrow[2]{s,l}{m} \node{}\arrow{s,r}{\Delta} \node{A^{\otimes 2}\{-1\}} 
\arrow{ne,l}{m} \arrow[2]{s,l,1}{m} \node{} \\
\node{} \node{A\{-2\}} \arrow{e,-}  \node{} \arrow{e,t}{m}
\node{A^{\otimes 2}\{-3\}} \\
\node{A^{\otimes 2}\{-1\}} \arrow[2]{e,b}{m} \arrow{ne,l}{m} \node{}
\node{A\{-2\},} \arrow{ne,r}{\Delta} \node{} 
\end{diagram}
$$
where each $m$, $\Delta$ is applied  according to the topology change: 
\psfrag{m}{$m$}
\psfrag{d}{$\Delta$}
\psfrag{A}{$A$}
\psfrag{AA}{$A^{\otimes 2}$}
$$
\raisebox{-0.6cm}{\psfig{figure=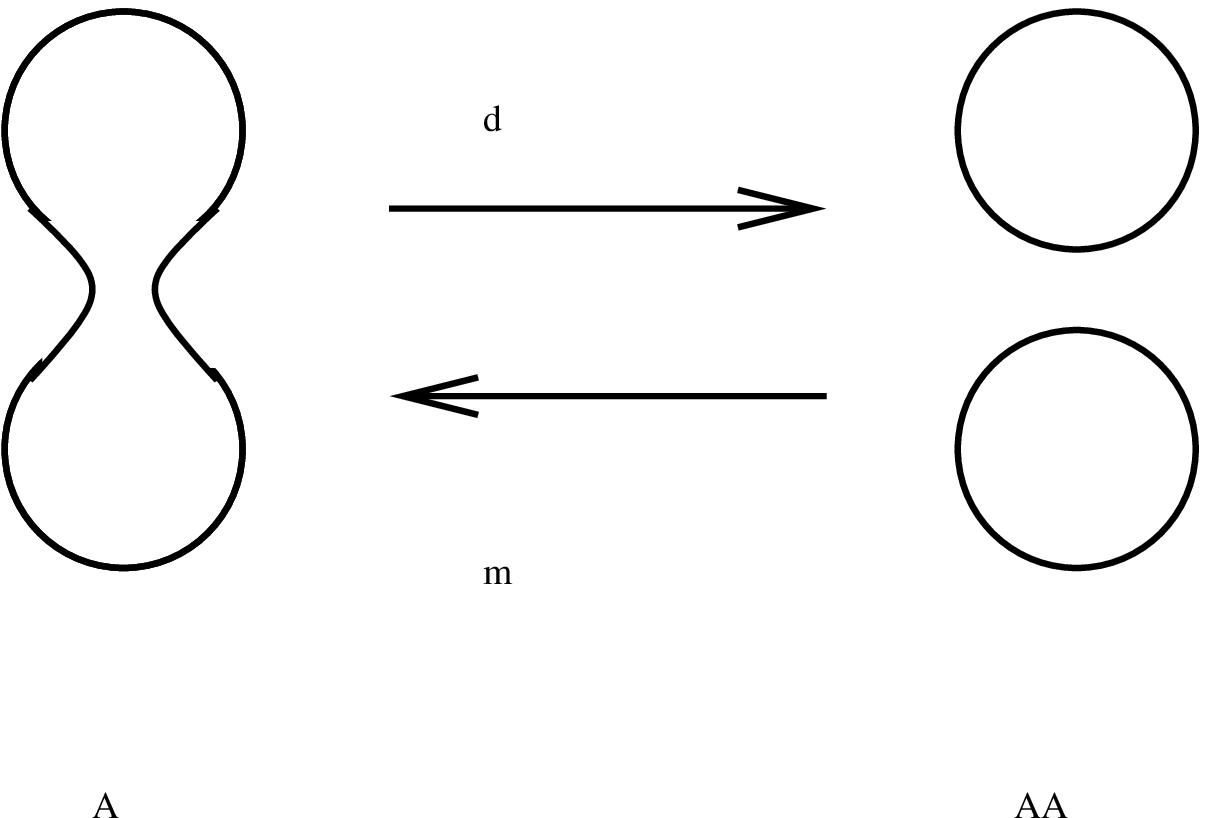,height=1.5cm}}
$$
For tensor factors for the circles that do not change, we apply $\mathrm{Id}.$ 
We add minus signs to make each square anticommutative and pass to the total complex 
of the cube. Due to anti-commutativity of each facet, $d^2=0$ holds. 
The total complex has the form 
$$
\overline{C}(D)=\left( 0 \to A^{\otimes 3} \stackrel{d}{\to} 
\begin{array}{c} 
 A^{\otimes 2} \{-1\} \\
 \oplus \\
  A^{\otimes 2} \{-1\}  \\
  \oplus \\
   A^{\otimes 2} \{-1\}  
    \end{array}
\stackrel{d}{\to}
\begin{array}{c} 
 A \{-2\} \\
 \oplus \\
  A \{-2\}  \\
  \oplus \\
    A \{-2\}   
    \end{array}
    {\to}
    A^{\otimes 3} \{-3\}  \stackrel{d}{\to} 0 \right). $$
 For an arbitrary oriented diagram $D$ we define the complex $C(D)$ by 
shifting the complex $\overline{C}(D)$ as in the formula (\ref{add-shift}). 
An even more elementary definition, avoiding explicit use of tensor powers, 
can be found in Viro~\cite{Vir}, together with other interesting observations. The complex 
$C(D)$ starts in homological degree $-x(D)$, where $x(D)$ is the number 
of negative crossings and ends in homological degree $y(D)$, the number of 
positive crossings of $D$. 
For $D$ in example 3, $x(D)=3$ and $y(D)=0$. 

Finally, define $H(D)$ as the cohomology of the complex $C(D)$. Note 
that $H(D)$ is bigraded and, from the construction, the Euler characteristic 
of $H(D)$ is the Kauffman bracket (the Jones polynomial) of $L$.  
 \begin{exercise} Compute $H(D)$ for the above diagram of the trefoil.  
 \end{exercise}
 The following holds \cite{Kh1}: 
 \begin{theorem}
 If two diagrams $D_1$ and $D_2$ are related by a chain of Reidemeister moves,
the complexes of graded abelian groups  $C(D_1)$  and $ C(D_2)$ are homotopy 
equivalent and homology groups $H(D_1)$ and $H(D_2)$ are isomorphic. 
\end{theorem} 

Define the link homology $H(L):=H(D)$ for a diagram $D$ of $L.$ 
Homology groups $H(L)$ are known as Khovanov homology. We have
 $$J(L)=\chi(H(L))=\sum_{i,j\in \Z} (-1)^i q^j {\rm rk} H^{i,j}(L).$$
Notice that the above theorem only says that the isomorphism class of $H(L)$ 
as a bigraded group is an invariant of $L.$ 
We'll discuss the issue of functoriality under link 
diffeomorphisms and, more generally, link cobordisms, in Lecture~5. 
 \begin{exercise}
 Let $D_1$, $D_2$ be diagrams of links $L_1$, $L_2$. 
 Check that $C(D_1 \sqcup D_2)\cong C(D_1) \otimes C(D_2),$
 and derive the K\"unneth formula for the homology of the disjoint 
union $L_1\sqcup L_2$. This formula categorifies the 
multiplicativity property of the Jones polynomial,  
$J(L_1\sqcup L_2) = J(L_1) J(L_2)$, in our normalization. 
 \end{exercise}


\subsection{Properties and examples} 

At least three programs, by D.~Bar-Natan~\cite{BN1}, A.~Shumakovitch~\cite{Shu}, 
D.~Bar-Natan and J.~Green~\cite{BN3} are available for computation of 
$H(L)$. The following tables, provided to us by A.~Shumakovitch, 
show Khovanov homology of several knots. 
Given a diagram $D$ of a knot, the complex $C(D)$ (and, hence, its homology)  
is nontrivial in odd internal degrees only. Thus, even 
$q$-degrees are not shown in the tables.
\begin{center} \hspace{-0.2cm} 
\includegraphics[scale=0.65]{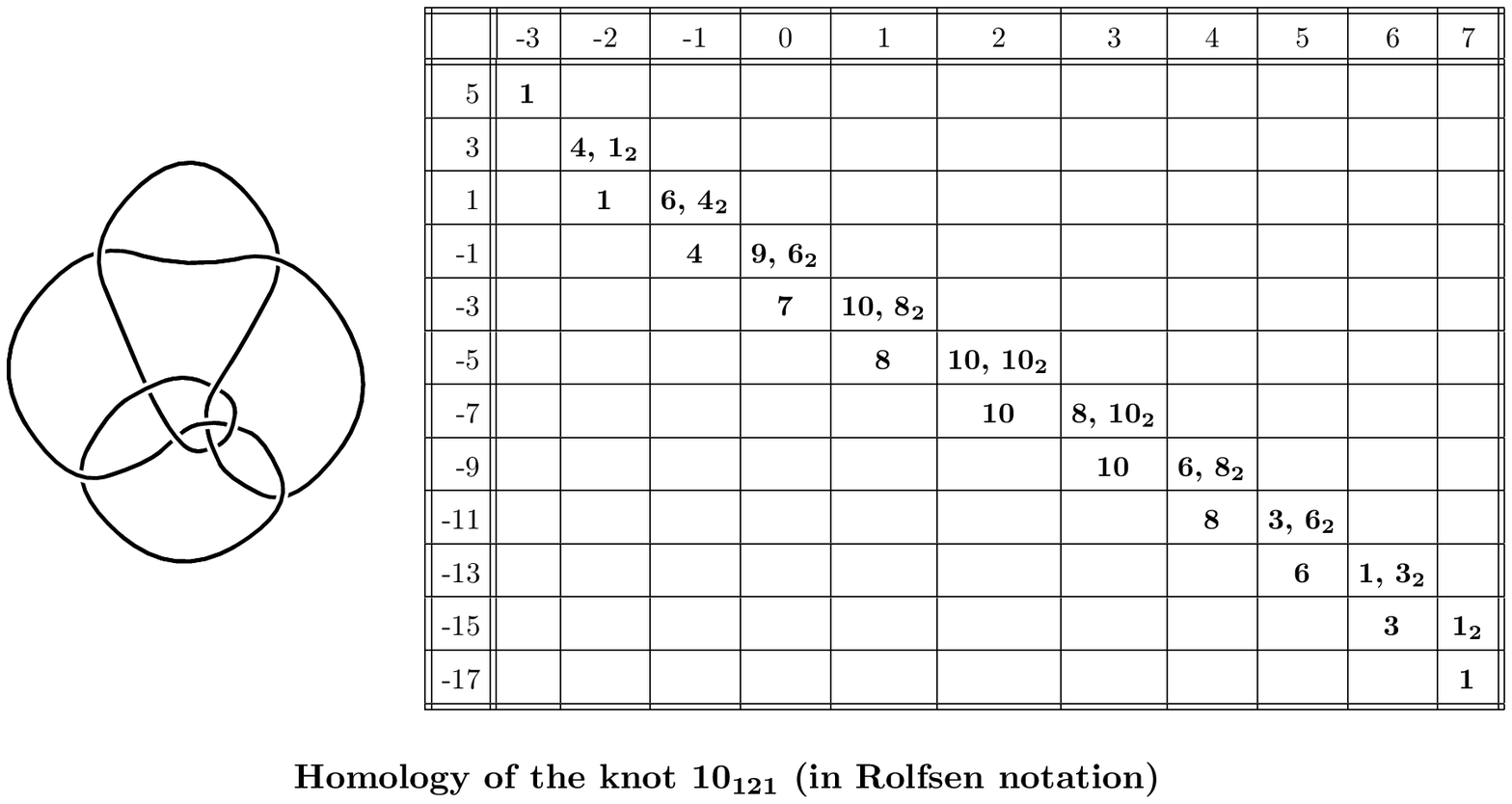}
\end{center} 
The first table displays the homology of an alternating 10-crossing knot $10_{121}$. 
A  single integer entry $a$ says that the homology group in that bidegree 
is free of rank $a$.  
Two integers $a,b$ separated by a comma indicate that the homology group is 
the direct sum of $\Z^a$ and $(\Z/2)^b$. For instance, $H^{5, -13}\cong \Z^6$ 
and $H^{1,-3}\cong \Z^{10}\oplus (\Z/2)^8.$
Calculations done 
by Shumakovitch~\cite{Shu2} show that the only torsion in the homology 
of knots with 14 or fewer crossings is $\Z/2$-torsion. Shumakovitch found several 
knots with 15 and 16 crossings whose homology contains a copy of $\Z/4$. 
One of these knots is the $(4,5)$-torus knot. For more information and results 
on torsion see~\cite{Shu2}, \cite{AP}. 

Looking at the table, we notice that the homology spans horizontal 
degrees from -3 to 7, and the knot has homological width 10, equal to the 
number of crossings of the diagram. By the \emph{homological width} of $L$ we 
mean the difference between the largest $i$ 
such that $H^{i,j}(L)\not= 0$ for some $j$ and the minimal $i$ with the same property. 
The homological width of $L$ gives a lower bound on the crossing number of $L$ 
(the smallest number of crossings in a planar diagram of $L$). 

\begin{exercise} Determine the geometric condition on $D$ which ensures that 
$C(D)$ has nontrivial homology in the leftmost and in the rightmost degrees
(such diagrams $D$ are called \emph{adequate}). 
\end{exercise} 

Thus, if $L$ has an $n$-crossing adequate diagram, the crossing number 
of $L$ is $n$. This was originally proved by Thistlethwaite~\cite{This} using 
the 2-variable Kauffman polynomial (not to be confused with the Kauffman 
bracket, which is a one-variable polynomial). One of the Tait conjectures 
about the crossing number of alternating links (originally proved with the help 
of the Jones polynomial~\cite{Kau}, \cite{Mu}) follows from this result. 
This alternative approach to the Thistlethwaite theorem does not 
use the internal grading on the homology, only the horizontal grading, 
almost invisible on the level of Euler characteristic. 

Apparently, the only known explicit relation between the 2-variable Kauffman polynomial 
and Khovanov homology is that they both specialize to the Jones polynomial. 
Yet, both can be used to prove the Thistlethwaite theorem and to give 
upper bounds on the Thurston-Bennequin number of Legendrian links, see 
L.~Ng~\cite{Ng} and references therein. 

Diagonal width of the homology gives a lower bound on the Turaev genus 
of a knot~\cite{Mnt}. For relations between homology and contact topology 
see~\cite{Plam} and references therein. 

The total rank of the complex $C(D)$ grows very fast as a function of the size of $D$. 
Thistlethwaite's spanning tree model for the Jones polynomial admits a categorification~\cite{Wh}, 
\cite{ChK}, giving a complex of much smaller rank which also computes $H(D)$, 
but no combinatorial formula for the differential of the resulting complex is 
known, except in very special cases. 

The homology of $10_{121}$ occupies two adjacent diagonals. It is easy to see
that the homology cannot lie on just one diagonal, so two is the minimum. 
E.~S.~Lee~\cite{Lee1} proved that the homology of any alternating knot $L$ lies 
on two adjacent diagonals consisting of $(i,j)$ with $2i-j=c\pm 1$ where $c$ 
is the signature of $L$. Computer calculations show this to be true for most 
non-alternating knots with 11 or fewer crossings as well~\cite{BN1}, \cite{Shu}. 
A partial explanation of this phenomenon was provided by C.~Manolescu 
and P.~Ozsv\'ath~\cite{MO} by extending Lee's result to quasi-alternating knots.

The next table shows the homology of the alternating knot $10_{123}$. 
This knot is amphicheiral, that is, isomorphic to its mirror image.  

\vspace{0.2in} 

\begin{center}\hspace{-0.2cm}
\includegraphics[scale=0.65]{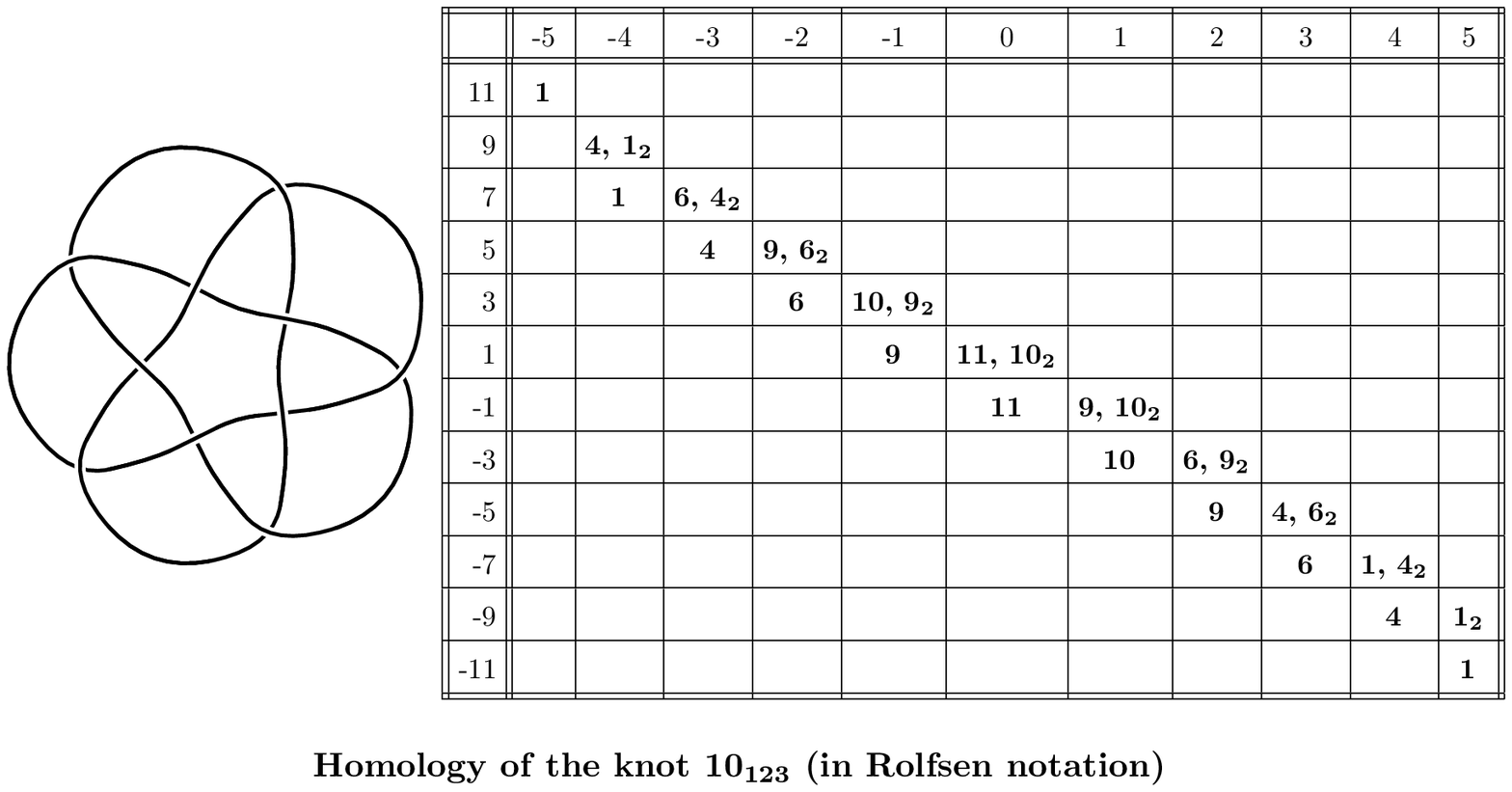}\end{center}

\vspace{0.2in} 

For a link $L$ the mirror image $L^!$ is given by reversing the orientation 
of the ambient $\R^3$ in which the 1-manifold is embedded. 
By reversing all crossings of a diagram $D$ of $L$ we obtain the diagram $D^!$ 
of $L^!$. 
 \begin{exercise}
Construct an isomorphism of complexes 
$$C(D^{!})\cong \Hom_{\Z} (C(D), \Z).$$
\end{exercise}
Thus, the dual $C(D)^{\ast}$ of the of complex $C(D)$ is naturally isomorphic to $C(D^!)$. 
This duality takes $C^{i,j}(D)$ to the dual of the free abelian 
group $C^{-i,-j}(D)$. Passing to homology, we see that the free part 
of $H^{i,j}(L)$ becomes the dual of the free part of $H^{-i,-j}(L^!)$. 
By the free part of a finitely-generated abelian group $G$ we mean the 
quotient $G/\text{Tor}(G)$ of $G$ modulo torsion. In particular, abelian 
groups $H^{i,j}(L)$ and  $H^{-i,-j}(L^!)$ have the same rank. 
The torsion of $H^{i,j}(L)$ is the dual of the torsion of $H^{-i+1,-j}(L^!)$
(notice the grading shift), in particular, the two torsion groups have the same rank. 
The reader can see this duality in the above homology table of 
$10_{123}\cong 10_{123}^!.$ 
The first integer in the $(i,j)$-entry equals the first integer in the $(-i,-j)$-entry. 
The torsion, only present on the upper diagonal, stays on the upper diagonal after 
dualization, due to the shift by $1$. For instance, $H^{5,-9}\cong \Z/2$ 
becomes, after dualization, the torsion subgroup $\Z/2$ of $H^{-4,9}$. 

Recall that, for the singular chain complex $C(X)$ of a topological space $X$ 
which computes the homology groups of $X$, the dual complex $C(X)^{\ast}$ 
computes the cohomology groups of $X$. In the link homology framework, 
the duality works in a different way--the underlying link is converted to its 
mirror image. This manifests our terminological imperfection in 
calling groups $H(L)$ the \emph{homology} groups of $L$. We would be equally 
justified in calling them \emph{cohomology} groups. In the next two lectures we'll 
discuss functoriality of $H$. To a link cobordism $S$ from $L_1$ to $L_2$ 
we assign a homomorphism $H(S): H(L_1) \lra H(L_2)$, which, over all $S$, 
gives a covariant functor. However, there is an equally natural construction 
that assigns to $S$ the homomorphism going in the opposite direction, 
producing a contravariant functor. We see that $H(L)$ exhibits both covariant 
and contravariant behaviour and, in a flexible terminological environment, 
we are free to call $H(L)$ either homology of cohomology. Another solution 
is to call $H(L)$ bivariant (co)homology groups.

\vspace{0.2in} 

\begin{center}\includegraphics[scale=0.65]{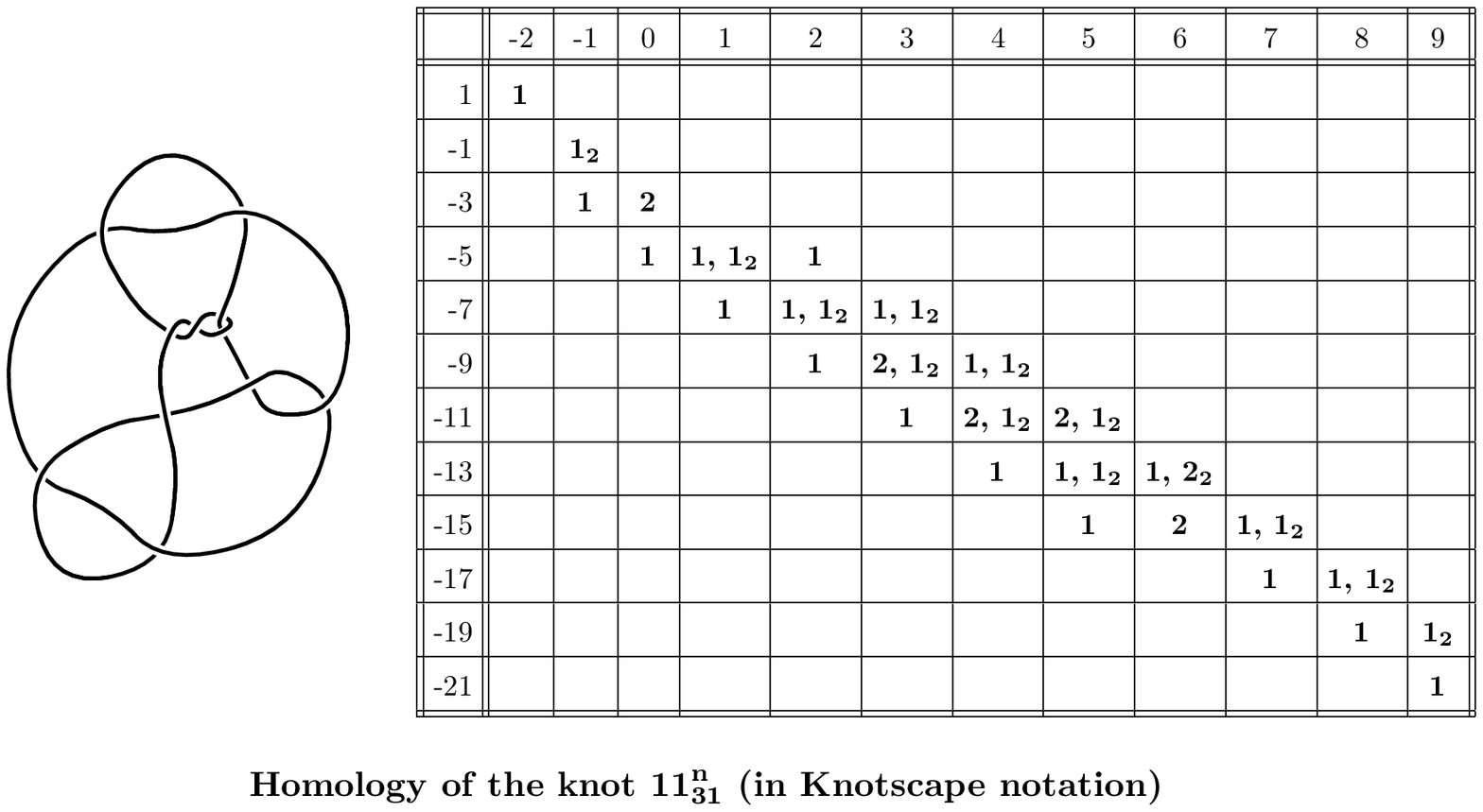}
\end{center}

\vspace{0.2in} 

In the third table we see an 11-crossing non-alternating knot whose homology occupies 3 
adjacent diagonals. This knot is adequate, and the width of the homology is 
11. This diagram $D$ has 2 negative and 9 positive crossings,  
and the groups are bounded by homological degrees $-2$ and $9$. 
Unlike the previous two examples, each homology group has small rank (at most 
two). 

\vspace{0.2in} 

\begin{center} \includegraphics[scale=0.75]{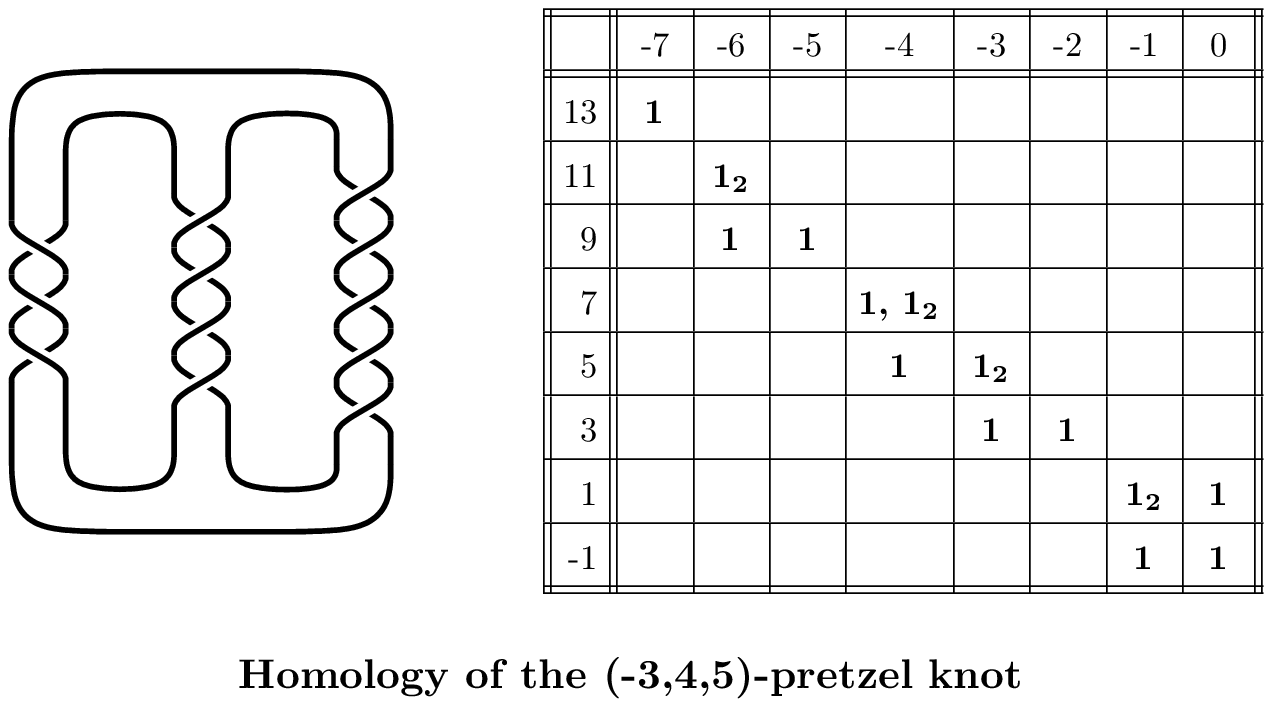}
\end{center} 

\vspace{0.2in} 

The fourth table shows the (-3,4,5)-pretzel knot and its homology. This time 
the rank of each group is at most one. 
Homological width of this knot equals 7, much less then 12, the crossing number 
of the knot.

\vspace{0.2in} 

The simplest knot known to have odd torsion in its homology is the (5,6)-torus 
knot, see the table below. It has a $\Z/3$-summand in bidegree (14,-43) and 
$\Z/5$-summands in bidegrees (11,-35) and (12,-49). 

This is a positive knot (all the crossings look like 
$\raisebox{-0.2cm}{\psfig{figure=lec3.2.eps,height=0.5cm}}$ in the 
planar diagram $D$ given by the closure of the braid $(\sigma_1 \sigma_2 
\sigma_3 \sigma_4 )^6$), so
it has homology groups in nonnegative homological degrees only. You can check 
by hand the correctness of the table in the  homological degree $0$ by computing the 
kernel of the leftmost differential 
$$ 0 \lra C^0(D) \stackrel{d}{\lra} C^1(D)\lra \dots  $$ 
in the diagram $D$.  For various results on homology of positive and torus knots and links 
see~\cite{Sto1}, \cite{Sto2}.

\begin{center}
\includegraphics[scale=0.75]{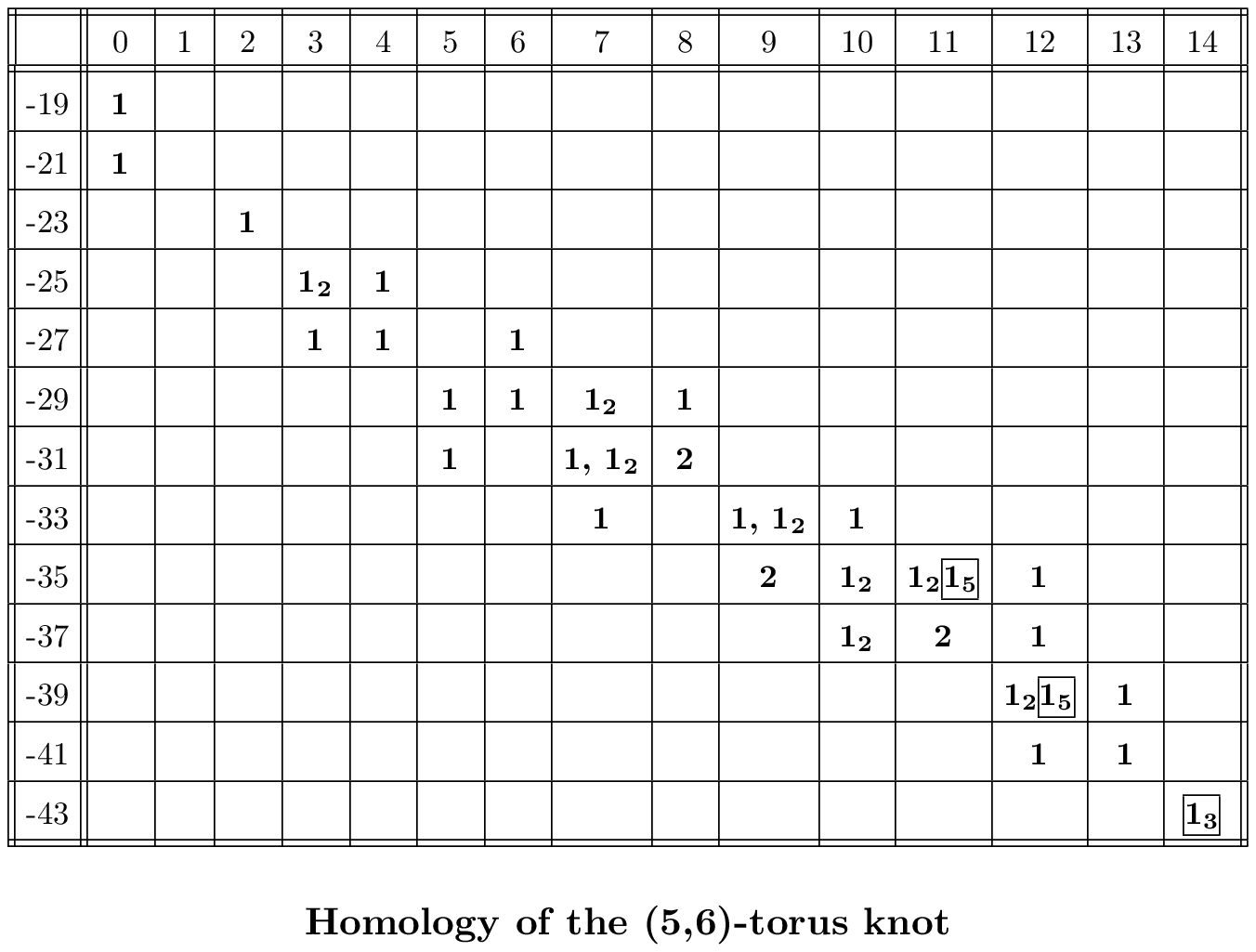}
\end{center}


\section{Flat tangles and bimodules}

 \subsection{Two-dimensional TQFTs and Frobenius algebras}
   \begin{definition}
A 2-dimensional TQFT (topological quantum field theory) is a tensor functor 
from the category of two-dimensional oriented cobordisms between oriented closed 
one-manifolds to an additive tensor category.  
 \end{definition}
 Such functor $F$ satisfies $F( X \sqcup Y) \cong F(X) \otimes F(Y)$ for 1-manifolds 
$X, Y$ and 
$F(f\otimes g) = F(f) \otimes F(g)$  for cobordisms $f,g$.  
 We won't define here additive tensor categories, but rather provide examples: 
\begin{itemize} 
\item the category of vector spaces over a field, 
\item the category of graded vector spaces over a field,  
\item the category of free modules over a commutative ring $R$, 
\item the category of complexes of free modules over a commutative ring 
$R$ modulo chain homotopies. 
\end{itemize} 
In these examples the tensor structure is the obvious one (the tensor product is taken over 
$R$ in the last two).  We restricted to free $R$-modules since, in homological algebra,
 for general $R$-modules 
the tensor product $M\otimes_R N$ must be redefined via a free or projective resolution 
of $M$ or $N$. 

A 2-dimensional TQFT $F$ with values in the category of free $R$-modules 
must assign $R$ to the empty 1-manifold, $F(\emptyset)=R$, and 
some free $R$-module $A$ to the circle, $F( \raisebox{-0.1cm}
{\psfig{figure=lec3.4.eps,height=0.6cm}})= A$. Then 
$$F(   \underbrace{  \raisebox{-0.1cm}{\psfig{figure=lec3.8.eps,height=0.5cm}} 
\cdots  \raisebox{-0.1cm}{\psfig{figure=lec3.8.eps,height=0.5cm}}}_{j})=A^{\otimes j}$$
since the functor $F$ is tensor.  To the identity cobordism $F$ assigns the identity map 
$$F\left(  \raisebox{-0.4cm}{\psfig{figure=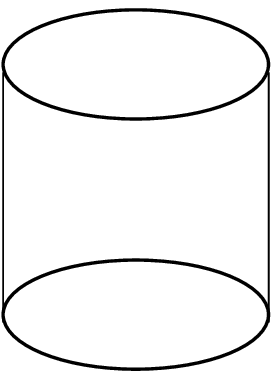,height=1cm}}\right)=
{\rm id}_A.$$
To the \emph{inverted pants} cobordism 
$$\raisebox{-0.8cm}{\psfig{figure=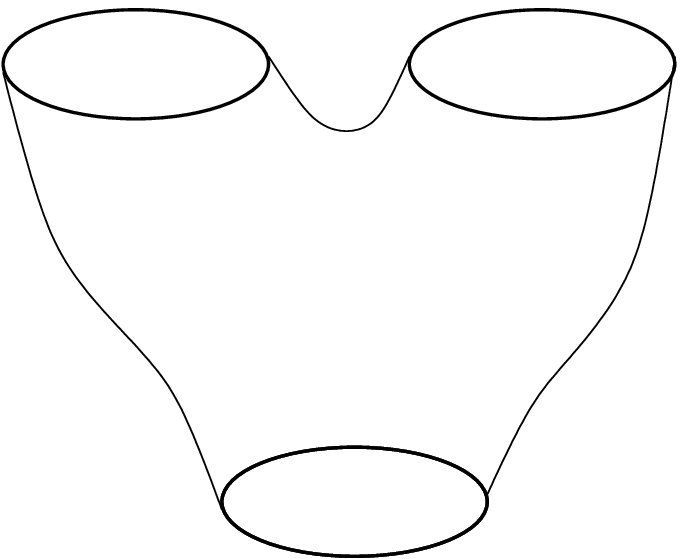,height=1.7cm}} 
 \begin{array}{l}
 A^{\otimes 2} \\
 \hspace{0.1cm} \downarrow m \\
 A
 \end{array}
 $$
$F$ assigns a homomorphism $m: A^{\otimes 2}\lra A$, associative 
due to the equality of cobordisms 
 $$
 \raisebox{-0.8cm}{\psfig{figure=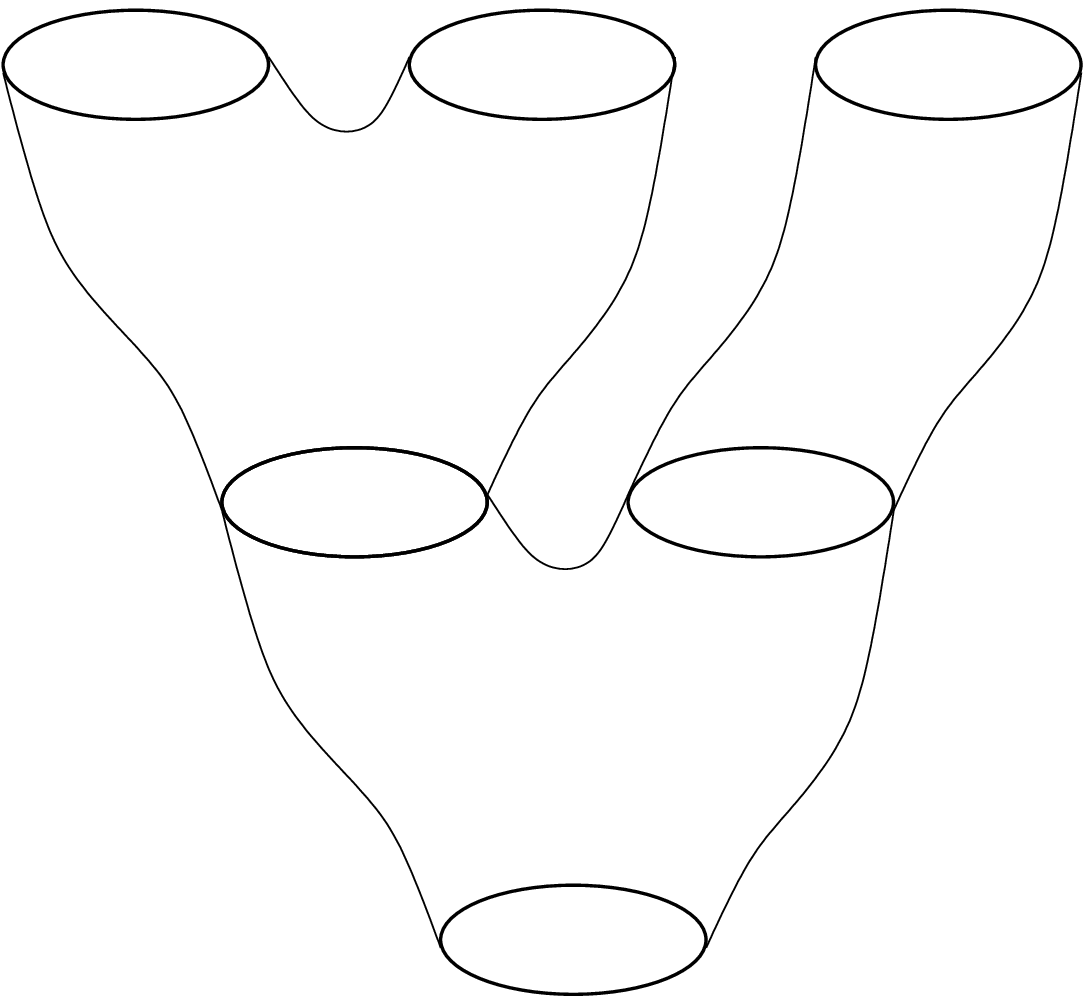,height=2.5cm}}  \cong
 \raisebox{-0.8cm}{\psfig{figure=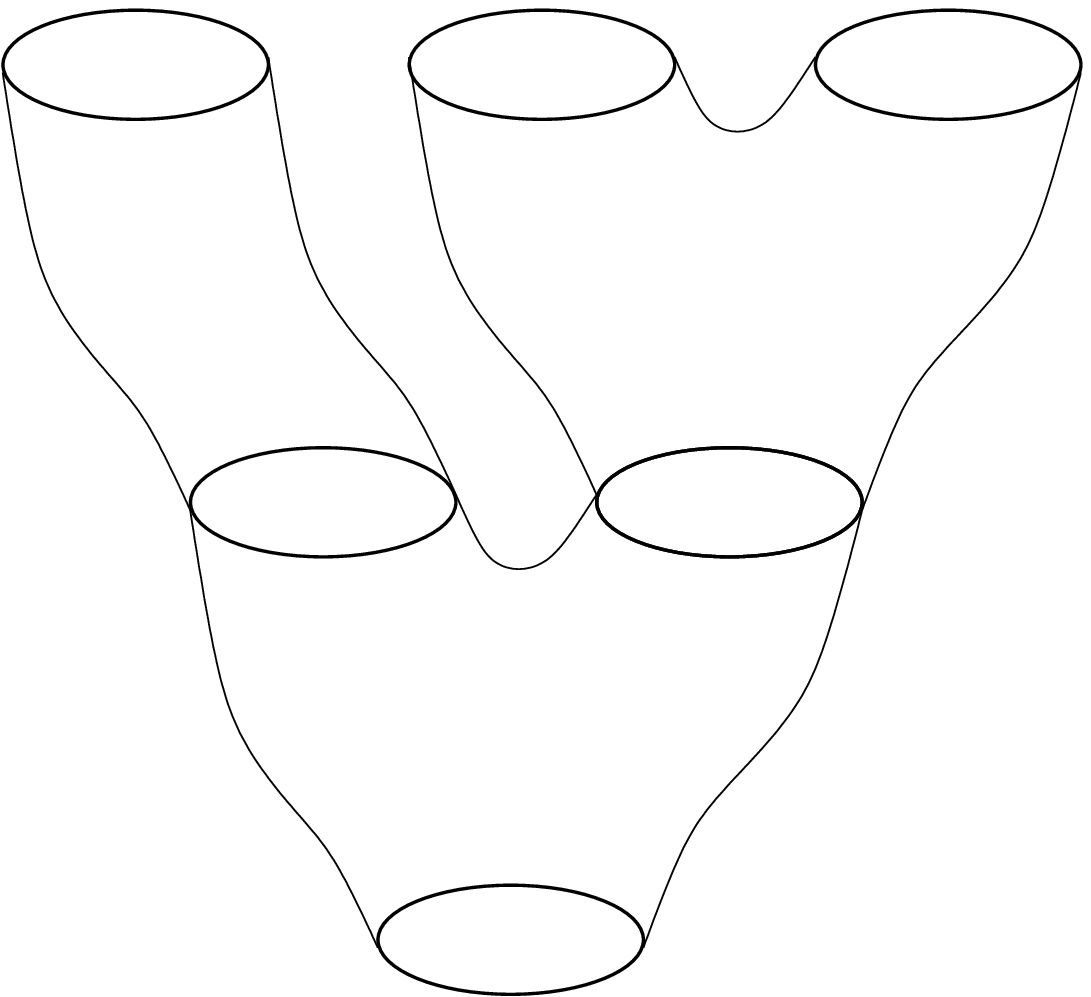,height=2.5cm}} 
 $$
\vspace{0.05in}
 The two upper legs of the multiplication cobordisms  may be permuted without 
changing the diffeomorphism type of the cobordism. Thus $m$ is commutative as well. 

The following three cobordisms induce three more maps between tensor power of $A$, 
denoted $\Delta, \iota$ and $\epsilon$, respectively, which make $A$ into a 
commutative Frobenius algebra over $R$. 
 \begin{center}
 \raisebox{-0.8cm}{\psfig{figure=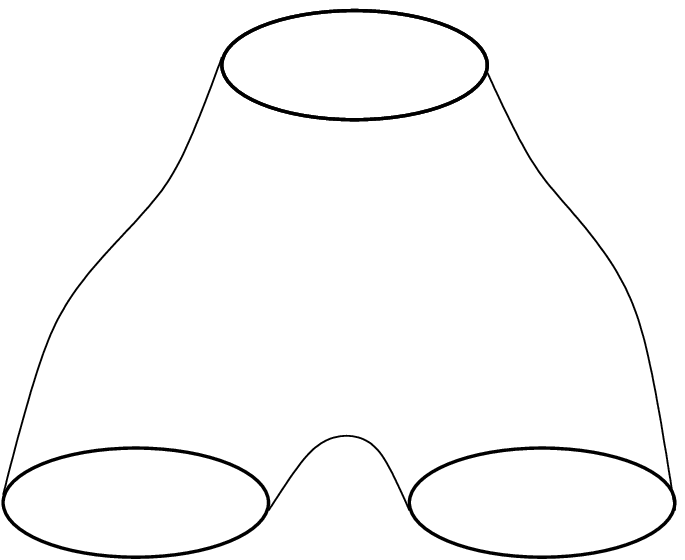,height=1.7cm}} 
$ \begin{array}{l}
 A \\
\hspace{0.1cm} \downarrow \Delta \\
 A^{\otimes 2}
 \end{array}$
 \hspace{1cm}
 \raisebox{-0.8cm}{\psfig{figure=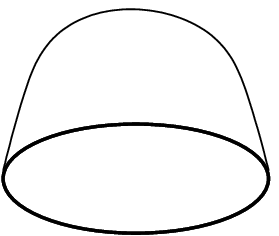,height=1cm}} 
$ \begin{array}{l}
 k   \\
 \downarrow  \iota \\
A 
 \end{array}$ 
 \hspace{1cm}
  \raisebox{-0.8cm}{\psfig{figure=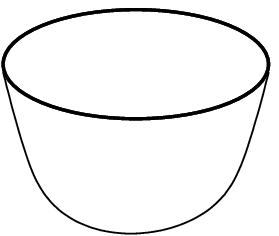,height=1cm}} 
$ \begin{array}{l}
A \\
 \downarrow {\epsilon} \\
k
 \end{array}$
 \end{center}
\vspace{0.05in}
It is well known \cite{Kock}, \cite[Section 4.3]{BK}, \cite{Kad}
that two-dimensional TQFTs with the target category being the 
category of free modules over a commutative ring $R$ are classified by 
commutative Frobenius $R$-algebras $A$. The Frobenius property means 
that 
$$A\cong A^{\ast}:=\Hom_R(A,R)$$ 
as an $A$-module. 
The isomorphism takes ${\bf 1}\in A$ to an $R$-linear trace map $\epsilon:A\lra R$ 
which is non-degenerate, meaning that the above map $a\longmapsto \epsilon(a\ast)$ from 
$A$ to $A^{\ast}$ is an isomorphism.  When $R$ is a field, $\epsilon$ is nondegenerate 
iff $\forall a\in A\setminus\{0\}\ \exists b$ such that $\epsilon(ab)\not= 0$.  Given 
$\epsilon$ as above, we can reconstruct $\Delta$ as the dual of $m$: 
$$ \Delta \ : \ A \cong A^{\ast}\stackrel{m^{\ast}}{\lra} 
A^{\ast}\otimes A^{\ast} \cong A \otimes A.$$
\begin{example}\label{ex-f-1} 
The direct sum of even-dimensional 
cohomology groups $H^{\text{even}}(M,R)$ of a closed oriented $2n$-dimensional manifold $M$ 
is a commutative Frobenius $R$-algebra, with the trace map given by the
integration over the fundamental $2n$-cycle.  
\end{example}
\begin{example}\label{ex-f-2}
 Let $R$ be a field and $f\in \C[x_1, \dots, x_m]$ a polynomial. 
If the quotient algebra $A$ of $\C[x_1, \dots, x_m]$ by the ideal generated by 
all partial derivatives 
$\frac{\partial f}{\partial x_1}, \dots, \frac{\partial f}{\partial x_m}$ 
is finite-dimensional then $A$ is Frobenius. This example 
comes up in singularity theory, see~\cite{AVGZ}. 
\end{example}
The above two types of Frobenius algebras have a nonempty intersection. 
For instance, if $f= x^{n+1}\in \Q[x]$ in the second example then 
$A= \Q[x]/(x^n)$, isomorphic to the cohomology ring of the complex 
projective space $\mathbb{CP}^{n-1}$. Notice that we only get 
a countable number of commutative Frobenius algebras (up to isomorphism) 
from Example~\ref{ex-f-1}, but an uncountable number from Example~\ref{ex-f-2}. 
Both examples are important sources of 2-dimensional TQFT's. 


\subsection{Algebras $H^n$}
Our goal in this lecture and the next one is to extend link homology to 
tangles and tangle cobordisms. We start with an arbitrary Frobenius algebra $A$ 
over $R$ and construct an invariant of flat (or crossingless) tangles. 

Consider $2n$ points on a horizontal line and denote by $B^n$ the set 
of crossingless matchings of these points by $n$ arcs lying in 
the lower half-plane. The cardinality of $B^n$ is the $n$-th Catalan number 
\begin{math}\frac{1}{n+1}\left( \begin{smallmatrix} 2n \\ n \end{smallmatrix} \right) .  
\end{math}
\begin{example} 
\label{crossingless}
The set $B^3$ has 5 elements 
\vspace{0.2in} 
\begin{center} \raisebox{-0.3cm}{\psfig{figure=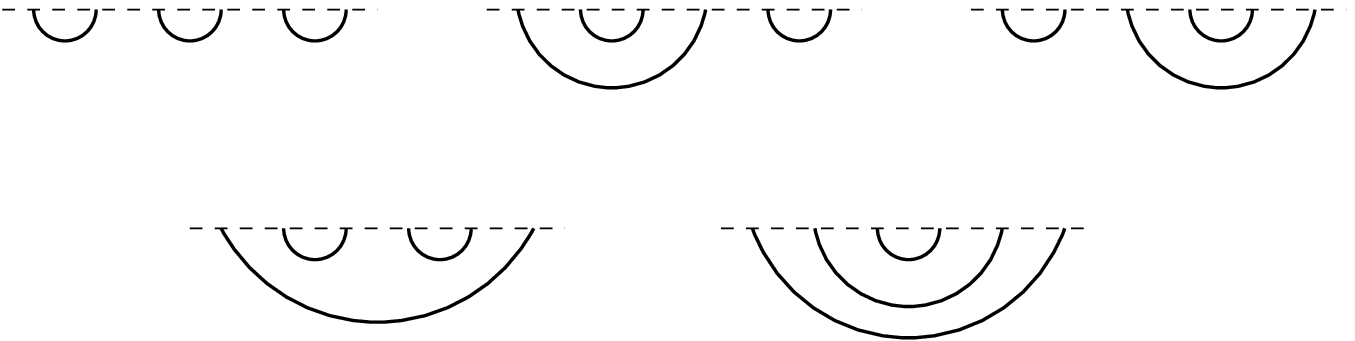,height=2.8cm}}\end{center}
\end{example}
Let $W(b)$ be the reflection of a matching $b$ along the horizontal line. For $a, b \in B^n$, 
the composition $W(b)a$ makes sense and can be viewed as a closed 1-manifold. 
\begin{center}\psfrag{C6-1}{{$b$}}\psfrag{C6-2}{{$W(b)$}}
\psfrag{C6-3}{{$a$}}\psfrag{C6-4}{{$W(b)a$}}
{\psfig{figure=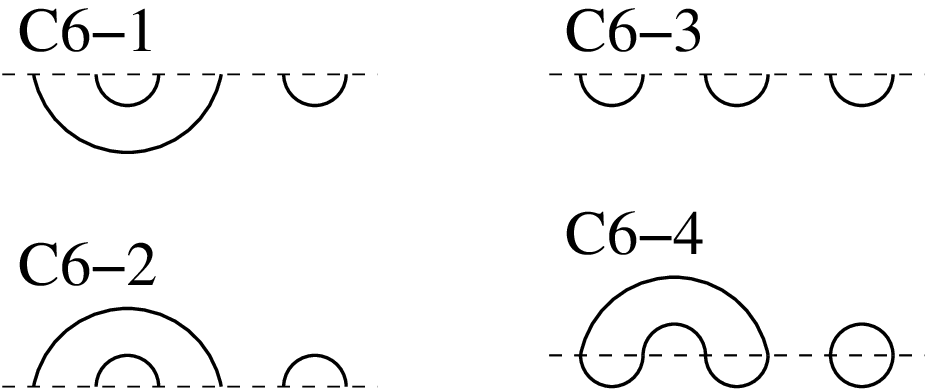,height=2.8cm}}\end{center}
Applying the functor $F$ to it, we get $F(W(b)a)$, which is a tensor power of $A$. 

For each $n\ge 0$ we define the ring $H^n$ by 
$$ H^n \ := \ \oplusop{a,b\in B^n} F(W(b)a).$$ 
The multiplication in $H^n$ is built out of compositions 
$$ F(W(c) b) \otimes F(W(b)a) \lra F(W(c)a)$$ 
induced by cobordisms from $W(c)bW(b)a$ to $W(c)a$ which contract $b$ with $W(b)$: 
$$
\begin{array}{ccc}
 \raisebox{-1cm}{\psfig{figure=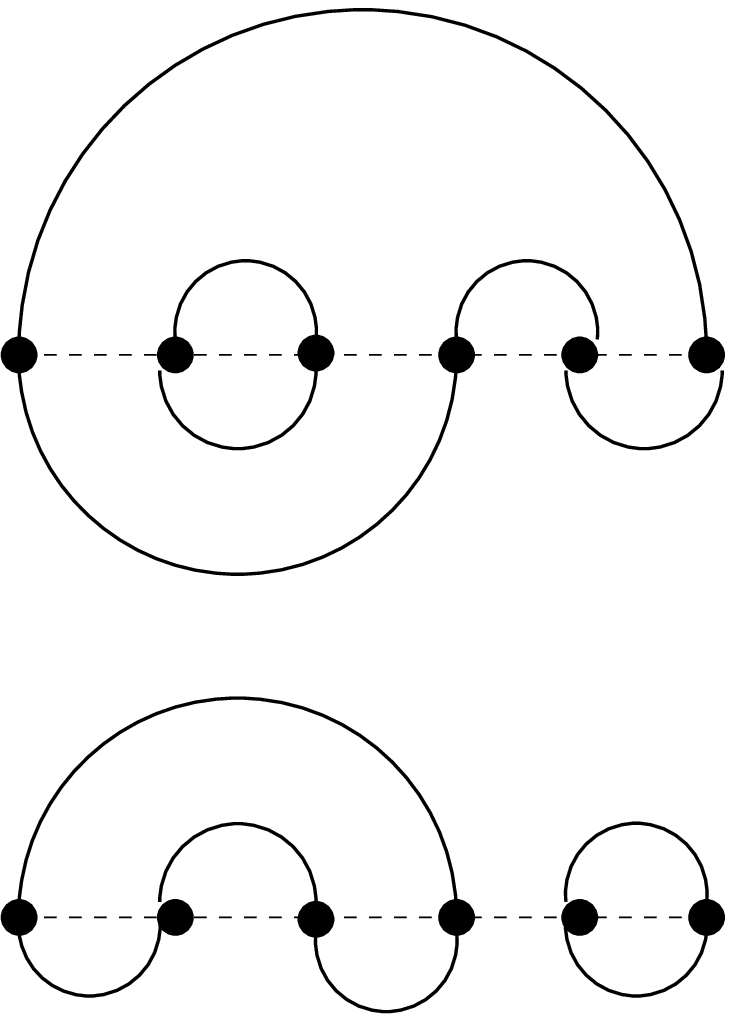,height=3cm}} 
\begin{array}{c}
W(c) \\ b \\ \    \   \\ W(b) \\ a
\end{array}
&
{\stackrel{\psfig{figure=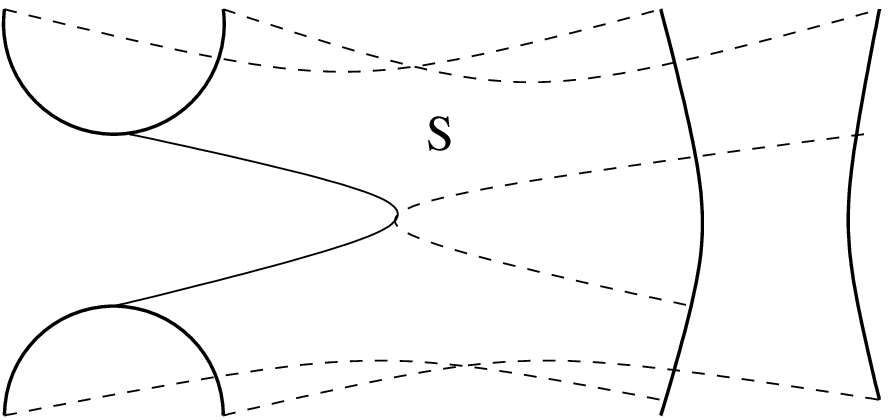,height=1cm}}{\longrightarrow}}
&
 \raisebox{-1cm}{\psfig{figure=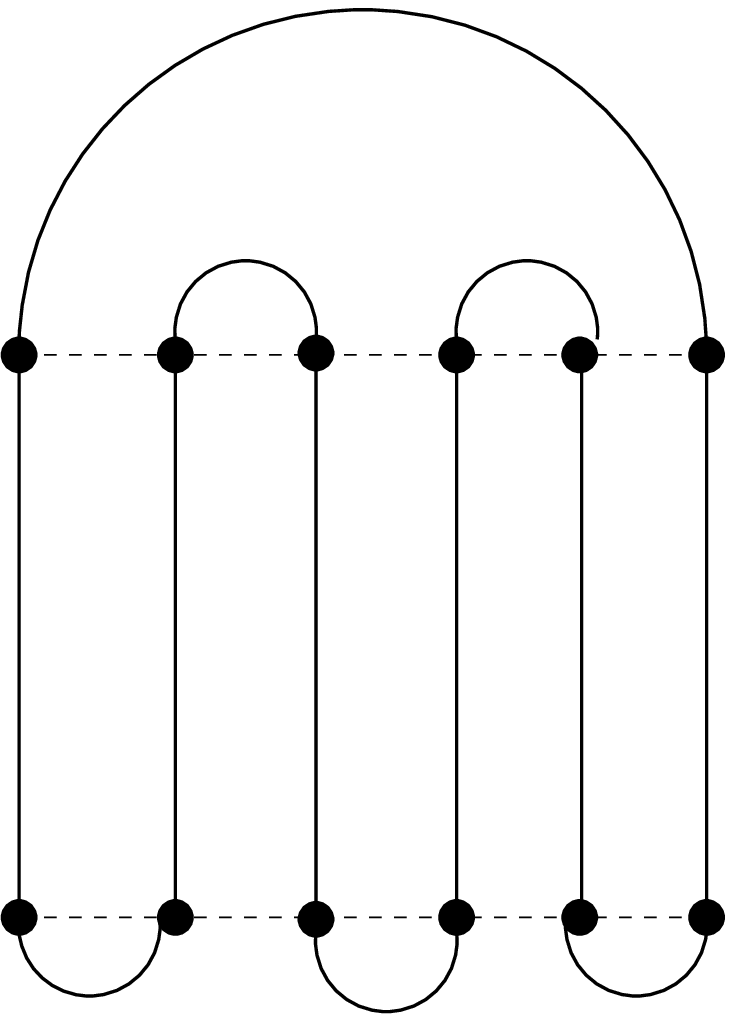,height=3cm}}  W(c)a \\
 & \Downarrow F & \\
F(W(c)b) \otimes F(W(b)a) & \longrightarrow & F(W(c)a) 
 \\
 \cap & & \cap \\
 H^n \otimes H^n & \stackrel{m}{\longrightarrow} & H^n \end{array}
 $$
For $x  \in F(W(d)c)$ and $y \in F(W(b)a)$ the product $xy=0$ if 
$c \neq b$. For $n=0$ the set $B^0$ contains only the empty diagram, and 
$H^0=R,$ the ground ring. 
 \begin{example} The set $B^1$ consists of the single diagram 
\{ \raisebox{-0.1cm}{\psfig{figure=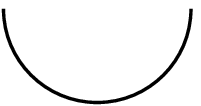,height=0.3cm}}\} which we denote $a$.  
 $H^1=F(W(a)a)=A.$ The product is given by $F(S): A^{\otimes 2} \to A$, where 
 \psfrag{waa}{$ {W(a)a}$}
 \psfrag{S}{$S$}
$$  \raisebox{-1cm}{\psfig{figure=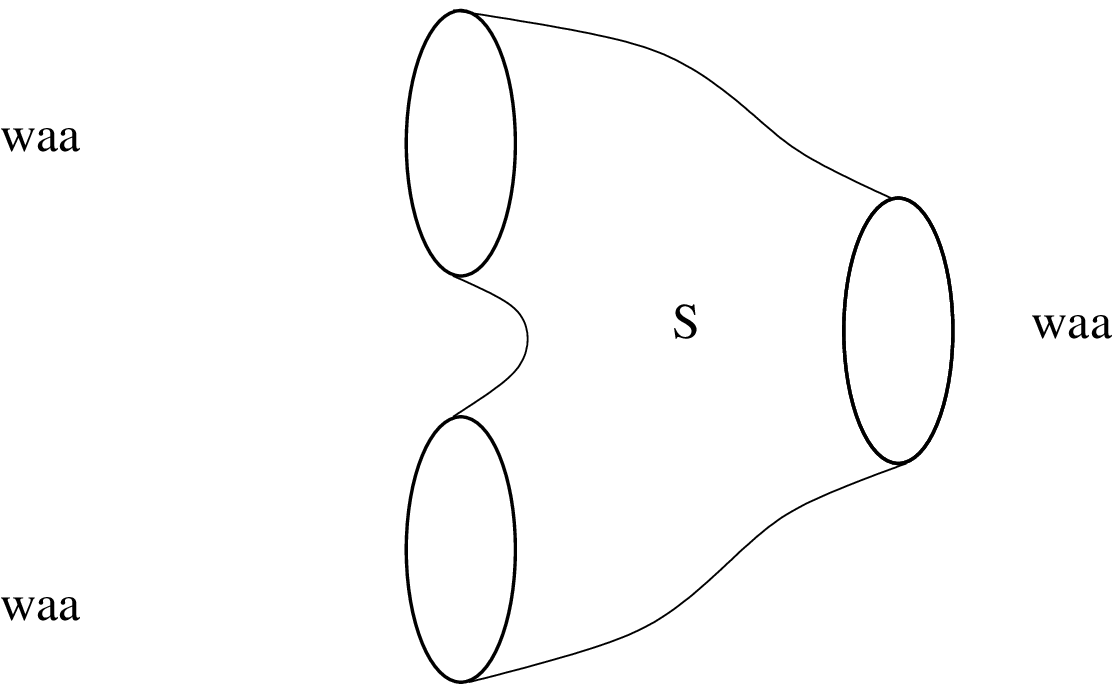,height=2cm}} $$ 
 Thus, the product in $H^1$ is the multiplication in $A$ and $H^1\cong A$ 
as an associative $R$-algebra. 
 \end{example}
\begin{example}
For $n=2$ the set $B^2$ has two elements 
$a= \raisebox{-0.1cm}{\psfig{figure=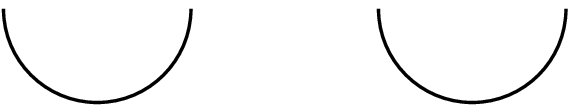,height=0.3cm}}$ and 
$b= \raisebox{-0.3cm}{\psfig{figure=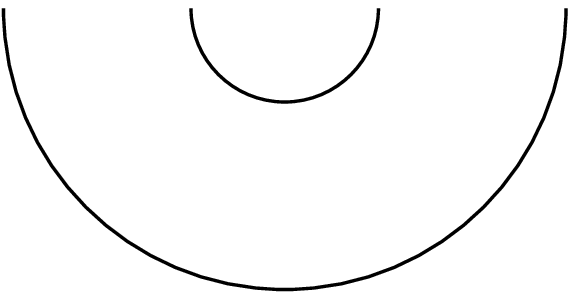,height=0.6cm}}$. 
$$
\begin{array}{cccccccc}
H^2=& F(W(a)a) &\oplus &F(W(b)a) &\oplus& F(W(a)b) &\oplus & F(W(b)b) \\
&  \raisebox{-0.1cm}{\psfig{figure=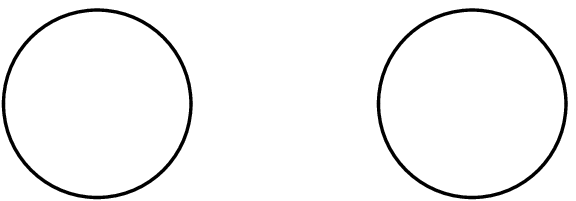,height=0.4cm}} \  A^{\otimes 2} &&  
\raisebox{-0.2cm}{\psfig{figure=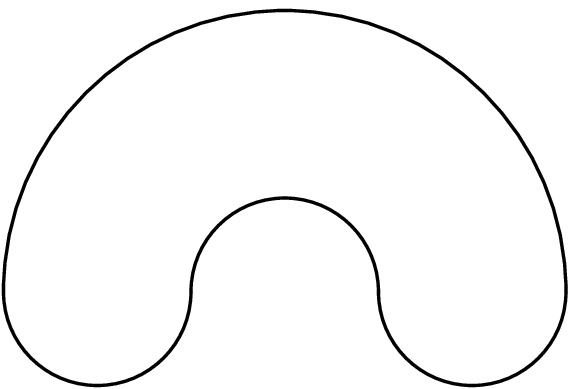,height=0.5cm}} \ A &&  \raisebox{-0.2cm}
{\psfig{figure=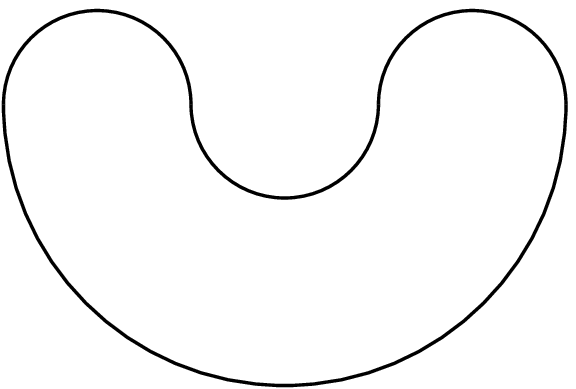,height=0.5cm}}  \ A &&
 \raisebox{-0.2cm}{\psfig{figure=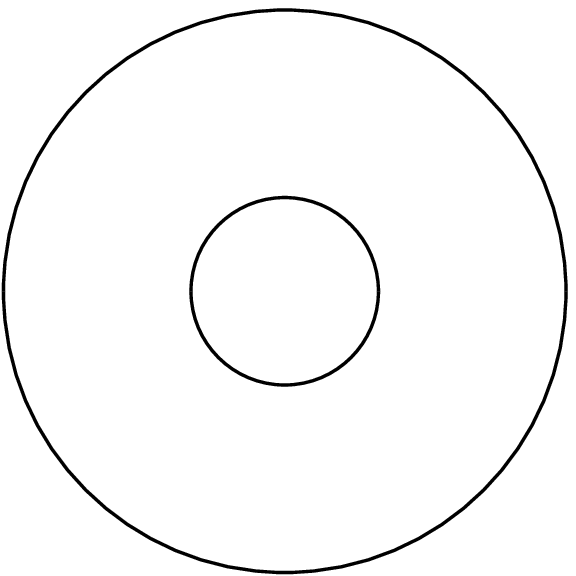,height=0.5cm}} \  A^{\otimes 2}
 \end{array}
 $$
 An example of multiplication
 $$ F(W(a)b) \otimes F(W(b)a) \lra F(W(a)a)$$ is the composition of 
morphisms: 
 \psfrag{Wa}{$W(a)$}
 \psfrag{b}{$b$}
 \psfrag{Wb}{$W(b)$}
 \psfrag{a}{$a$}
 $$
\vspace{0.06in}
 \begin{array}{ccccc}
  \raisebox{-1cm}{\psfig{figure=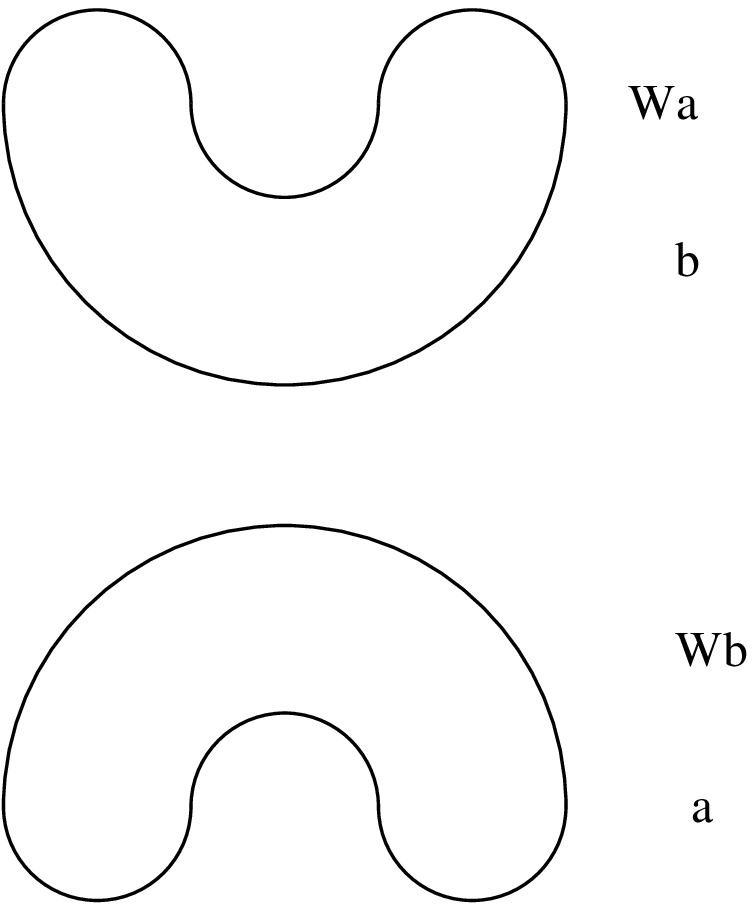,height=2cm}}   & \longrightarrow & 
   \raisebox{-1cm}{\psfig{figure=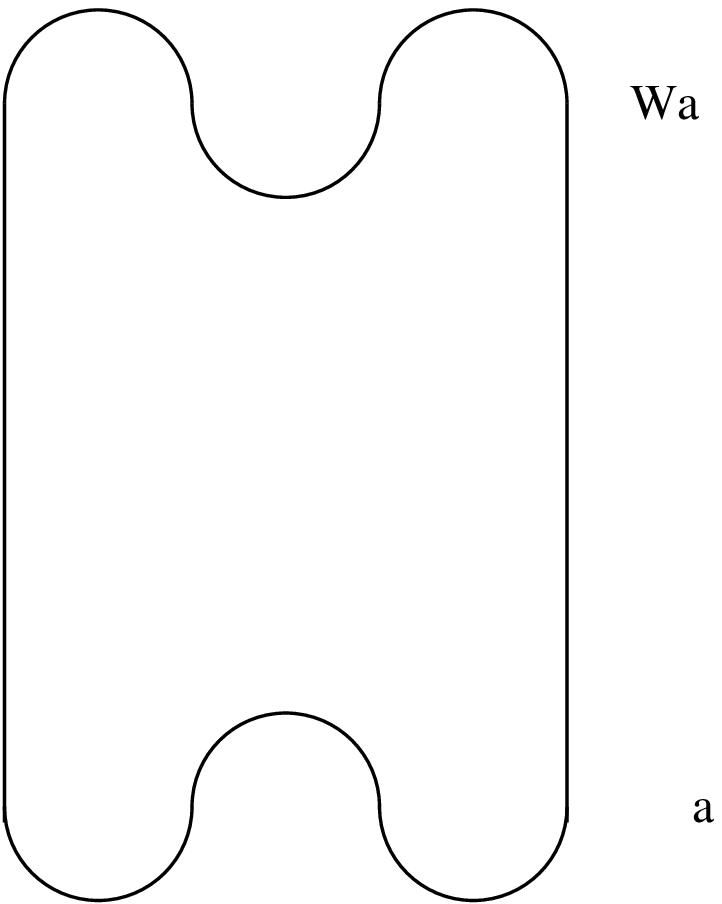,height=2cm}}  & \longrightarrow &
    \raisebox{-1cm}{\psfig{figure=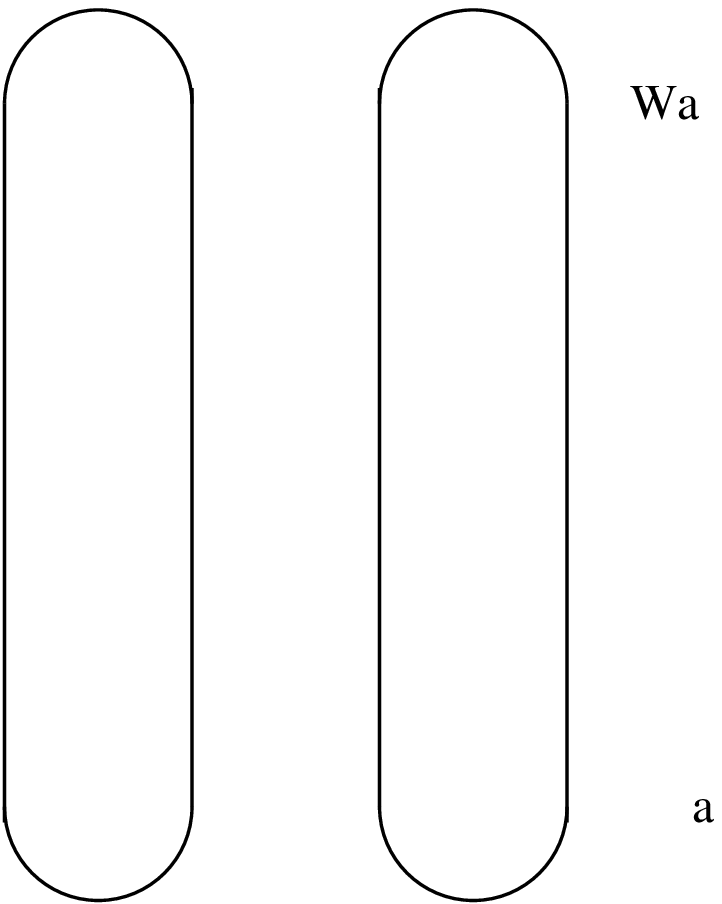,height=2cm}} 
  \\
  A^{\otimes 2} & \stackrel{m}{\longrightarrow} & A & \stackrel{\Delta}{\longrightarrow} 
&A^{\otimes 2} 
  \end{array}
  $$
\end{example}
\begin{exercise} Functor $F$ is defined on oriented cobordisms only. Find 
a consistent way to equip 1-manifolds $W(b)a$ and the multiplication cobordisms 
with orientations to make legitimate the above application of functor $F$. 
\end{exercise} 
\begin{exercise} Use the functoriality of $F$ to show that the multiplication in 
$H^n$ is associative. 
\end{exercise} 
\begin{exercise} Define $1_a\in H^n$ as ${\bf 1}^{\otimes n}\in A^{\otimes n} \cong 
F(W(a)a).$ Show that $x 1_a =x$ for any $x\in F(W(b)a)$ and 
$1_a y = y$ for any $y\in F(W(a)b).$ Check that $\{1_a\}_{a\in B^n}$ 
are mutually orthogonal idempotents and 
$$1= \sum_{a\in B^n} 1_a$$ 
is the unit element of $H^n$. 
\end{exercise} 
  \begin{remark} \label{rem-pr} 
Observe the similarity with the setting of the ring $A_n$ from Lecture 1, 
with idempotents $1_a$ of $H^n$ analogous to idempotents $(i)$ of $A_n$. 
The latter were used to define projective $A_n$-modules $P_i =A_n (i)$.  
Likewise, any $a \in B^n$ produces a  
left projective $H^n$-module
$$P_a:= H^n 1_a=\oplusop{b \in B^n} F(W(b)a)$$
and a right projective $H^n$-module
$${}_a P := 1_a H^n = \oplusop{b \in B^n} F(W(a)b).$$
\end{remark}

To summarize, $H^n$ is a unital associative $R$-algebra built out of a commutative 
Frobenius $R$-algebra $A$. For $n>1$ the algebra $H^n$ is noncommutative.  


\subsection{Flat tangles and their cobordisms}
By a \emph{flat} or \emph{crossingles} tangle we mean a finite 
collection of arc and circles properly embedded in $\R\times [0,1]$. 
 \psfrag{T-1}{$1$} \psfrag{T-2}{$2n$} \psfrag{T-3}{$1$}\psfrag{T-4}{$2$}
 \psfrag{T-5}{$2m$}
$$  \raisebox{-0.4cm}{\psfig{figure=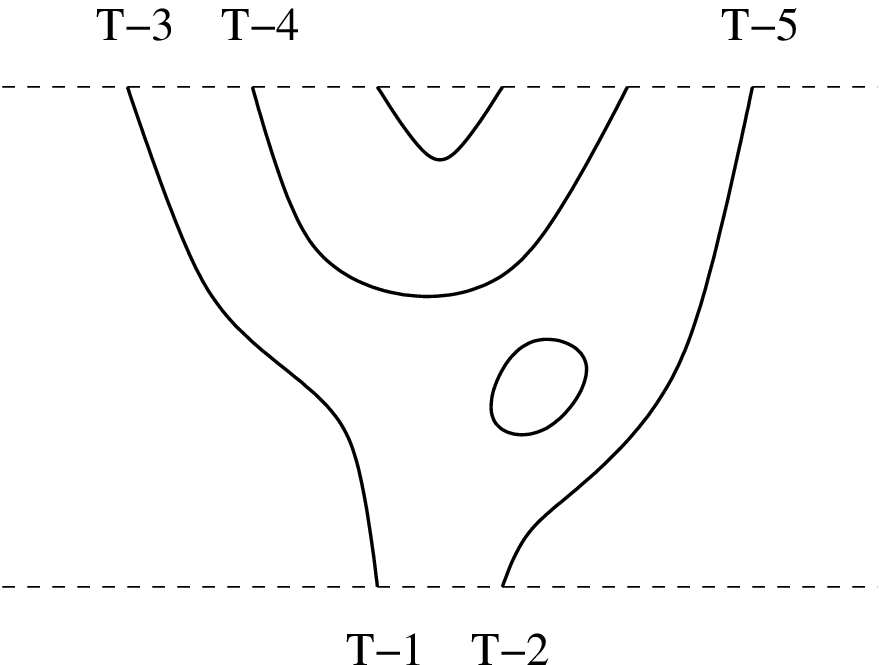,height=3.5cm}} $$ 
We require that the number of top endpoints be even, which implies that the 
number of bottom endpoints is even as well, and call a flat tangle with $2m$ 
top and $2n$ bottom endpoints a flat $(m,n)$-tangle. We also fix once and for 
all the position of $2n$ points on $\R$, to make flat tangles easy to 
compose. The composition of a flat $(k,m)$-tangle and a flat $(m,n)$-tangle is a 
flat $(k,n)$-tangle. 

By a cobordism $S$ between flat $(m,n)$-tangles $T, T'$ we mean a surface properly embedded in 
$\R\times [0,1]\times [0,1]$ with the boundary comprised of $T, T'$ and the product 
1-manifold $\partial T \times [0,1]\cong \partial T'\times [0,1].$ We think of 
$\partial T\times \{0\}$ as the boundary of $T$ and $\partial T' \times \{1\}$ as the boundary  
of $T'$. 
\psfrag{T-6}{$T$} \psfrag{T-7}{$T'$} \psfrag{T-8}{$S$}
$$  \raisebox{-0.3cm}{\psfig{figure=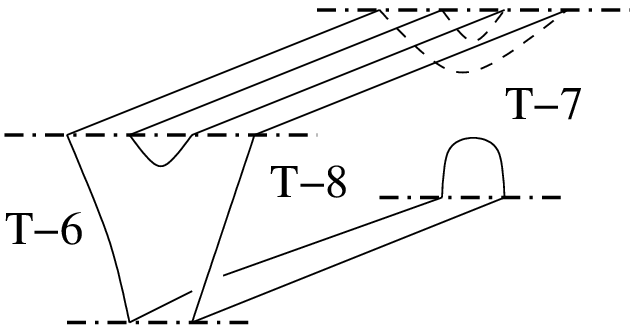,height=4cm}} $$ 
The upper part of $T'$ is shown by dashed lines; the corner $\R$'s of the 
3-manifold $\R\times [0,1]\times [0,1]$ are indicated by dashed-dotted lines. 
Surface $S$ has $4n+4m$ corner points. 

There are two ways to compose these cobordisms. If $S_1$ is a cobordism 
from $T$ to $T'$ and $S_2$ a cobordism from $T'$ to $T''$, we can glue 
them along $T'$ to produce the cobordism $S_2 \circ S_1$ from $T$ to $T''$. 
If $S_1$ is a cobordism between flat $(m,n)$-tangles $T_1$ and $T_1'$ and $S_2$ a cobordism 
between flat $(k,m)$-tangles $T_2$ and $T_2'$, we can compose $S_2$ and $S_1$ 
along the one-manifold $\{\text{$2m$ points}\}\times [0,1]$ to get the 
cobordism $S_2 S_1$ from $T_2 T_1$ to $T_2' T_1'$. 

This structure can be encoded into the 2-category of flat tangle cobordisms. 
The objects of this 2-category are nonnegative integers $n$, one-morphisms from 
$n$ to $m$ are flat $(m,n)$-tangles $T$, two-morphisms from $T$ to $T'$ are 
isotopy classes rel boundary of flat tangle cobordisms $S$. 

\vspace{0.1in} 

To a given flat $(m,n)$-tangle  $T$ we assign the $R$-module 
$$ F(T) \ := \ \oplusop{a\in B^n, b\in B^m} F(W(b)Ta). $$ 
In other words, we consider all possible ways to close up $T$ by crossingless 
matchings $a$ and $b$ at the bottom and the top, respectively, 
 to produce a closed 1-manifold $W(b)Ta$, and apply the functor $F$ to each closure. 
\psfrag{wb}{$W(b)$}
\psfrag{a}{$a$}
\psfrag{T}{$T$}
$$
\psfig{figure=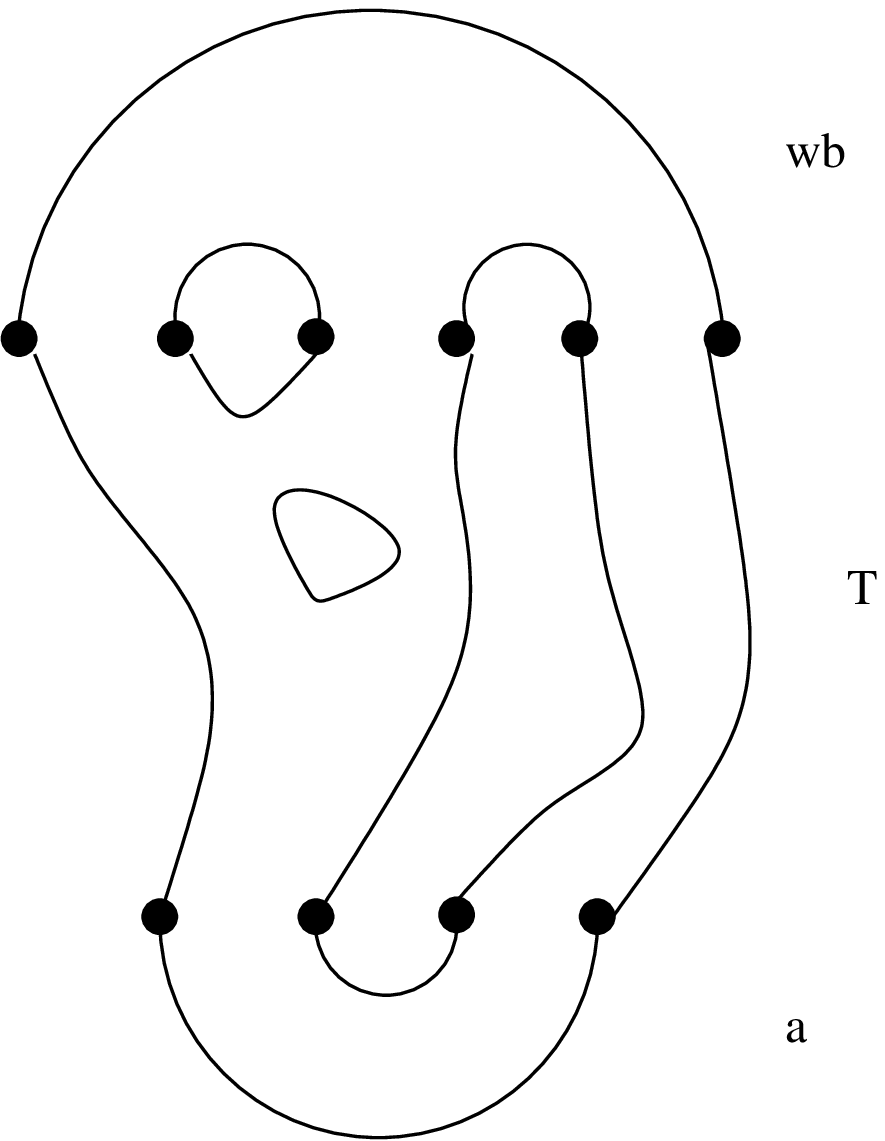,height=2.5cm}
$$
$F(T)$ is actually an $(H^m, H^n)$-bimodule, that is, it has a right $H^m$-action and 
a commuting left $H^n$-action. In the rest of the notes, we call an $(H^m, H^n)$-bimodule
simply an $(m,n)$-bimodule, and assume that in these bimodules the left and the right 
action of $R$ are equal. The action of $H^m$ on $F(T)$ comes from maps 
$$
\begin{array}{ccccc} 
 F(W(c)b) & \times & F(W(b)Ta) & \lra & F(W(c)Ta) \\
  \cap       &            &   \cap         &        &   \cap \\
   H^m    &  \times &   F(T)        & \lra &   F(T) 
\end{array}
$$
\begin{example} $F$ applied to the identity flat $(n,n)$-tangle produces $H^n$
viewed as an $H^n$--bimodule: 
$$F\left(\raisebox{-0.4cm}{\psfig{figure=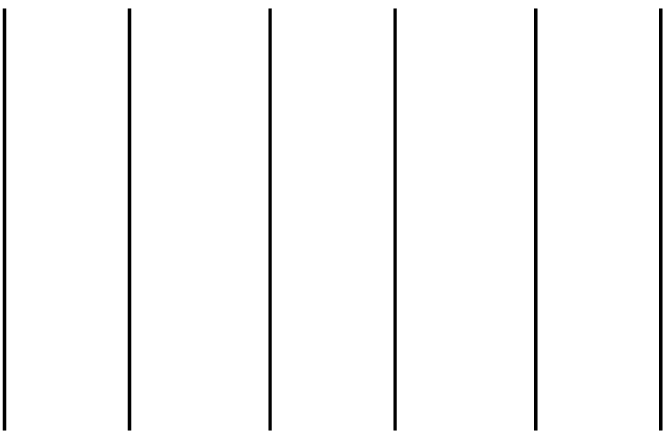,height=1cm}}\right)  \cong H^n.$$
\end{example}
\begin{example}
A crossingless matching $a\in B^m$ is a flat $(m,0)$-tangle and the bimodule 
$F(a)$ is simply the left $H^m$-module $P_a$, see Remark~\ref{rem-pr}
 (notice that $(m, 0)$-bimodules are just left $H^m$-modules, since $H^0= R$, 
the ground ring). Likewise, $F(W(a))\cong {}_aP$ is a right projective $H^m$-module.  
\end{example}
\begin{example} 
A flat $(0,0)$-tangle $T$ is a closed 1-manifold embedded in the plane, 
$(0,0)$-bimodules are just $R$-modules, and  
$F(T)=A^{\otimes r}$, where $r$ is the number of components of $T$. 
\end{example}
The composition of flat tangles $T_2 T_1$ corresponds to the tensor product of 
bimodules 
$$ F(T_2 T_1) \cong F(T_2) \otimes_{H_m} F(T_1) $$ 
for a flat $(k,m)$-tangle $T_2$ and a flat $(m,n)$-tangle $T_1$, 
see~\cite[Theorem 1]{FVIT}. 

If we fix $n$ and consider only flat $(n,n)$-tangles $T$ and $H^n$-bimodules 
$F(T)$ we get a functor realization of the Temperley-Lieb algebra $TL_{2n}$. 
In Lecture 1 we already constructed a realization of $TL_{n+1}$ by $A_n$-bimodules, 
with the generators
$$u_i=\raisebox{-1.3cm}{\psfig{figure=lec1.14.eps,height=2.5cm}}$$
of the Temperley-Lieb algebra represented by $A_n$-bimodules $U_i$. 
Recall that $U_i \otimes U_j=0$ if $|i-j| > 1$, while 
$u_i u_j \neq 0$ in $TL$ algebra, so our bimodule realization was degenerate. 
On the other hand, $F(T)$ is a non-trivial bimodule for 
any flat $(n,n)$-tangle $T$. Also, in the current setting a closed loop evaluates 
to the $R$-module $A$. For instance, $F(u_i u_i) \cong F(u_i)\otimes_R A.$ 
When we pass to ranks, the value of the closed loop becomes the rank of 
$A$ as a free $R$-module, a positive integer. 

To get more general values for the closed loop we extend the 
framework of commutative Frobenius $R$-algebras $A$ and rings $H^n$ to 
the graded case, by requiring that $A$ be a graded $R$-algebra and 
the morphism $F(S)$ associated with a 2-dimensional cobordism $S$ be 
homogeneous of degree proportional to the Euler characteristic of $S$. 
Under these assumptions, the rings $H^n$ and bimodules $F(T)$ become 
graded~\cite{FVIT}. A closed loop still corresponds to $A$. Upon decategorification, 
closed loop evaluates to the graded rank of $A$ as 
$R$-module and takes value in $\mathbb{N}[q,q^{-1}]$. In the simplest nontrivial 
case of the graded pair $(R,A)$ described in Lecture 3 the value of the closed loop is 
$q+q^{-1}$, the standard value of the loop in the Temperley-Lieb algebra. 

Let $S$ be a cobordism in $\R \times [0,1] \times [0,1]$ 
between flat $(m,n)$-tangles $T_1$ and $T_2$. To $S$ we assign a homomorphism 
$$ F(S) \ : \ F(T_1) \lra F(T_2) $$ 
between $(m, n)$-bimodules. For $a\in B^n$ and $b\in B^m$ we 
can compose $S$ with the identity cobordism $\Id_a$ from $a$ to $a$ and the 
identity cobordism $\Id_{W(b)}$ from $W(b)$ to $W(b)$ to get 
the cobordism $\Id_{W(b)} S \Id_a$ from the closed 1-manifold 
$W(b)T_1 a$ to $W(b)T_2 a$. 
\psfrag{T1-1}{$T_1$}\psfrag{T1-2}{$T_2$}\psfrag{T1-3}{$S$} 
\psfrag{T1-4}{$\Id_{W(b)}$}\psfrag{T1-5}{$\Id_a$}\psfrag{T1-6}{$\Id_{W(b)}S\Id_a$}
$$\raisebox{-0.5cm}{\psfig{figure=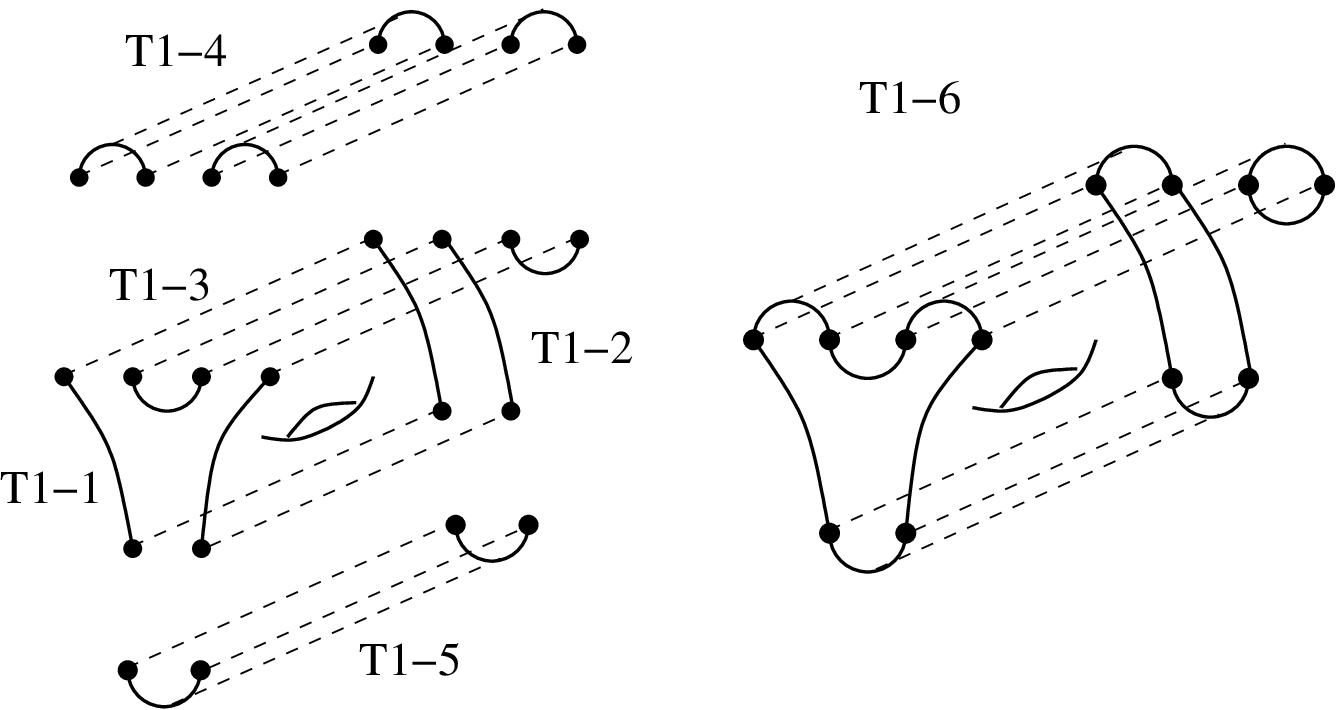,height=5cm}}$$
This cobordism induces a homomorphism $F(W(b)T_1 a) \lra F(W(b)T_2 a)$. 
Summing over all $a$ and $b$ we get a homomorphism of $(m,n)$-bimodules 
$F(S) : F(T_1)\lra F(T_2)$. 
It is straightforward to check that $F(S)$ is natural with respect to both types 
of compositions of cobordisms. When all the properties  are written down, we
 discover that $F$ becomes a 2-functor from the 
2-category of flat tangle cobordisms to the 2-category of 
$(m,n)$-bimodule homomorphisms. Objects of the latter 2-category
 are nonnegative integers $n$, 1-morphisms from $n$ to $m$ are $(m,n)$-bimodules and 
2-morphisms are bimodule homomorphisms. Composition of bimodules is 
given by the tensor product. The 2-functor $F$ is the identity on objects,
$n \longmapsto n$, takes a flat tangle $T$ to the bimodule $F(T)$,  
 and flat tangle cobordism $S$ to the 
bimodule homomorphism $F(S)$. 
This 2-functor converts topological information about 1-manifolds embedded in 
the plane and 2-manifolds embedded in $\R^3$ into the algebraic information 
provided by bimodules and bimodule homomorphisms. 

The following table summarizes the construction of the 2-functor $F$. 
$$
\begin{array}{llll}
& \begin{array}{c}{2\mbox{-category of flat}}\\  \mbox{tangle cobordisms} \end{array}& 
\stackrel{F}{\Longrightarrow} & 
\begin{array}{c}{2\mbox{-category of bimodule}}\\  \mbox{homomorphisms}\end{array}\\
 &    &    &    \\
\mbox{objects}& n=0,1,2, \dots & \longmapsto &n=0,1,2, \dots \\
1\mbox{-morphisms } & \mbox{flat $(m,n)$-tangles } T  & \longmapsto & (m, n)
\mbox{-bimodules  $F(T)$} \\
2\mbox{-morphisms}  & \mbox{flat tangle cobordisms $S$} &\longmapsto &  
\mbox{bimodule homomorphisms  $F(S)$}
\end{array}
$$


\section{A homological invariant of tangles and tangle cobordisms}


\subsection{An invariant of tangles}
 In the previous lecture we described a 2-functor from the 2-category of cobordisms 
between flat tangles to the 2-category of bimodule maps. Such functor exists for 
any Frobenius $R$-algebra $A$.  In this lecture, for a specific $(R,A)$, we extend 
the construction to a 2-functor  from the 2-category of tangle cobordisms to 
the 2-category of homomorphisms between complexes of bimodules (up to 
chain homotopy). Going one dimension up, from flat tangles, which are 2-dimensional 
objects, to tangles, which are 3-dimensional, exactly correponds to passing 
from the abelian category of bimodules to the triangulated category of 
complexes of bimodules. Ditto for cobordisms, which are 3-dimensional between 
flat tangles and 4-dimensional between tangles. Thus, in the framework described 
here, the key transformation in 
algebra 
\begin{center} abelian categories \ \ \  $\Longrightarrow$ \ \  \  triangulated categories
\end{center} 
is mirrored in low-dimensional topology by the transformation 
\begin{center} 
 (2+1)-dimensional structures  \ \ \ $\Longrightarrow$ \  \ \ (3+1)-dimensional structures. 
\end{center}     

We specialize the construction of the previous lecture to $R=\Z$ and 
$A=\Z[X]/(X^2)$ with the trace $\epsilon(X)=1,$ $\epsilon({\bf 1})=0$. 
The ring $H^n$ is made graded by defining 
$H^n=\oplusop{a,b \in B^n} F(W(b)a)\{n\}$. 
The multiplication in $H^n$ is grading-preserving and $\deg(1_a)=0$. 
 For a flat $(m,n)$-tangle $T$ the  $(m, n)$-bimodule $F(T)$ is graded. 

Start with an oriented $(m,n)$-tangle $T$. Just as in the flat 
case, we assume that $T\subset \R^2\times [0,1]$ has $2n$ bottom endpoints 
placed in the standard position on the plane $\R^2\times \{0\}$ and 
$2m$ top endpoints in standard position on  $\R^2\times \{1\}$. 
Oriented tangles can be composed in the same way as flat tangles, assuming 
that the orientations at the endpoints match. A tangle cobordism $S$ between 
$(m,n)$-tangles $T_1, T_2$ is an oriented smooth surface in 
$\R^2\times [0,1]^2$ subject to the boundary conditions mirroring 
those for flat tangles. In particular, the boundary of $S$ consists of 
four pieces, two of which are $T_1, T_2$ and the other two are product 
1-manifolds. Let $CobT$ be the 2-category of tangle cobordisms. Its 
objects are finite sequences of pluses and minuses, its 1-morphisms 
are tangles with prescribed orientations at the endpoints, and 
2-morphisms are tangle cobordisms up to rel boundary isotopies. 

We choose a diagram $D$ of a tangle $T$, a generic projection of $T$ onto 
the plane, with the endpoints projecting in the standard way. 
Assume that $D$ has a single crossing. Let $D^0$, $D^1$ be unoriented 
flat tangles obtained by resolving the crossing of $D$. 
\psfrag{T2-1}{$D$}\psfrag{T2-2}{$D^0$}
\psfrag{T2-3}{$D^1$}\psfrag{T2-4}{Cobordism $S$} 
$$\raisebox{-1.5cm}{\psfig{figure=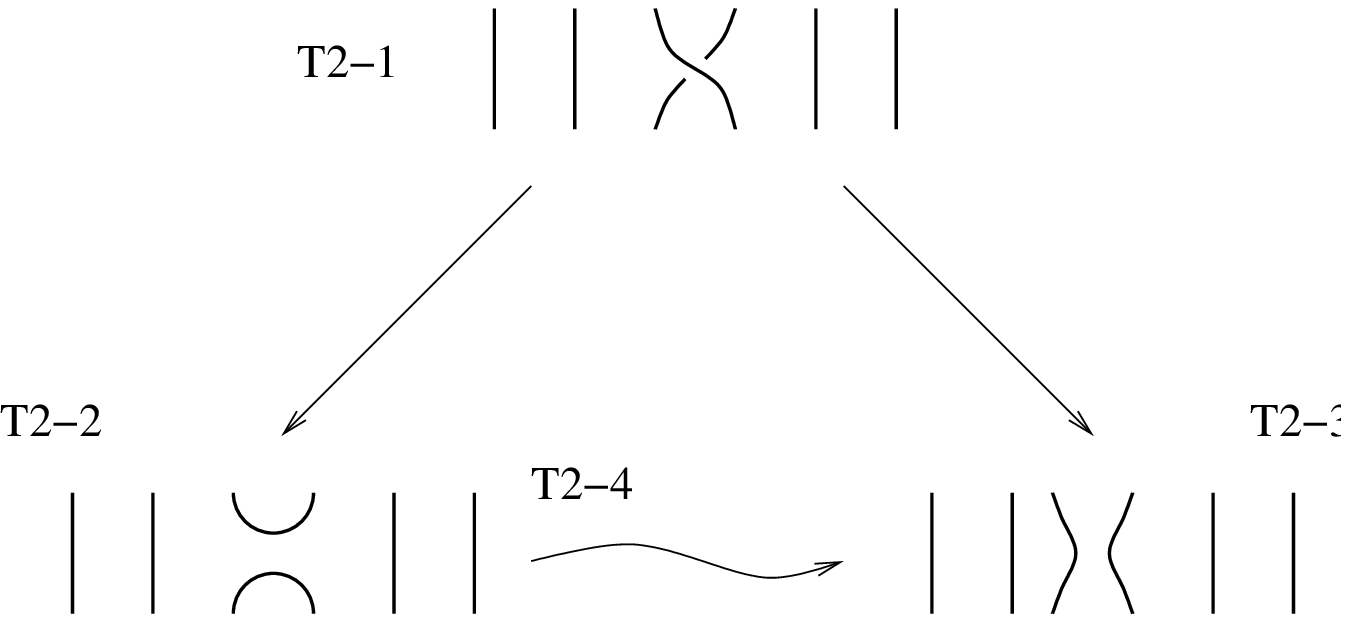,height=5cm}}$$
Let $S$ be the simplest cobordism between $D^0$ and $D^1$; it has one saddle point
for the projection $S\subset \R\times[0,1]^2\lra [0,1]$. Cobordism $S$ 
induces a degree $1$ morphism $F(S):F(D^0)\lra F(D^1)$ between graded 
$(m,n)$-bimodules. We define $F'(D)$ as the 
chain complex 
$$F(T):=(0 \to F(D^0) \stackrel{F(S)}{\to} F(D^1)\{-1\} \to 0)$$
of graded bimodules with a grading-preserving differential, with 
$F(D^0) $ in homological degree 0. Recalling that $D$ is a diagram of an 
oriented tangle $T$, we set 
$$ F(D) := F'(D)[x(D)]\{2x(D)-y(D)\}, $$ 
where, as before, $x(D)$ and $y(D)$ is the number of negative and positive crossings of 
$D$. 

If $D$ is an arbitrary tangle diagram, decompose $D$ as the composition of diagrams 
with at most one crossing each, $D=D_k \cdots D_2 D_1$
\psfrag{T3-1}{$D_1$}
\psfrag{T3-2}{$D_2$}
\psfrag{T3-3}{$D_k$}
$$\raisebox{-1cm}{\psfig{figure=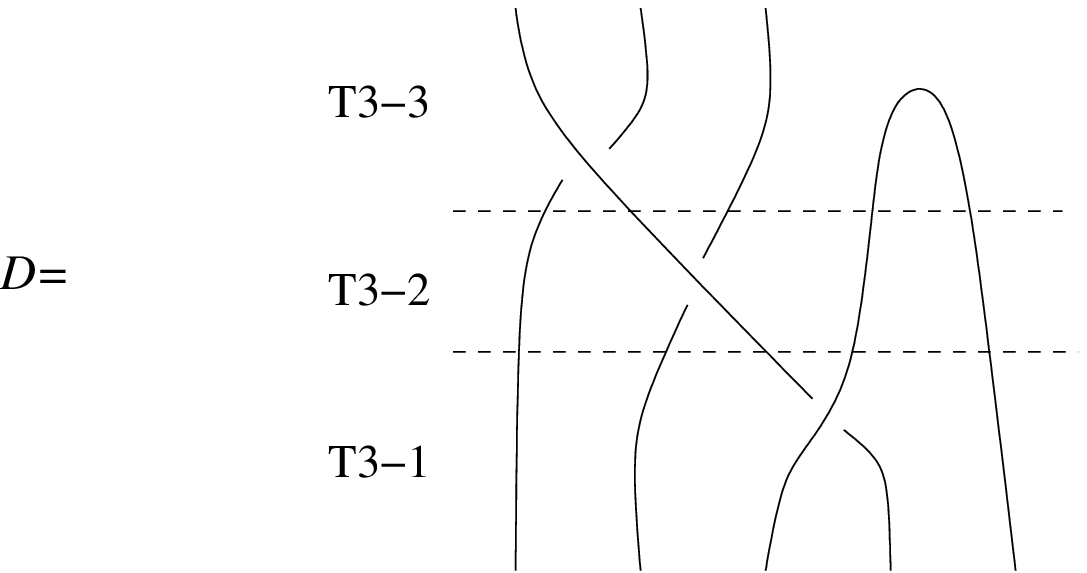,height=4cm}}$$
and define
$$F(D):=F(D_k)\otimes \cdots \otimes  F(D_2) \otimes F(D_1), $$
with the tensor product over rings $H^k$, for suitable $k$'s equal to half 
the number of top/bottom endpoints for intermediate diagrams $D_i$. 
This is a complex of graded $(m,n)$-bimodules. Let $\mc{C}_{m,n}$ be the 
category of complexes of graded $(m,n)$-bimodules up to chain homotopies. 

\begin{theorem} If diagrams $D_1, D_2$ of an $(m,n)$-tangle $T$ 
are related by a chain of Reidemeister moves, then $F(D_1) \cong F(D_2)$ in $\mc{C}_{m,n}$. 
\end{theorem}
The proof consists of constructing an explicit homotopy equivalence 
$F(D_1)\cong F(D_2)$ for two diagrams related by a Reidemeister move~\cite{FVIT}. 
\begin{example} $n=m=0$. In this case the tangle $T$ is a link, and the ring  
$H^0=\Z$. $F(D)\cong C(D)$ is a complex of graded abelian groups, and 
its homology groups $H(F(D))$ coincide with the link homology $H(T)$. 
\end{example}


\subsection{Tangle cobordisms}

Let $S$ be a tangle cobordism from tangle $T_0$ to tangle $T_1.$ 
The two tangles can be represented by their planar diagrams $D_0, D_1$. 
Likewise, we can represent $S$ by a sequence of planar diagrams of  
intersections of $S$ and $\R^2 \times [0,1] \times \{t\}$ for various $t \in [0,1]$ . 
Such representation is called a {\it movie} of $S$. Consecutive diagrams $D_0=D^0, D^1,
\dots, D^m=D_1$ in a movie 
differ by either a Reidemeister move or a move that corresponds to going through 
a critical point of the projection $S\lra [0,1]$. There are 3 types of critical point 
moves. The saddle move $ \raisebox{-0.2cm}{\psfig{figure=lec3.6.eps,height=0.6cm}} 
\leftrightarrow  \raisebox{-0.2cm}{\psfig{figure=lec3.7.eps,height=0.6cm}}$
correspond to passing through an index 1 critical point. Creation and 
annihilation moves $ \phi \leftrightarrow  \raisebox{-0.2cm}
{\psfig{figure=lec3.8.eps,height=0.6cm}}$ correspond to critical points of 
index 0 and 2, respectively. 
 An example of a movie is as follows.
  \psfrag{T0=T}{$D_0=D^0$}
  \psfrag{T1}{$D^1$}
  \psfrag{T2}{$D^2$}
  \psfrag{T3}{$D^3$}
  \psfrag{T4=T'}{$D^4=D_1$}
   $$S=\raisebox{-1.2cm}{\psfig{figure=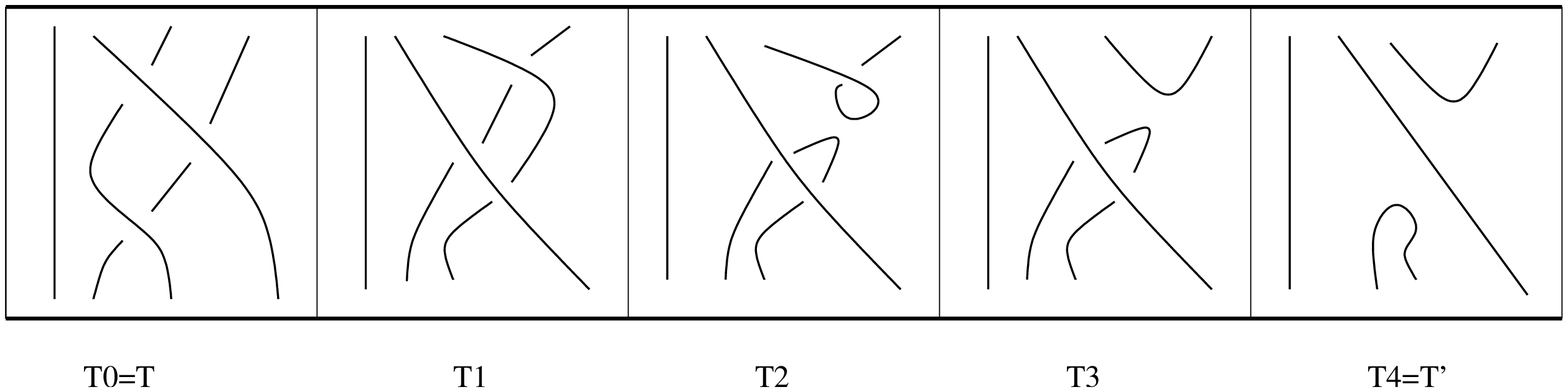,height=2.5cm}} $$
We often denote a movie representing a 
cobordism $S$ by $S$ as well. 
To a Reidemeister move between $D^i$ and $D^{i+1}$ we assign the 
isomorphism $F(D^i)\stackrel{\cong}{\lra} F(D^{i+1})$ mentioned 
earlier. To a degree $0$ critical point move from $D^i$ to $D^{i+1}$ we 
assign the map 
$$ F(D^i) \cong F(D^i) \otimes \Z \stackrel{\Id \otimes \iota}{\lra} 
F(D^i)\otimes A \cong F(D^{i+1})$$ 
induced by the unit map $\iota: \Z \lra A$ from $F(\emptyset)$ to 
$F(\text{circle})$. To a degree $2$ critical point move we assign the  
map induced by the trace homomorphism $\epsilon: A \lra \Z$. 
The map assigned to a degree $1$ critical point move was essentially 
described earlier in the lecture, as the bimodule homomorphism  induced 
by the standard cobordism between two resolutions of a crossing. 
For instance, in $S$ shown above,  diagrams $D^1$ and $D^2$ are related 
by such a move (saddle move). We can decompose each $D^1$ and $D^2$ 
as the composition of 3 diagrams, with only the middle diagrams being different 
\psfrag{Tt}{$$}
 \psfrag{Ttt}{$$}
  $$ \raisebox{-1.2cm}{\psfig{figure=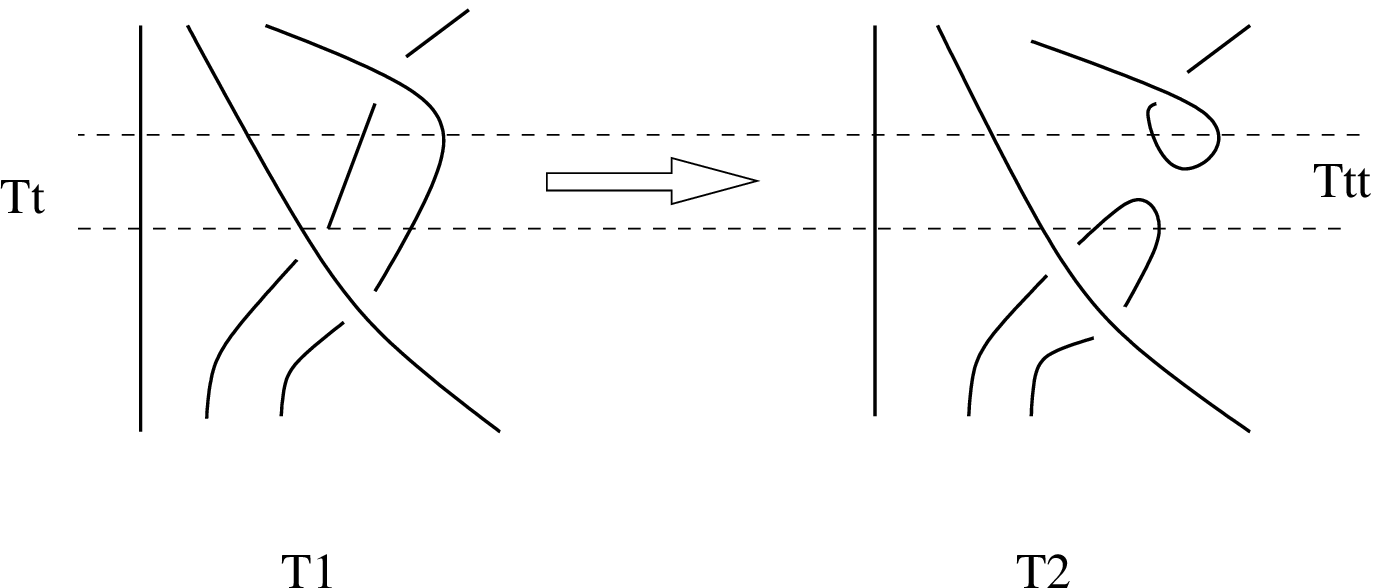,height=2cm}} $$
and define $F(D^1)\lra F(D^2)$ as the composition of the identity map 
on the first and the third terms and the homomorphism of the middle terms 
induced by the saddle point cobordism between two crossingless tangles. 

We define $F(S)$ as the composition of homomorphisms $F(D^i)\lra F(D^{i+1})$ 
associated to frame changes $D^i \lra D^{i+1}$. Notice that homomorphisms 
associated to Reidemeister moves are invertible in $\mc{C}_{m,n}$, 
the category of complexes of graded $(m,n)$-bimodules modulo chain homotopies. 
For instance, for the movie $S$ drawn above, the homomorphism 
$F(D^1) \lra F(D^2)$ is the only non-invertible one. 

It is a theorem of S.~Carter and M.~Saito~\cite{CS} that 
two movies represent the same tangle cobordism $S$ 
iff they can related by a sequence of {\it movie moves}. A movie move converts 
a certain sequence of frames to another sequence of frames representing 
the same cobordism. Here's an example: 
\psfrag{T4-1}{$D_0$} \psfrag{T4-2}{$D_1$}
\begin{eqnarray*}
 S_0&=& \raisebox{-1.2cm}{\psfig{figure=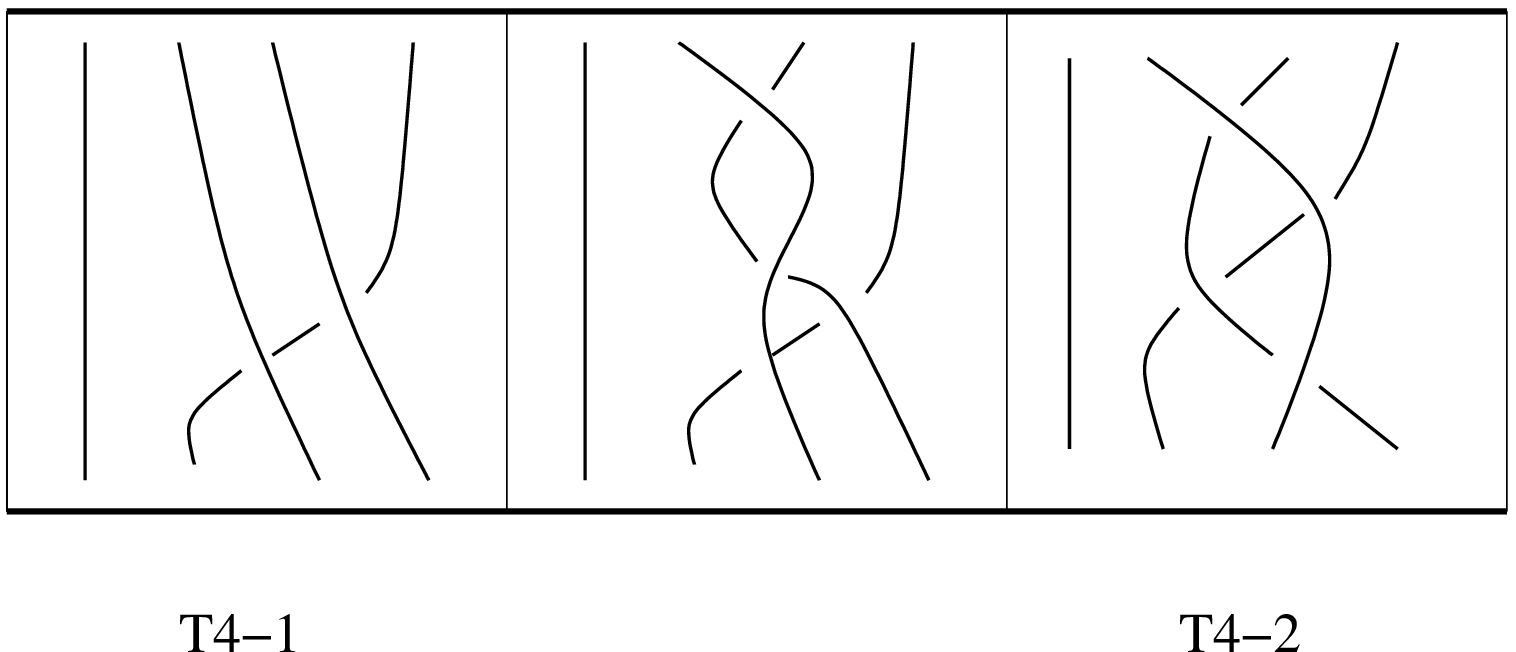,height=2cm}} \\
 S_1&=&  \raisebox{-1.2cm}{\psfig{figure=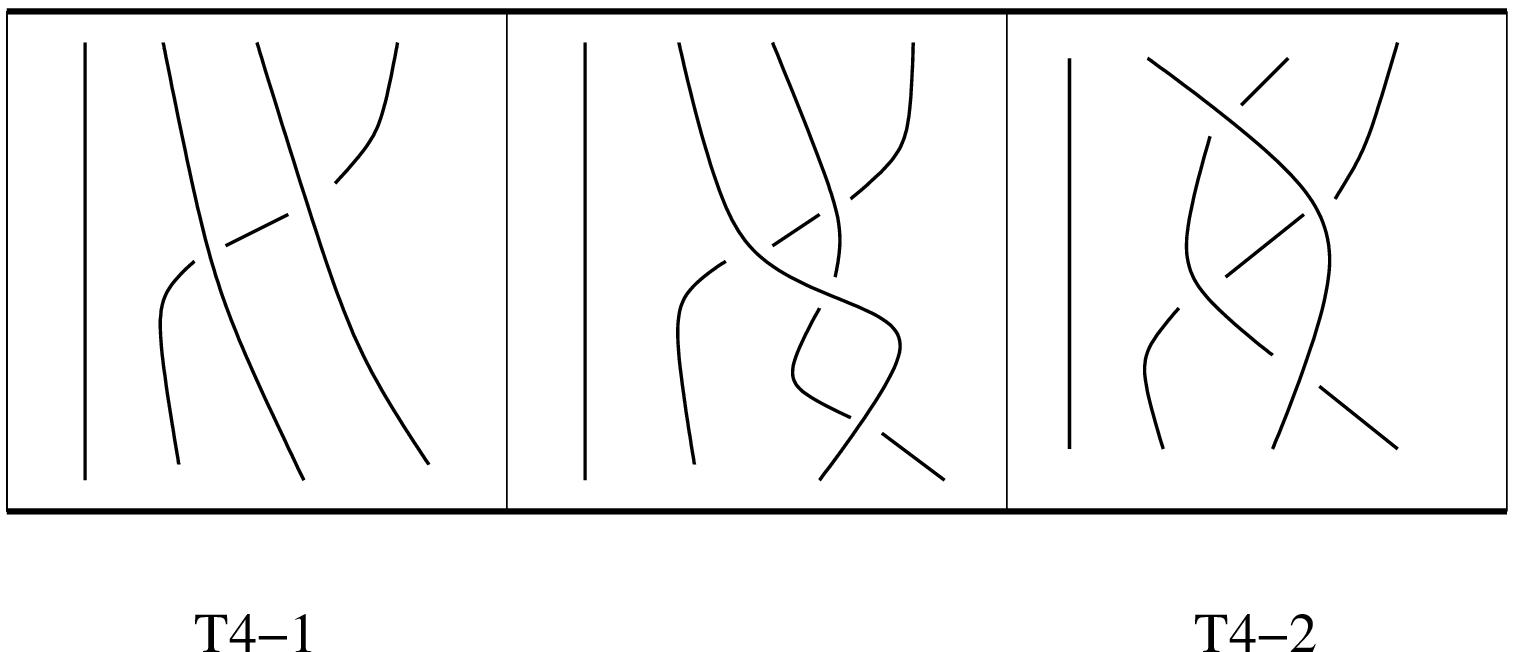,height=2cm}}
 \end{eqnarray*}
We refer the reader to~\cite{CS} or \cite{Kh2} for a complete list of movie moves. 

 \begin{theorem} \label{th-cobs-i}
 If movies $S_0$ and $S_1$ represent isotopic tangle cobordisms, 
$F(S_0)=\pm F(S_1)$ in $\mc{C}_{m,n}$.
 \end{theorem}
 Proving this statement amount to checking the invariance of $F$ up to a sign for each 
movie move. We only explain it for the above example of a movie move. 
We have two morphisms: $F(S_0), F(S_1): F(D_0) \to F(D_1)$. They are isomorphisms 
of complexes of bimodules in the homotopy category $\mc{C}_{n,n}$, since both 
movie moves $S_0, S_1$ are compositions of Reidemeister moves. 
Then $f=F(S_1)^{-1}F(S_0)$ is an isomorphism of $F(D_0)$. Note that 
$F(D_0)$ is an invertible complex of bimodules, since it is represented by a braid. 
Tensoring it with the complex associated with the inverse braid will give 
us the identity bimodule $H^n$. Therefore, the group of automorphisms of $F(D_0)$ 
in $\mc{C}_{n,n}$ is isomorphic to the group of automorphisms of $H^n.$ 
Automorphisms of the $H^n$-bimodule $H^n$ are multiplications by invertible 
central elements of $H^n.$ Moreover, automorphisms in the category $\mc{C}_{n,n}$ 
should preserve the internal grading of $H^n$. It is a simple exercise to check that 
there are only two such automorphisms, $\pm \Id$. Thus, $f= \pm \Id$ and 
$F(S_1) = \pm F(S_0)$. 

This argument works for any movie move where each movie is a sequence of 
Reidemeister moves. Taking care of other movie moves is only slightly more 
complicated~\cite{Kh2}.

Putting everything together, we get a (projective) $2$-functor 
from the 2-category of oriented tangle cobordisms $TC$ to the 2-category $\mc{C}$. 
The objects of the latter are nonnegative integers, 1-morphisms are 
complexes of graded $(m,n)$-bimodules, and 2-morphisms are homogeneous 
homomorphisms of complexes modulo chain homotopies. The word \emph{projective} refers 
to the sign indeterminancy in Theorem~\ref{th-cobs-i}. 
Objects of $TC$ are even length sequences of pluses and minuses. 
One-morphisms are tangles with prescribed orientations at the top and 
bottom endpoints; 2-morphisms are isotopy classes of tangle cobordisms.  
The 2-functor $F$ assigns $n$ to a signed sequence of length $2n$, 
complex of graded bimodules $F(T)$ to a tangle $T$, and homomorphism 
$\pm F(S)$ to a cobordism $S$. 

Specializing from tangles to links, we get a (projective) functor 
from the category of link cobordisms to the category of bigraded abelian groups
\cite{Jac}, \cite{Kh2}. 
To a link $L$ it assigns Khovanov homology $H(L)$, to a link cobordism $S$ 
between $L_0$ and $L_1$ it assigns a homomorphism $\pm H(S): H(L_0) \lra H(L_1)$ 
of bidegree $(0, - \chi (S))$.  The sign indeterminacy in Theorem~\ref{th-cobs-i} has 
been eliminated by D.~Clark, S.~Morrison and K.~Walker~\cite{CMW} 
and C.~Caprau~\cite{Ca} at the cost 
of certain decorations of tangles and cobordisms. 

For more information on rings $H^n$ and 
related topological invariants see~\cite{CK}, \cite{Str}, \cite{SW} and references therein.

\subsection{Equivariant versions and applications} 

The ring $H^n$ and complexes $C(D)$ 
can be defined for any commutative Frobenius $R$-algebra $A$. Requiring that 
$C(D)$ be invariant under the first Reidemeister move implies the condition that 
$A$ is a rank two free $R$-module~\cite{Kh-e}. It turns out that there are 
many examples of rank two Frobenius pairs $(R,A)$ giving rise to link homology, 
tangle and tangle cobordism invariants. They were originally described by 
D.~Bar-Natan~\cite{BN2} in a more categorical language, avoiding the use of 
$H^n$. 

One of these rank two Frobenius pairs is given by $R=\Q[t]$ and $A=\Q[X],$ with the 
condition that $X^2=t$, making $A$ an $R$-algebra. The trace map is 
$\epsilon(X)=1, \epsilon(1)=0$. It is natural to think of $R$ as the 
$SU(2)$-equivariant cohomology of a point and $A$ as the $SU(2)$-equivariant 
cohomology of the 2-sphere, with $SU(2)$ acting via the surjective homomorphism 
onto $SO(3)$: 
\begin{eqnarray*} 
R & = & \mathrm{H}^{\ast}_{SU(2)}(\mathrm{pt}, \Q) = \mathrm{H}^{\ast}(BSU(2), \Q)
 = \mathrm{H}^{\ast}(\mathbb{HP}^{\infty}, \Q)\cong \Q[t], \\
A & = & \mathrm{H}^{\ast}_{SU(2)}(\mathbb{S}^2, \Q)\cong 
   \mathrm{H}^{\ast}_{SO(2)}(\mathrm{pt}, \Q)\cong \Q[X].   
\end{eqnarray*} 
Construction of Lecture 3, done for this $(R,A)$, produces a functorial  link homology 
theory, which we denote $H_t$. It extends to tangles and tangle cobordisms via 
the framework described in Lecture 4 and the first two sections of Lecture 5.
A cobordism $S$ between links $L_0, L_1$ induces 
a homomorphism $H_t(S): H_t(L_0)\lra H_t(L_1)$, well-defined up to overall minus sign. 
Groups $H_t(L)$ are bigraded finitely-generated 
$\Q[t]$-modules, with the multiplication by $t$ shifting the bigrading by $(0,4)$. 
Let $Tor(L)\subset H_t(L)$ be the torsion submodule; it consists of elements of $H_t(L)$ 
annihilated by some power of $t$. Homomorphisms $H_t(S)$ take $Tor(L_1)$ into 
$Tor(L_2)$, thus $Tor$ is a functorial subtheory of $H_t$. Let 
$H'(L) = H_t(L)/Tor(L)$. The quotient theory  $H'$ is functorial with respect to link 
cobordisms, and each $H'(L)$ is a free bigraded $\Q[t]$-module. It follows 
from~\cite{Lee2} that $H'(L)$ has rank $2^m$, where $m$ is the number 
of components of $L$. Moreover, when $L$ is a knot, $H'(L)$ lives in 
cohomological degree $0$ and 
$$ H'(L) \cong \Q[X]\{-s(L)-1\} $$ 
for some even integer $s(L)$ called the Rasmussen invariant of $L$. 
The Rasmussen invariant tells us where the internal 
grading of $H'(L)$ starts: in $q$-degree $-s(L)-1$. 
By writing $\Q[X]=\Q[t]\cdot 1 \oplus \Q[t]\cdot X$, we see two 
copies of $\Q[t]$ with the relative grading shift by $2$.  

Rasmussen~\cite{Ras1} showed that, given a connected cobordism $S$ between 
knots $L_0, L_1$, the induced homomorphism of $\Q[X]$-modules 
$H'(S): H'(L_0) \to H'(L_1)$ is nontrivial. Moreover, this homomorphism 
has degree $-\chi(S)$. Since 
$$ H'(L_0) \cong \Q[X]\{-s(L_0)-1\}, \ \ \ 
     H'(L_1) \cong \Q[X]\{-s(L_1)-1\}, $$ 
nontriviality of the homomorphism implies that the absolute value of 
the difference $s(L_0) - s(L_1)$ is bounded by twice the genus of $S$, 
 $$ | s(L_0) - s(L_1) | \le 2g(S)= -\chi(S). $$
In particular, $L\longmapsto s(L)$ descends to a homomorphism from 
the knot concordance group to $2\Z$. Specializing to cobordisms from the
trivial knot to $L$, one gets a lower bound on the slice genus of $L$: 
$$  | s(L) |  \le 2 g_4(L) .$$ 
The slice genus $g_4(L)$ of a knot $L$ is the minimum genus of a smooth 
oriented surface in  the four-ball $D^4$ that bounds $L\subset S^3= \partial D^4$.
It is also the minimum genus of a cobordism between the trivial knot and $L$.

Generally, both the slice genus $g_4(L)$ and the Rasmussen invariant 
$s(L)$ are very difficult to compute. 
For positive knots, however, the computation of $s(L)$ is straightforward. 
If  $L$ is a positive knot with a positive diagram $D$, the complex $C_t(D)$ 
starts in cohomological degree $0$, and  $H_t^0(D)$ is the kernel of 
the differential $C^0_t(D) \lra C^1_t(D)$. 
This allows us to determine 
 $H_t^0(D)$, its quotient $H'(D)$, and find the Rasmussen invariant $s(L)$. 
It is equal to $n+1-c$, where $n$ is the number of 
crossings of $D$, and $c$ is the number of Seifert circles. At the same time, 
Seifert's algorithm gives a Seifert surface for $L$ of genus $\frac{n+1-c}{2}$, 
so the ordinary genus $g(L) \le \frac{n+1-c}{2}$. 
The chain of inequalities 
$$ g(L) \ge g_4(L) \ge \frac{|s(L)|}{2}=\frac{n+1-c}{2} \ge  g(L)$$
implies that all of them are equalities, and 
the slice genus of $L$ is $\frac{n+1-c}{2}$. 
As a special case, this argument proves the Milnor conjecture, also known as 
Kronheimer-Mrowka theorem, that the slice genus 
of the torus knot $T_{p,q}$ is $\frac{(p-1)(q-1)}{2}.$ The first proof of the Milnor 
conjecture  was given by P.~Kronheimer and T.~Mrowka~\cite{KM} via 
Donaldson theory. The above much more recent proof, due to Rasmussen, is 
algebraic. The sketch, presented here, uses the graded theory $H_t$ instead 
of its filtered version, utilized by Rasmussen~\cite{Ras1}. 

Two 1's in the zero column of the homology table for $T_{5,6}$,  
depicted at the end of Lecture 3, are all that's left of $H'(T_{5,6})$ in 
the homology groups $H(T_{5,6})$, where $t=0$ and $\Q[X]$ becomes 
$A$. The Rasmussen invariant of this torus knot equals $20$. 

\vspace{0.1in} 

Extension of the link homology to tangles, in addition to giving an easy proof of 
functoriality, also helps with computing link homology, as demonstrated by 
D.~Bar-Natan and J.~Green, who produced a fast program for computing 
Khovanov homology~\cite{BN3}. One puts a link $L$ in 
"thin" position, namely a position minimizing  the number of intersection points of 
horizontal planes with $L$, as illustrated below. 
\psfrag{dt}{$\vdots$}
\psfrag{T1}{$T_1$}
\psfrag{T2}{$T_2$}
\psfrag{T3}{$T_3$}
\psfrag{Tk}{$T_k$}
$$\raisebox{-2cm}{\psfig{figure=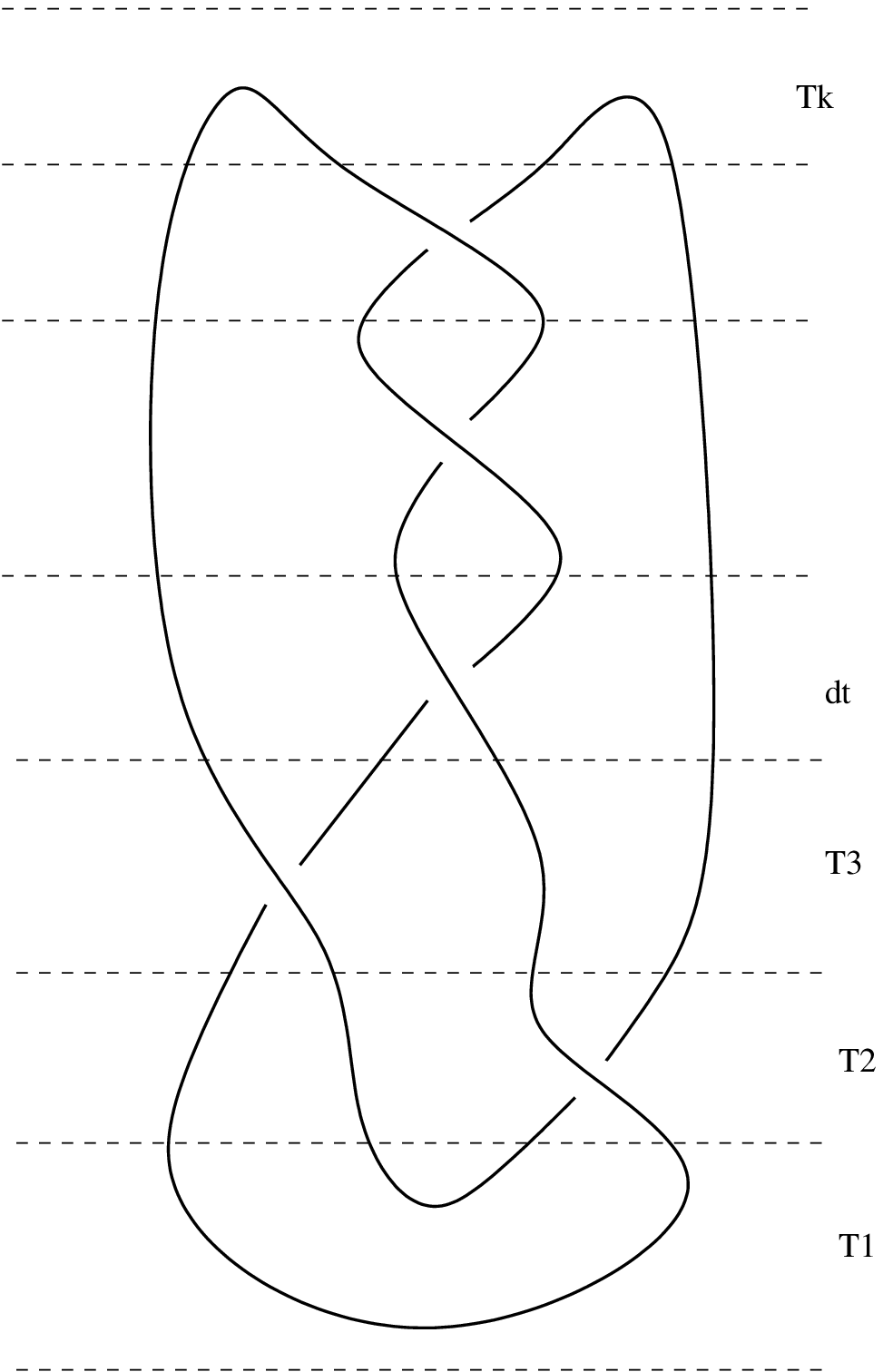,height=4cm}}$$
Write $L$ as the product of tangles with at most one crossing each, 
$L=T_k \dots T_2 T_1$. 
To compute $H(L)=F(L)$, one starts with $F(T_1)$, then 
computes $F(T_2T_1), $ $F(T_3 T_2 T_1)$, etc. Each $F(T_i \dots T_2 T_1)$ 
is a complex of projective $H^m$-modules, where $2m$ is the number of 
top endpoints of $T_i$. At each step, one simplifies $F(T_i\dots T_2 T_1)$ as much as possible 
by removing null-homotopic components, isomorphic to 
$0 \to P_a \stackrel{1}{\to} P_a \to 0$, where $P_a$ are projective $H^m$-modules
described earlier and labelled by crossingless matchings $a$. 
After that, the reduced complex is tensored with $F(T_{i+1})$,  
and the simplification procedure is repeated. This method gives the most 
efficient algorithm at present for computing link homology.


\section{Categorifications of the HOMFLY-PT polynomial}

\subsection{The HOMFLY-PT polynomial and its generalizations} 
The HOMFLY-PT polynomial~\cite{HOMFLY}, \cite{PT} 
is a generalization of the Jones polynomial which is determined by 
the skein relation 
\begeq
a  P \left( \raisebox{-0.3cm}{\psfig{figure=lec3.1.eps,height=0.8cm}}\right) 
- a^{-1} P \left( \raisebox{-0.3cm}{\psfig{figure=lec3.2.eps,height=0.8cm}}\right)
&=& (q-q^{-1})P\left( \raisebox{-0.3cm}{\psfig{figure=lec3.3.eps,height=0.8cm}}\right)
\eneq	
and its value on the unknot 
\begeq 
P \left(  \raisebox{-0.3cm}{\psfig{figure=lec3.4.eps,height=0.8cm}}\right) &=&  
\frac{a-a^{-1}}{q-q^{-1}}.
\eneq
For $P(L)$ to really be a (Laurent) polynomial, one should change variables in the 
above formulas by introducing $b=q-q^{-1}$, for then $P(L)\in \Z[a^{\pm 1}, b^{\pm 1}]$. 
Variables $a,q$ are natural from the representation-theoretical viewpoint, though, 
since the one-variable specialization $P_n(L):= P_{a=q^n}(L)$ of $P(L)$ for $n>0$ 
can be extended~\cite{RT} 
to an invariant of tangles via representation theory of $U_q(sl(n))$, a 
Hopf algebra deformation of the universal enveloping algebra $U(sl(n))$. 
For the first few values of $n$, the polynomial $P_n(L)$ is as follows: 
\begin{itemize} 
\item $P_0(L)$ is the Alexander polynomial of $L$. 
\item $P_1(L)=1$ for all $L$ is a trivial invariant. 
\item $P_2(L)=J(L)$ is the Jones polynomial of $L$. 
\end{itemize} 
We already discussed a categorification of $P_2(L)$. A categorification 
of $P_0(L)$ (the Alexander polynomial) has been constructed by P.~Ozsv\'ath, 
Z.~Szab\'o~\cite{OS1} and, independently, J.~Rasmussen~\cite{RasT}. 
It is a bigraded homology theory, 
known as the knot Floer homology, which comes in several versions and has 
found a multitude of applications in low-dimensional topology, see~\cite{OS-M} and 
references therein. The invariant $P_1(L)$ is trivial, so there's nothing to categorify. 
Polynomial $P_n(L)$ was categorified in~\cite{Kh3} for $n=3$ (see~\cite{MV}, 
\cite{MN} for extension to tangles and for equivariant versions) and in~\cite{KR1} 
for all $n>1$. Both constructions employ a generalization of the Kauffman's bracket 
decomposition of the crossing, which for $n>2$ takes the form 
$$\raisebox{-0.3cm}{\psfig{figure=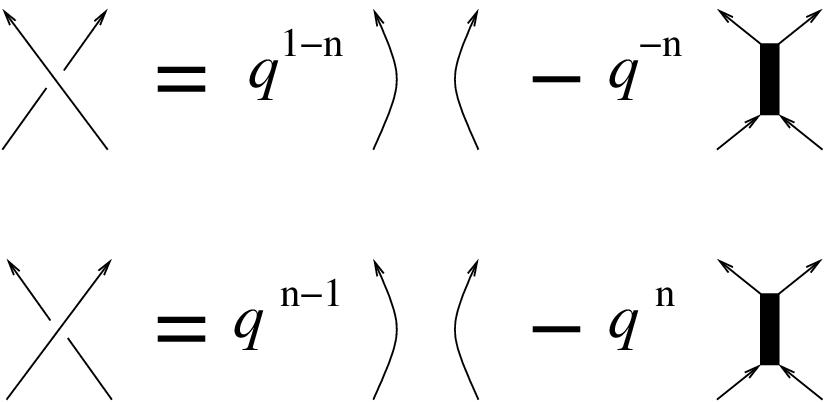,height=3.5cm}}$$
Each crossing of a diagram is resolved in two possible ways, and a complete 
resolution of a diagram produces a planar graph of a particular type. The 
invariant $P_n(L)$ extends to these graphs and has a positive evaluation on 
each of them: $P_n(\Gamma)\in \mathbb{N}[q,q^{-1}]$ for a planar graph $\Gamma$. 
To categorify $P_n(L)$ one first categorifies $P_n(\Gamma)$, which become 
graded dimensions of graded $\Q$-vector spaces $H_n(\Gamma)$ (that's 
the trickiest part of the construction). These vector 
spaces are then put into the vertices of an $m$-dimensional cube, where $m$ 
is the number of crossings of the diagram, and for each edge of the cube one 
constructs a linear map between the spaces. Homology of $D$ is defined as 
the homology of the total complex of the cube and checked to be invariant under 
the Reidemeister moves. The result~\cite{KR1} is a family of bigraded link homology 
theories 
$$ H_n(L) = \oplusop{i,j\in \Z} H^{i,j}_n(L), \ \ n>1, $$ 
with each group $H_n^{i,j}(L) $ a finite-dimensional $\Q$-vector space, 
and the Euler characteristic 
$$ P_n(L) = \sum_{i,j\in \Z} (-1)^i q^j \dim H^{i,j}_n(L)$$ 
(when $n=2$, we recover Khovanov homology, tensored with $\Q$.) 
These invariants extend to tangles and tangle cobordisms in a way conceptually 
similar to the one described in Lectures 4 and 5 for $H(L)$. 

There are now several strikingly different constructions of graded~\cite{Man} and 
bigraded~\cite{Sus}, \cite{CaKa2}, \cite{MSV}, \cite{MS} homology theories of links 
which exhibit behavior similar to $H_n$ and might all be 
isomorphic to $H_n$ or mild modifications of the latter (see also~\cite{SS}, \cite{CaKa1}, 
\cite{Str} for prior constructions in the $n=2$ case).

In the rest of the lecture we discuss a triply-graded homology theory~\cite{KR2}, 
\cite{Kh4}  which categorifies the 2-variable HOMFLY-PT polynomial $P(L)$. 


\subsection{Hochschild homology}
Let $R$ be a ring, and $M$ (resp. $N$) be a right (resp. left) $R$-module. 
The tensor product $M \otimes_R N$ is an abelian group. 
In homological algebra various functors need to be redefined to make them 
exact in a suitable category. For general $M$, tensor product with $M$ is 
only a right exact functor on the category of left $R$-modules. To convert it 
into an exact functor (on a bigger category, say the derived category of the 
original category),  we consider a  
projective resolution of $M$, that is, a chain complex of projective $R$-modules 
$(P_i, \varphi_i )$ so that the vertical chain map in the following diagram 
is a quasi-isomorphism (i.e. isomorphism on homology). 
$$
\begin{array}{ccccccccc}
\cdots & \to & P_2 & \stackrel{\varphi_2}{\to} & P_1 &  \stackrel{\varphi_1}{\to}  & 
P_0 & \to&  0 \\
&&&&\downarrow && \downarrow && \\
&&&& 0 & \to & M & \to & 0
\end{array}
$$
This simply means that the chain complex is exact everywhere except the last term, 
where $M\cong P_0/{{\rm Im} \varphi_1}$. 

Consider the chain complex $M \stackrel{L}{\otimes} N:= (P_i, \varphi_i ) 
\otimes_R N=(P_i \otimes_R N, \varphi_i  \otimes {\rm id})$. We call the 
$i$-th homology of this complex {\it the $i$-th derived tensor product} of $M$ and $N$. 
 It is known that derived tensor products do not depend on the choice of 
projective resolution, and that if we projectivize $N$ instead of $M$, we get 
the same answer as well. 

\begin{exercise} Determine the derived tensor product of $\Z$-modules 
$\Z_n$ and $\Z_m$. 
\end{exercise}
We represent a right $R$-module $M$ 
and a left $R$-module $N$ graphically as 
\psfrag{M}{$M$}
\psfrag{N}{$N$}
$$
\psfig{figure=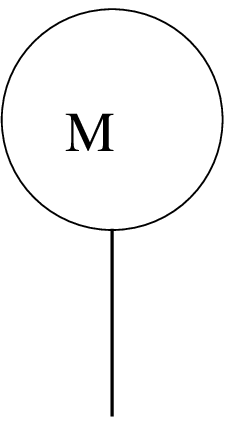,height=1.5cm}
\hspace{1cm}
 \psfig{figure=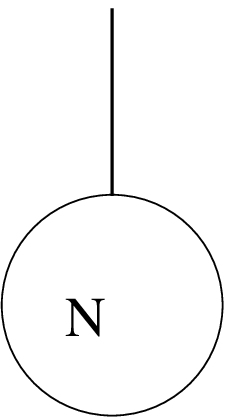,height=1.5cm}
 $$
 $R$-action is depicted by wires, and "left" is seen as "up", "right" is "down". 
Turn your head $90$ degrees clockwise or, alternatively, turn 
the paper $90$ degrees counterclockwise to see the match.  
We represent $M \stackrel{L}{\otimes} N$ by the following picture:
 \psfrag{M}{$M$}
\psfrag{N}{$N$}
$$ \psfig{figure=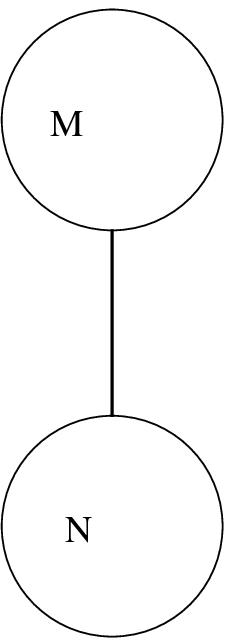,height=2cm} $$
If $M$ and $N$ are $R$-bimodules, we depict them and their 
derived tensor product $M \stackrel{L}{\otimes} N$ as follows:
 \psfrag{M}{$M$}
\psfrag{N}{$N$}
$$
\psfig{figure=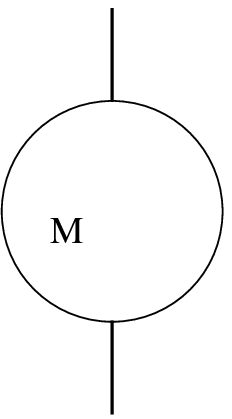,height=1.5cm}
\hspace{1cm}
 \psfig{figure=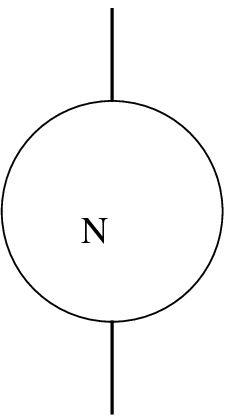,height=1.5cm}
 \hspace{1cm}
 \psfig{figure=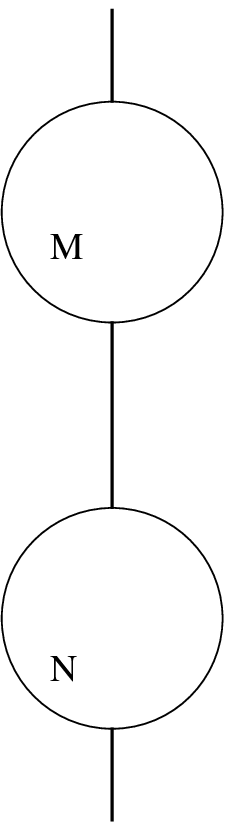,height=2cm}
 $$
(compare with Lecture 4). The top and the bottom wires of the diagram 
for a bimodule $M$ indicate left and right actions of $R$.  A single vertical 
undecorated wire denotes $R$ viewed as an $R$-bimodule. 

The space of $R$-coinvariants of an $R$-bimodule $M$ is 
$$ M_R := M/[R,M],$$ 
the quotient of $M$ by the abelian subgroup generated by expressions 
$rm-mr$ over all $r\in R, m\in M$. The functor  $M \mapsto M_R$ 
is right exact. 

\begin{remark} ``Quotient object'' functors (such as the $R$-coinvariants functor) 
 between abelian categories are often right exact. 
 The ``subobject'' functors tend to be left exact. 
An example of a ``subobject'' functor is $M \longmapsto M^R$, which 
to a bimodule $M$ assigns its $R$-invariants
 $$M^R := \{ m\in M | rm=mr \ \text{for all} \ r\in R\}.$$
\end{remark} 
 
Graphically, passing from $M$ to its $R$-coinvariants should correspond 
to joining the two wires of $M$. If we imagine elements of $R$ moving 
along the wires, the equations $rm=mr$ that hold in $M_R$ for all $r$ and $m$ 
mean that $r$ can jump from the top to the bottom wire and back without changing 
the value of the diagram. The easiest way to achieve this geometrically is by 
closing off the two ends of the diagram: 
\psfrag{M}{$M$}
\begin{equation} \label{eq-sharp} 
 \ \ \ \raisebox{-1cm}{\psfig{figure=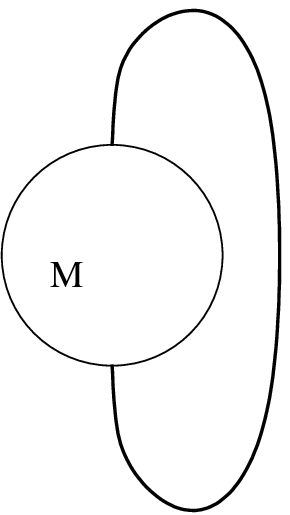,height=2cm}}
=M/[R,M]=M_R.
\end{equation}
This is only an approximation, though, since $M\longmapsto M_R$ is 
not left exact and we need to form its derived functor, known as 
Hochschild homology. Notice that $R$-bimodules are the same 
as left (or right) modules over the ring $R^e := R\otimes R^{op}$. If we are viewing 
$R$ as a $k$-algebra, for some commutative ring $k$, often a field, the 
tensor product in the definition of $R^e$ should be taken over $k$. 
The group $M_R$ equals the tensor product $M \otimes_{R^e} R $ 
of a right $R^e$-module $M$ and a left $R^e$-module $R$ (right and left here 
can be transposed). We define the $i$-th Hochschild homology 
of $M$ as the $i$-th derived functor of the tensor product: 
\begin{eqnarray*} 
\mathrm{HH}_i(R, M ) & := & \mathrm{H}_i (M \stackrel{L}{\otimes}_{R^e} R), \\
  \mathrm{HH}_{\ast}(R, M ) & := & \oplusop{i\ge 0 } \  \mathrm{HH}_i(R, M ). 
\end{eqnarray*}  
Going back to diagrammatics, we should interpret the closure of a bimodule diagram 
as taking the entire Hochschild homology of $M$ rather than just $M_R$, its 
degree $0$ part: 
 $$\raisebox{-1.2cm}{\psfig{figure=lec6.8.eps,height=2.5cm}}={\rm HH}_*(R, M).$$
 The Hochschild homology exhibits "tracial" behaviour, since there are (functorial in $M$ 
and $N$) isomorphisms
 $${\rm HH}_{\ast} (R,M \stackrel{L}{\otimes} N)\cong {\rm HH}_{\ast} 
(R,N \stackrel{L}{\otimes} M).$$
These isomorphisms acquire topological interpretation 
 $$\raisebox{-1.2cm}{\psfig{figure=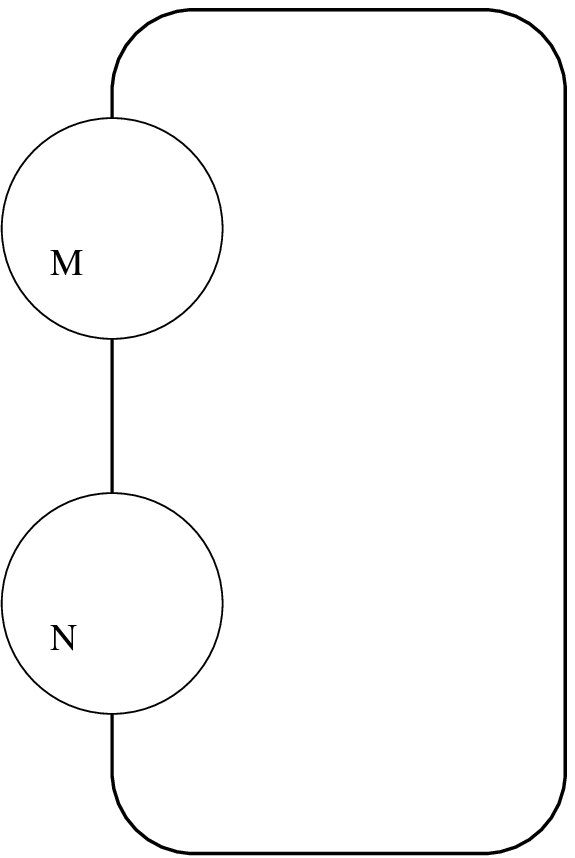,height=2.5cm}} \ \Leftrightarrow \ 
 \raisebox{-1.2cm}{\psfig{figure=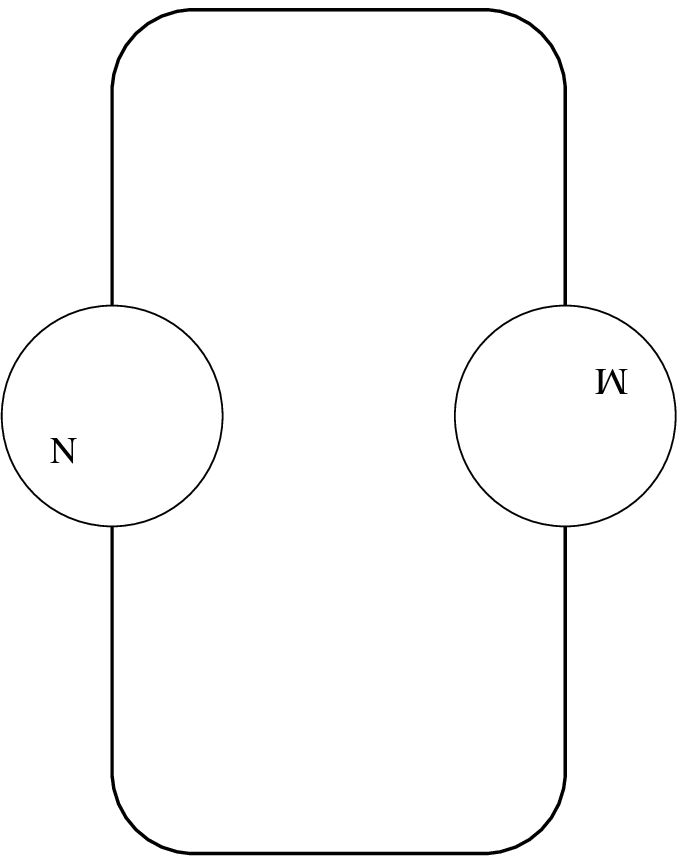,height=2.5cm}}  \ \Leftrightarrow \
 \raisebox{-1.2cm}{\psfig{figure=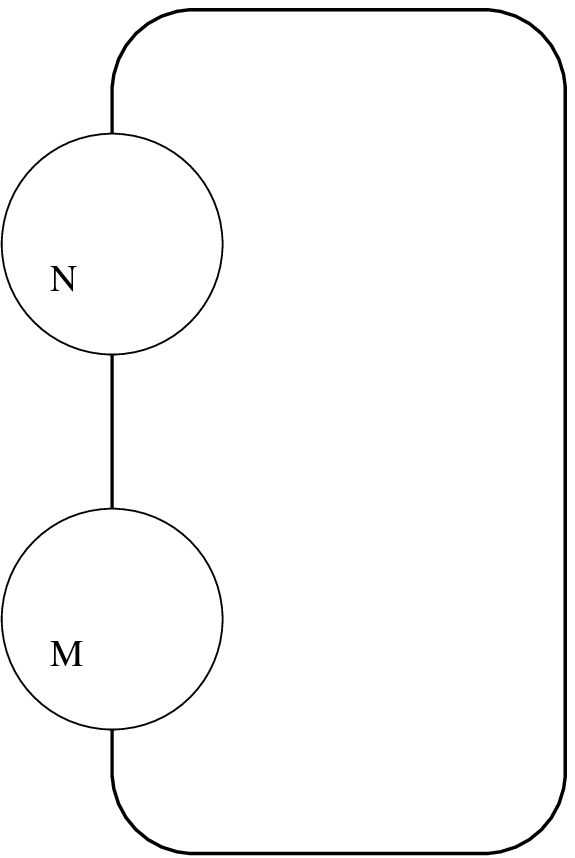,height=2.5cm}},
 $$
that is, bimodule boxes can be dragged along the wires. 
Hochschild homology of a bimodule can be viewed as a categorification of the trace of a 
linear operator. 

 To compute the Hochshild homology of $M$, it suffices to construct a projective resolution of
 $R$ as an  $R$-bimodule, tensor with $M$ over $R^e$ and take homology.  
 \begin{example} Let $R=\Q[x]$, viewed as a $\Q$-algebra. 
Note that $R^{\rm op}\cong R$. The following is a free  
resolution of $R$ as $R\otimes R$-module called \emph{Koszul resolution} 
$$
0\to \Q[x] \otimes \Q[x] \stackrel{\varphi}{\to} \Q[x] \otimes \Q[x] \to 0, 
$$
where $\varphi$ is the $R\otimes R$-module map determined by $\varphi(1\otimes1)=
x \otimes 1 - 1\otimes x.$ Tensoring with $M$, we get
$$(M \otimes_{R^e} (0 \to R^e \stackrel{\varphi}{\to}  R^e \to 0))  \cong 
(0 \to M  \stackrel{\psi}{\to}  M \to 0), $$
where $\psi(m)=xm-mx$ (notice that $M\otimes_{R^e} R^e = M$). Therefore, 
$$ {\rm HH}_0(R,M)\cong M_R, \ \ \   {\rm HH}_1(R,M)\cong M^R, $$ 
and all higher Hochschild homology groups vanish.

 \end{example}
 \begin{example} \label{ex-poly} 
$R=\Q[x_1, \dots, x_n].$ We have a resolution of $R$ by 
free $R^e$-modules 
 \begin{eqnarray*}
 &&\bigotimes_{i=1}^n (0 \to \Q[x_i] \otimes \Q[x_i] \stackrel{\varphi_i}{\to} 
\Q[x_i] \otimes \Q[x_i] \to 0) \\
 &=& (0 \to R \otimes R \to \cdots (R \otimes R)^{\oplus  \binom{n}{k} }\to \cdots 
\to R\otimes R \to 0),
  \end{eqnarray*}
  where $\varphi_i(1\otimes1)=x_i \otimes 1 -1 \otimes x_i.$
  One may consider the $k$-th term as
  $$(R \otimes R)^{\oplus  \binom{n}{k} } \cong R\otimes R \otimes \wedge^k V,   $$
  where $V:=\spa_\Q \{y_1, ..., y_n\},$ and the differential is given by 
\begin{eqnarray*}   
& & d(z_1\otimes z_2\otimes y_{r_1} \wedge \cdots \wedge y_{r_k} )= \\
& &  \sum_{j=1}^k  (-1)^j (x_{r_j}z_1 \otimes z_2 -z_1 \otimes z_2 x_{r_j}) 
  y_{r_1} \wedge \cdots  \widehat{y_{r_j}} \cdots \wedge y_{r_k}. 
\end{eqnarray*}
 Tensoring with $M$, we get the complex  
  $ 0 \to M \to \cdots \to M^{\oplus \binom{n}{k}} \to \cdots \to M \to 0$
which computes Hochschild homology of $M$. For the boundary terms, we get 
$${\rm HH}_0(R,M)=M_R, \ \  {\rm HH}_n(R,M)=M^R. $$ 
 \end{example}


 \subsection{A categorification of the HOMFLY-PT polynomial}
 We use $R=\Q[x_1, \dots, x_n]$ as in Example~\ref{ex-poly}. 
For the transposition $s_i=(i,i+1) $ in the symmetric group $S_n$ let 
 $$R_i:=R^{s_i}=\Q[x_1, ..., x_{i-1}, x_i+x_{i+1}, x_i x_{i+1}, x_{i+2}, ... , x_n]\subset 
R,$$
be the space of $s_i$-invariants under the permutation action of $S_n$ on $R$. 
As an $R_i$-module, $R$ is free of rank $2$ and can be written as 
$R=R_i \cdot 1 \oplus R_i \cdot x_i$. We set the degree of $x_i$ to $2$; this 
makes $R$ and $R_i$ into graded rings. 
Then $B_i:=R \otimes_{R_i} R\{-1\}$ is a graded $R$-bimodule, and we have
$$B_i \otimes_R B_i\cong R \otimes_{R_i} ( R_i \cdot 1 \oplus R_i \cdot x_i) 
\otimes_{R_i} R\{-1\}\cong B_i\{1\} \oplus B_i \{-1\}.$$
Recall that bimodules $U_i$ in Lecture 1 satisfy the same relation. 
In Lecture 1 we formed chain complexes $(0 \to U_i \to A_n \to 0)$ and 
$(0 \to A_n \to U_i \to 0)$ corresponding to the braid $\s_i$; these complexes 
gave rise to a braid group action in the homotopy category of complexes. An 
analogous theory exists for bimodules $B_i$. Form bimodule complexes  
\begeq
C_i &:=&(0\to B_i\{1\}=R \otimes_{R_i} R\stackrel{m}{\lra} R \to 0), \\
C'_i &:=&(0\to R  \stackrel{\psi}{\to} B_i \{-1\} \to 0),
\eneq
where $\psi(1):=(x_i -x_{i+1})\otimes 1 + 1 \otimes (x_i -x_{i+1})$ and 
in both complexes $R$ sits in cohomological degree $0$. 
\begin{theorem} (R.Rouquier \cite{Rou1})
In the category of complexes of graded $R$-bimodules modulo homotopic to 
zero morphisms there are the following isomorphisms:
\begeq
C_i \otimes C'_i &\cong& R \  \cong \ C'_i \otimes C_i, \\
C_i \otimes C_{i+1}  \otimes C_i &\cong& C_{i+1} \otimes C_i \otimes C_{i+1}, \\
C_i \otimes C_j &\cong& C_j \otimes C_i, \ {\rm if} \ |i-j| >0.
\eneq
\end{theorem}
The theorem say that there is a weak braid group action on the homotopy 
category of complexes of graded $R$-modules. Rouquier also showed that this weak  
action lifts to a genuine action. Moreover, just like in the example of Lectures 1 and 2, 
this braid group action extends to an action of the category of braid cobordisms~\cite{KT}. 

This is reminiscent of our previous categorifications, via rings $A_n$ and $H^n$:  
$$
 \begin{array}{cccc} 
 A_n: & \mbox{Braid cobordisms} & \stackrel{\mbox{categorification}}{\Leftarrow} & 
\mbox{Burau representations} \\
 H^n: & \mbox{Tangle cobordisms} & \stackrel{\mbox{categorification}}{\Leftarrow} & 
\mbox{Jones polynomials}\end{array}
 $$

It turns out that the braid group action via complexes of $R$-bimodules leads to a 
categorification of the HOMFLY-PT polynomial of braid closures. 

Starting with an arbitrary braid word $\sigma$, which is a product of $\sigma_i$'s and their 
inverses, form the corresponding complex of graded $R$-bimodules, denoted $C(\sigma)$, the 
tensor product of  $C_i$'s and $C_i'$'s: 
  $$C(\sigma): \hspace{0.2in}\cdots \dto C^j (\sigma) \dto C^{j+1}(\sigma) \dto \cdots.$$
  Now, take the Hochschild homology of each term. The differential map $d$ induces a 
mapping of Hochschild homology groups, and we obtain the chain complex
  $$\cdots \to {\rm HH}(R, C^j(\s)) \stackrel{{\rm HH}(d)}{\to}  
{\rm HH}(R, C^{j+1}(\s)) \to \cdots. $$
  Each ${\rm HH}(R, C^j(\s)) $ is a $\Q$-vector space with two gradings: the Hochschild 
grading and the internal grading (the grading of $R$ by deg$x_i=2$). Therefore, 
taking the homology, we get a triply-graded vector space  
  $${\rm H}({\rm HH}(R, C(\s)), {\rm HH}(d))={\rm HHH}(\s).$$
This triply-graded vector space needs an overall shift, as desribed by 
Hao Wu~\cite{Wu}. With it in place, we have 
\begin{theorem} Triply-graded homology groups  ${\rm HHH}(\s)$ is an invariant 
of the link ${\hat \s}$, the closure of braid $\sigma$. The Euler characteristic of 
${\rm HHH}(\s)$ is the HOMFLY-PT polynomial 
of the link ${\hat \s}$. 
\end{theorem}
This construction simplifies a categorification of the HOMFLY-PT polynomial in~\cite{KR2}. 
There are some problems with this homology theory: ${\rm HHH}$ is not functorial 
under link cobordisms,for instance due to the theory being infinite-dimensional on 
non-empty links. It might be possible 
to make it finite-dimensional by setting several $x_i$'s in $R$ to $0$, one for each 
component of the link. It is not clear, though, why the finite-dimensional version should be 
functorial under tangle cobordisms. One would also like to define ${\rm HHH}$ in a 
more natural way and extend it to tangles. 
Overall, it is an open problem to develop the homology theory  ${\rm HHH}$ 
and make it as aethetically pleasing as the one described in lectures 3-5. 

There are ideas on how to generalize this theory to the so-called colored HOMFLY 
polynomials and relate it to topological strings~\cite{GSV}, \cite{GIKV}.  
The homology was computed for a number of knots by Rasmussen~\cite{Ras2}, and earlier, 
via a computer program, by Ben Webster. 

The ring $R$ has a geometric interpretation as the $GL(n)$-equivariant 
cohomology of the variety of full flags in $\C^n$. This interpretation 
was extended to the Hochschild homology of indecomposable summands of $C(\sigma)$ 
in~\cite{WW}. 

Similarities between the Hochschild homology and link homology were originally
observed in~\cite{Pr}. 

\vspace{0.1in} 

To summarize, we've seen three categorifications in these lectures.
\begin{itemize}  
\item Rings $A_n$ lead to a categorification of 
 the reduced Burau representation of the braid group. 
Braids act by complexes of $A_n$-bimodules. The theory can be extended to give invariants 
of braid cobordisms via homomorphisms of complexes of bimodules. 
\item 
Bimodules over rings $H^n$ give a categorification of the Temperley-Lieb algebra. 
Complexes of bimodules produce an invariant of tangles; homomorphisms of 
complexes--invariants of tangle cobordisms. The construction specializes to a 
bigraded homology theory of links categorifying the Jones polynomial. 
\item 
Suitable bimodules over polynomial rings $R$ give rise to a braid group action. 
The action extends to braid cobordisms. Taking Hochschild homology produces 
a triply-graded link homology theory categorifying the HOMFLY-PT polynomial. 
\end{itemize} 

\begin{exercise} Choose an integral structure that appears in combinatorics, or 
algebra, or topology, etc. and categorify it. 
\end{exercise}


\vspace{0.2in}  

\noindent
M.A.:  { \sl \small Department of Mathematics, University of California, Riverside, CA 92521}
\noindent
  {\tt \small email: marta@math.ucr.edu}

\noindent 
M.K.: {\sl \small School of Mathematics, Institute for Advanced Study, Princeton, NJ 08540, and}
{ \sl \small Department of Mathematics, Columbia University, New York, NY 10027}

\noindent
  {\tt \small email: khovanov@math.columbia.edu}

\end{document}